\documentclass{article}
\usepackage{amsfonts}
\usepackage{amssymb,latexsym,amsmath,amscd,epsfig}
\usepackage{multirow}

\newtheorem{eg}{Examples}[section]
\renewcommand{\theeg}{}
\newtheorem{definition}{Definition}[section]

\newtheorem{lem}{Lemma}[section]
\newtheorem{prop}{Proposition}[section] 
\newtheorem{thm}{Theorem}[section]
\newtheorem{remark}{Remark}[section]
\newcommand{\pa}{\partial}

\newcommand{\nab}{\nabla}
\newcommand{\ka}{K\"{a}hler }

\newcommand{\BR}{\mathbb{R}}
\newcommand{\BC}{\mathbb{C}}
\newcommand{\BP}{\mathbb{P}}
\newcommand{\BZ}{\mathbb{Z}}

\newcommand{\CY}{Calabi--Yau }

\newcommand{\we}{\wedge}
\newcommand{\oa}{\omega}
\newcommand{\oal}{\oa^{(1,1)}}
\newcommand{\Oa}{\Omega}
\newcommand{\ga}{\gamma}
\newcommand{\Ga}{\Gamma}
\newcommand{\si}{\sigma}
\newcommand{\Si}{\Sigma}
\newcommand{\da}{\delta}

\newcommand{\ta}{\theta}

\newcommand{\za}{\zeta}
\newcommand{\ep}{\epsilon}
\newcommand{\al}{\alpha}
\newcommand{\ba}{\beta}
\newcommand{\kp}{\kappa}
\newcommand{\lda}{\lambda}
\newcommand{\Lda}{\Lambda}
\newcommand{\vi}{\varphi}
\newcommand{\Ups}{\Upsilon}

\newcommand{\pf}{\emph{Proof. }}
\newcommand{\SM}{S^1 \times M}
\newcommand{\SMt}{S^1 \times M_t}

\numberwithin{equation}{subsection}

\begin{document}

\title{\large \bf SIMULTANEOUS DESINGULARIZATIONS OF CALABI--YAU AND SPECIAL LAGRANGIAN 3-FOLDS WITH CONICAL SINGULARITIES}
\author{ \mdseries{YAT-MING CHAN} \\[0.5cm] (\emph{Mathematical Institute, 24-29 St Giles, Oxford OX1 3LB}) }
\date{} 
\maketitle
\bigskip
\noindent
\small \textsc{Abstract.} \ This paper is a follow-up to an earlier paper \cite{Chan1} on desingularizations of \CY 3-folds with a conical singularity. In \cite{Chan1} we study \CY 3-folds $M_0$ with a conical singularity at $x$ modelled on some \CY cone $V$, and construct a desingularization of $M_0$ by gluing in an Asymptotically Conical (AC) \CY 3-fold $Y$ to $M_0$ at $x$. In this paper, we shall investigate a similar desingularization problem on special Lagrangian (SL) 3-folds in the corresponding \CY 3-folds. More precisely, suppose $M_0$ is now a \CY 3-fold with finitely many conical singularities at $x_i$ modelled on \CY cones $V_i$ for $i=1, \dots,n$, and $N_0$ an SL 3-fold in $M_0$ with conical singularities at the same points $x_i$ modelled on SL cones $C_i$ in $V_i$. Let $Y_i$ be an AC \CY 3-fold modelled on the \CY cones $V_i$, and $L_i$ an AC SL 3-fold in $Y_i$ modelled on the SL cones $C_i$. We then simultaneously desingularize $M_0$ and $N_0$ by gluing in rescaled $Y_i$ and $L_i$ at each $x_i$. The construction is achieved by applying Joyce's analytic result \cite[Thm. 5.3]{Joycedesing3} on deforming Lagrangian submanifolds to nearby special Lagrangian submanifolds. As an application, we take $M_0$ to be the orbifold $T^6/\BZ_3$ and construct some singular SL 3-folds $N_0$ in $M_0$ and AC SL 3-folds $L_i$ in the corresponding $Y_i$, and glue them together to obtain examples of nonsingular SL 3-folds in the desingularized \CY 3-folds.      \\

\bigskip


\large

\centerline{1\ \textsc{\ Introduction}}
\bigskip 

\normalsize

\par Special Lagrangian (SL) submanifolds in \CY manifolds play an important role in the explanation of a phenomenon in physics called Mirror symmetry. They are examples of \emph{calibrated submanifolds}, appearing in Harvey and Lawson \cite{HarLaw}, which generalizes the concept of volume-minimizing property of complex submanifolds of \ka manifolds. Let $(M,J,\oa,\Oa)$ be a \CY manifold of complex dimension $m$. Then Re($\Oa$) is a calibrated form whose calibrated submanifolds are real $m$-dimensional special Lagrangian submanifolds (SL $m$-folds). Special Lagrangian submanifolds attracted much interest in connection with the \emph{SYZ conjecture} proposed by Strominger, Yau and Zaslow \cite{StromYauZas} in 1996, which explains Mirror symmetry between \CY 3-folds. Roughly speaking, the conjecture asserts that the mirror $\check{M}$ of a \CY 3-fold $M$ can be obtained by some suitable compactification of the dual of the SL $T^3$-fibration on $M$. Therefore to find a compactification and understand the relations with the Mirror symmetry one should understand the singularities of the moduli space of SL $m$-folds.  \\

\par Perhaps the simplest singularities to understand are isolated singularities modelled on SL cones. A lot of SL cones in $\BC^m$ have been constructed explicitly using various techniques, hence providing local models for conical singularities in SL $m$-folds in general \CY $m$-folds. For example, Joyce \cite{JoyceSLsym} used large symmetry groups to construct SL cones. Haskins \cite{Haskins} focused on dimension three and explored examples of SL cones in $\BC^3$. \\

\par Recently, Joyce has developed a comprehensive programme on the desingularization of SL $m$-folds with conical singularities in (almost) \CY manifolds and their deformation theory in his series of papers \cite{Joycedesing1}, \cite{Joycedesing2}, \cite{Joycedesing3}, \cite{Joycedesing4} and \cite{Joycedesing5}. The SL $m$-folds with conical singularities are desingularized by gluing in at the singular points some nonsingular SL $m$-folds in $\BC^m$ which are asymptotic to SL cones at infinity. \\ 


\par In \cite{Chan1} we obtained a desingularization of \CY 3-folds with a (or finitely many) conical singularity for the case $\lda < -3$ where $\lda$ denotes the \emph{rate} at which the Asymptotically Conical (AC) \CY 3-fold converges to the \CY cone. The result on this case can be applied to desingularizing \CY 3-orbifolds with isolated singularities. In the present paper we deal exclusively with the special Lagrangians inside the corresponding Calabi--Yau's, and simultaneously desingularize the \CY and SL 3-folds with conical singularities by gluing in AC \CY and AC SL 3-folds at each singular point. \\

\par Throughout this paper $M_0$ will denote a \CY 3-fold with finitely many conical singularities at $x_i$ with rate $\nu$ modelled on \CY cones $V_i$ for $i=1, \dots ,n$, and $N_0$ an SL 3-fold in $M_0$ with conical singularities at the same points $x_i$ with rate $\mu$ modelled on SL cones $C_i$ in $V_i$. We denote by $Y_i$ an AC \CY 3-fold with rate $\lda_i$ modelled on the \CY cones $V_i$, and $L_i$ an AC SL 3-fold with rate $\kp_i$ in $Y_i$ modelled on the SL cones $C_i$. When we glue in a rescaled $Y_i$ to $M_0$ at each $x_i$ as in \cite{Chan1}, we also glue in a rescaled $L_i$ to $N_0$ at each $x_i$. The idea of the special Lagrangian desingularization is to construct a family of nonsingular 3-folds $N_t$ in the family of nonsingular \emph{nearly \CY 3-folds} (see \cite[\S3]{Chan1}) $(M_t, \oa_t, \Oa_t)$ so that $N_t$ is \emph{Lagrangian} in $(M_t, \oa_t, \Oa_t)$. Then Theorem 5.2 of \cite{Chan1} gives genuine \CY 3-folds $(M_t, \tilde{J}_t, \tilde{\oa}_t, \tilde{\Oa}_t)$. Choose a suitable coordinate/diffeomorphism $\psi_t$ on $M_t$ so that $\oa_t = \psi_t^* (c_t \tilde{\oa}_t)$, the pullback $(\hat{J}_t, \hat{\oa}_t, \hat{\Oa}_t)$ of $(\tilde{J}_t, \tilde{\oa}_t, \tilde{\Oa}_t)$ under $\psi_t$ is also a genuine \CY structure on $M_t$, and $N_t$ is Lagrangian in the \CY 3-folds $(M_t, \hat{J}_t, \hat{\oa}_t, \hat{\Oa}_t)$ as well. The main analytic result we need in this paper is adapted from Joyce \cite[Thm. 5.3]{Joycedesing3}, in which he shows that when $t$ is sufficiently small, we can deform the Lagrangian $m$-fold $N_t$ to a compact nonsingular SL $m$-fold. The hypotheses in Joyce's theorem involves estimates of various kinds of norms of $\textnormal{Im}\,(\hat{\Oa}_t)|_{N_t}$. This suggest us to compute the term $\textnormal{Im}\,(\hat{\Oa}_t)$ restricting on different regions of $N_t$. \\


\par We shall introduce the notion of SL cones in \S\ref{sec:SpecialLagrangianConesAndTheirLagrangianNeighbourhoods}, SL $m$-folds with conical singularities in \S\ref{sec:SLMFoldsWithConicalSingularities} and AC SL $m$-folds in \S\ref{sec:ACSLMFolds}. Some examples of AC SL $m$-folds in AC \CY $m$-folds will be given in \S\ref{sec:SomeExamplesOfACSLMFolds}. 
We then proceed to \S\ref{sec:JoyceSDesingularizationTheory} by establishing necessary notations and discussing Joyce's analytic result. In \S\ref{sec:ConstructionOfNT} we construct Lagrangian 3-folds $N_t$ by gluing in $L_i$ to $N_0$ at each $x_i$. Then in \S\ref{sec:EstimatesOfImTildeOaTNT} we compute the estimates of the size of $\textnormal{Im}\,(\hat{\Oa}_t)|_{N_t}$. We divide the whole computation into three components, as given in equation (\ref{eqn:basicestimate}). Using the concept of \emph{local injectivity radius}, together with a kind of \emph{isoperimetric inequality} and the elliptic regularity result, we finally verify all the conditions, and so we are able to prove a result on SL desingularizations in \S\ref{sec:DesingularizationsOfN0}. In the last section, \S\ref{sec:ApplicationsOfTheDesingularizationTheory}, we illustrate our main result by taking an example of \CY 3-folds with conical singularities, namely the orbifold $T^6/\BZ_3$, and perform the desingularization simultaneously for both \CY and SL 3-folds with conical singularities. We use AC SL 3-folds from \S\ref{sec:SomeExamplesOfACSLMFolds} for gluing, thus obtaining various kinds of nonsingular SL 3-folds in the nonsingular \CY 3-folds. \\

\par Results related to the desingularization of SL $m$-folds can also be found, for instance, in Butscher \cite{But}, Joyce \cite{Joycedesing3}, \cite{Joycedesing4}, \cite{Joycedesing5}, and Lee \cite{Leeyi}. Butscher shows existence of SL connected sums of two compact SL $m$-folds (in $\BC^m$ with boundary) at one point by gluing in Lawlor necks (for definition, see \cite{Lawlor}). Joyce proves a desingularization result of SL $m$-folds with conical singularities in nonsingular \emph{almost} \CY $m$-folds by gluing in AC SL $m$-folds in $\BC^m$. He also studies desingularizations in families of almost \CY $m$-folds. Lee considers a compact, connected, immersed SL $m$-fold in a \CY $m$-fold, whose self-intersection points satisfy an angle criterion. She uses Lawlor necks for gluing at the singular points. \\

\noindent
\textsc{Acknowledgements.} \ The author is very grateful to Dominic Joyce for suggesting this problem and for his continued support during the author's D.Phil study in Oxford. The author would also like to thank Nigel Hitchin and Mark Haskins for their helpful advice. This work appeared as a section of the author's D.Phil thesis, and was supported by the Croucher Foundation in Hong Kong.  \\

\begin{center}
\renewcommand{\thesubsection}{\normalfont{2}}
\subsection{\normalfont{\textsc{Review of the main results in \cite{Chan1}}}}
\label{sec:review}
\end{center}


\centerline{2.1\ \textsc{\ Nearly \CY structure}}
\bigskip

\par Let $M$ be an oriented 6-fold. A nearly \CY structure on $M$ consists of a real closed 2-form $\oa$, and a complex closed 3-form $\Oa$ on $M$. Basically, the idea of a nearly \CY structure $(\oa, \Oa)$ is that it corresponds to an SU(3)-structure with ``small torsion", and hence approximates a genuine \CY structure, which is equivalent to a torsion-free SU(3)-structure. Here is the definition taken from \cite[Def. 3.1]{Chan1}: \\

\renewcommand{\thedefinition}{2.1}
\begin{definition}
\label{def:NCY}
\quad \textnormal{Let $M$ be an oriented 6-fold, and let $\oa$ be a real closed 2-form, and $\Oa$ a complex closed 3-form on $M$. Write $\Oa = \ta_1 + i\ta_2$, where $\ta_1$ and $\ta_2$ are both real closed 3-forms. Let $\ep_0 \in (0,1]$ be a fixed small constant. Then $(\oa, \Oa)$ constitutes a \emph{nearly \CY structure} on $M$ if}

\begin{itemize}
	\item[\textnormal{(i)}] \textnormal{the real closed 3-form $\ta_1$ has stabilizer SL(3,\,$\BC$) at each point of $M$,
\\[0.05cm] then we can associate a unique almost complex structure $J'$ and a unique real 3-form $\ta_2'$ such that $\Oa' = \ta_1 + i\ta_2'$ is a (3,0)-form with respect to $J'$; \\[-0.3cm]}
	
	\item[\textnormal{(ii)}] \textnormal{the (1,1)-component $\oal$ of $\oa$ with respect to $J'$ is positive, \\[0.05cm] then we can associate a Hermitian metric $g_M$ on $M$ from $J'$ and the rescaled (1,1)-part $\oa' := f^{-\frac{1}{3}}\,\oal $ of $\oa$, where $f : M \longrightarrow (0, \infty)$ is a smooth function defined by $(\oal)^3 = f \cdot \frac{3}{2}\,\ta_1 \we \ta_2'$; \\[-0.3cm]}
	
	\item[\textnormal{(iii)}] \textnormal{the following inequalities hold for some $\ep$ with $0< \ep \leq \ep_0$:}
	\begin{align} 
	\left| \ta_2 - \ta_2' \right|_{g_M} &< \ep	\, , \tag{2.1} \label{eqn:ncyforta2} \\ 
	| \oa^{(2,0)} |_{g_M} &< \ep	\, , \textnormal{and}  \tag{2.2} \label{eqn:ncyforoa20} \\
	\left| \oa^3 - \frac{3}{2}\, \ta_1 \we \ta_2 \right|_{g_M} &< \ep \tag{2.3} \label{eqn:ncyforoata} 
	\end{align} 
\textnormal{where the norms $|\cdot|_{g_M}$ are measured by the metric $g_M$. }	 \\[-0.5cm]
\end{itemize}
\end{definition} 

\par If $(\oa, \Oa)$ is a nearly \CY structure on $M$, one can show that the function $f$ satisfies $| f - 1 | < C_0 \ep$ for some constant $C_0 > 0$, i.e. $f$ is approximately equal to 1 for sufficiently small $\ep$, as we would expect. \\

\par From the above definition we know that $(J', \oa', \Oa')$ gives rise to an SU(3)-structure with metric $g_M$ on $M$, that is, $\oa'$ is of type (1,1) w.r.t. $J'$ and is positive, $\Oa'$ is of type (3,0) w.r.t. $J'$ and the normalization formula $\oa'^3 = \frac{3i}{4}\,\Oa' \we \bar{\Oa}'\, (= \frac{3}{2}\,\ta_1 \we \ta_2')$ holds. \\ 

\par The next result \cite[Prop. 3.2]{Chan1} shows that if $M$ admits a genuine \CY structure, then any real closed 2-form $\oa$ and complex closed 3-form $\Oa$ on $M$ which are sufficiently close to the genuine \CY structure gives a nearly \CY structure. We will use the notation $g^{-1}$ to denote the induced metric on the cotangent spaces.  \\

\renewcommand{\theprop}{2.2} 
\begin{prop}
\label{prop:closetogenuineisncy}
\quad There exist constants $\ep_1, C, C' >0$ such that whenever $0 < \ep \leq \ep_1$, the following is true. \\[-0.5cm]
\par Let $M$ be an oriented 6-fold. Suppose $(\tilde{J}, \tilde{\oa}, \tilde{\Oa})$ is a \CY structure with \CY metric $\tilde{g}$, $\oa$ a real closed 2-form, and $\Oa = \ta_1 + i \ta_2$ a complex closed 3-form on $M$, satisfying 
	\begin{align}
	| \tilde{\oa} - \oa |_{\tilde{g}} < \ep \quad and \quad | \tilde{\Oa} - \Oa |_{\tilde{g}} < \ep, \tag{2.4}
	\label{eqn:condforclosetogenuineisncy}
	\end{align}
then $(\oa, \Oa)$ is a nearly \CY structure on $M$ with metric $g_M$ satisfying 
	\begin{align}
	| \tilde{g} - g_M |_{\tilde{g}} < C \ep \quad and \quad | \tilde{g}^{-1} - g^{-1}_{M} |_{\tilde{g}} < C' \ep, \tag{2.5}
	\label{eqn:resultforclosetogenuineisncy}
	\end{align}	
\end{prop}

\hfill

\par Let $(\oa, \Oa)$ be a nearly \CY structure on $M$. Then $(J', \oa', \Oa')$ is an SU(3)-structure with metric $g_M$ on $M$. Now we consider the 7-dimensional manifold $\SM$, with $s$ being a coordinate on $S^1$. From the SU(3)-structure $(J', \oa', \Oa')$ one can induce a $G_2$-structure (with torsion) on $\SM$ by the following construction. Define a 3-form $\vi'$ and a metric $g'$ on $\SM$ by
	\begin{align} 
	\vi' = ds \we \oa' + \ta_1 \quad \textnormal{and} \quad g' = ds^2 + g_M.  \tag{2.6}
	\label{eqn:vi'g'}
	\end{align}
For an introduction to $G_2$-structures, see for example \cite[Chapter 10]{Joyce1}. 
The associated 4-form $\ast_{g'} \vi'$ on $\SM$ is then given by 
	\begin{align} 
	\ast_{g'} \vi' = \frac{1}{2}\,\oa' \we \oa' - ds \we \ta_2'. \tag{2.7}
	\label{eqn:*vi'}
	\end{align}
Also, we can construct from the nearly \CY structure $(\oa, \Oa)$ another 3-form $\vi$ and 4-form $\chi$ on $\SM$ by 
	\begin{align} 
	\vi = ds \we \oa + \ta_1 \quad \textnormal{and} \quad \chi = \frac{1}{2}\,\oa \we \oa - ds \we \ta_2. \tag{2.8}
	\label{eqn:vichi}
	\end{align}

The next lemma \cite[Lem. 3.3]{Chan1} shows that the forms in (\ref{eqn:vichi}) are close to the $G_2$-forms $\vi'$ and $\ast_{g'} \vi'$ if we take $\ep_0$ in the definition of nearly \CY manifolds to be sufficiently small. \\

\renewcommand{\thelem}{2.3} 
\begin{lem}
\label{lem:vigisG2}
\quad There exist constants $C_1, C_2, C_3$ and $C_4 > 0$ such that if $\ep_0$ in Definition \ref{def:NCY} is chosen sufficiently small, then the following is true. \\[-0.2cm]

\par Let $(\vi', g')$ be the $G_2$-structure given by (\ref{eqn:vi'g'}), $\ast_{g'} \vi'$ the associated 4-form given by (\ref{eqn:*vi'}), $\vi$ the 3-form and $\chi$ the 4-form given by (\ref{eqn:vichi}) on $\SM$. Then
  \begin{align}	
	|\vi - \vi'|_{g'} < C_1 \ep \tag{2.9}
	\label{eqn:vi-vi'}
	\end{align}
where $\ep \in (0, \ep_0]$ is the small constant in (iii) of Definition \ref{def:NCY}. Hence $\vi$ is also a positive 3-form on $\SM$, and it defines another $G_2$-structure $(\vi, g)$. Moreover, the associated metric $g$ and the 4-form $\ast_{g}\vi$ satisfy 
	\begin{align} 
	|g - g' |_{g'} < C_2 \ep, \quad |g^{-1} - g'^{-1} |_{g'} &< C_3\,\ep \ \ and \tag{2.10} \label{eqn:resultg-g'} \\
	|\!\ast_{g}\!\vi - \chi |_{g'} < C_4\,\ep. \tag{2.11} \label{eqn:result*vi-chi}
	\end{align} \\[-0.5cm]
\end{lem}

\begin{center}
\renewcommand{\thesubsubsection}{\normalfont{2.2}}
\subsubsection{\normalfont{\textsc{An existence result for \CY structures on $M$}}}
\end{center}

\hfill

\par We present here the main analytic result \cite[Thm. 3.10]{Chan1} on deforming the nearly \CY structure $(\oa, \Oa)$ to a genuine \CY structure on $M$, given that $\ep_0$ is small enough. The proof is based on an existence result for torsion-free $G_2$-structures given by Joyce \cite[Thm. 11.6.1]{Joyce1}, which shows using analysis that any $G_2$-structure on a compact 7-fold with sufficiently small torsion can be deformed to a nearby torsion-free $G_2$-structure. The following \cite[Thm. 3.9]{Chan1} is a modified version of Joyce's result \cite[Thm. 11.6.1]{Joyce1}, giving the existence of an $S^1$-invariant torsion-free $G_2$-structure on the 7-fold $\SM$: \\

\renewcommand{\thethm}{2.4} 
\begin{thm}
\label{thm:torsionfreeG2onS1xM}
\quad Let $\kappa >0$ and $D_1, D_2, D_3 > 0$ be constants. Then there exist constants $\ep \in (0,1]$ and $K > 0$ such that whenever $0 < t \leq \ep$, the following is true. \\[-0.2cm]
\par Let $M$ be a compact 6-fold, and $(\vi, g)$ an $S^1$-invariant $G_2$-structure on $\SM$ with $d\vi = 0$. Suppose $\psi$ is an $S^1$-invariant smooth 3-form on $\SM$ with $d^* \psi = d^* \vi$, and
\begin{itemize}
	\item[\textnormal{(i)}] $\| \psi \|_{L^2} \leq D_1 t^{3+\kappa}$, $\| \psi \|_{C^0} \leq D_1 t^{\kappa}$\ \,and\ \,$\| d^* \psi \|_{L^{12}} \leq D_1 t^{-\frac{1}{2}+\kappa}$,
	\item[\textnormal{(ii)}] the injectivity radius $\da(g)$ satisfies $\da(g) \geq D_2 t$, and
	\item[\textnormal{(iii)}] the Riemann curvature $R(g)$ satisfies $\| R(g) \|_{C^0} \leq D_3 t^{-2}$. \\[-0.5cm]
\end{itemize}
Then there exists a smooth, torsion-free $S^1$-invariant $G_2$-structure $(\tilde{\vi}, \tilde{g})$ on $\SM$ such that $\| \tilde{\vi} - \vi \|_{C^0} \leq K t^{\kappa}$\, and\, $[\tilde{\vi}] = [\vi]$ in $H^3 (\SM, \BR)$. \\ 
\end{thm}

\par The proof of Theorem \ref{thm:torsionfreeG2onS1xM} depends upon the two results in \cite[Thm. 3.7 \& 3.8]{Chan1}. We mention here that we shall need improved versions of those two results to obtain a better estimate of $\| \tilde{\vi} - \vi \|_{C^0}$ for the special Lagrangian desingularization (see \S\ref{sec:EstimatesOfImTildeOaTNT}). We refer to the spaces $L^q$, $L^{q}_{k}$, $C^{k}$ and $C^{k, \al}$ as the Banach spaces defined in \cite[\S1.2]{Joyce1}. \\

\par With the help of Theorem \ref{thm:torsionfreeG2onS1xM}, we present the following existence result \cite[Thm. 3.10]{Chan1} for genuine \CY structures. The basic idea is that we start with the nearly \CY structure $(\oa, \Oa)$ on $M$, then one can induce a $G_2$-structure $(\vi, g)$ on $\SM$. Under appropriate hypotheses on the nearly \CY structure $(\oa, \Oa)$, the induced $G_2$-structure satisfies all the conditions in Theorem \ref{thm:torsionfreeG2onS1xM}, and therefore can be deformed to have zero torsion. The \CY structure on $M$ can then be obtained by pulling back this torsion-free $G_2$-structure. \\

\renewcommand{\thethm}{2.5} 
\begin{thm}
\label{thm:existencegenuineCY}
\quad Let $\kappa > 0$ and $E_1, E_2, E_3, E_4 > 0$ be constants. Then there exist constants $\ep \in (0,1]$ and $K > 0$ such that whenever $0 < t \leq \ep$, the following is true. \\[-0.4cm]
\par Let $M$ be a compact 6-fold, and $(\oa, \Oa)$ a nearly \CY structure on $M$. Let $\oa'$, $g_M$ and $\ta_2'$ be as in Definition \ref{def:NCY}. Suppose
\begin{itemize}
	\item[\textnormal{(i)}] $\| \oa - \oa' \|_{L^2} \leq E_1 t^{3+\kappa},\ \,\| \oa - \oa' \|_{C^0} \leq E_1 t^{\kappa},\ \| \ta_2 - \ta_2' \|_{L^2} \leq E_1 t^{3+\kappa}$, \\[0.1cm]
and\ \,$\| \ta_2 - \ta_2' \|_{C^0} \leq E_1 t^{\kappa}$,
	\item[\textnormal{(ii)}] $\| \nab (\oa - \oa') \|_{L^{12}} \leq E_1 t^{-\frac{1}{2}+\kappa} $,\ \,$\| \nab \oa \|_{L^{12}} \leq E_1 t^{-\frac{1}{2}+\kappa}$, \\[0.1cm]
and\ \,$\| \nab \ta_1 \|_{L^{12}} \leq E_1 t^{-\frac{1}{2}+\kappa}$,  
	\item[\textnormal{(iii)}] $\| \nab (\oa - \oa') \|_{C^0} \leq E_2 t^{\kappa - 1}$,\ \,and\ \,$\| \nab^{2} (\oa - \oa') \|_{C^0} \leq E_2 t^{\kappa - 2}$,
	\item[\textnormal{(iv)}] the injectivity radius $\da(g_M)$ satisfies $\da(g_M) \geq E_3 t$, and
	\item[\textnormal{(v)}] the Riemann curvature $R(g_M)$ satisfies $\| R(g_M) \|_{C^0} \leq E_4 t^{-2}$. \\[-0.5cm]
 \end{itemize} 
Then there exists a \CY structure $(\tilde{J}, \tilde{\oa}, \tilde{\Oa})$ on $M$ such that $\| \tilde{\oa} - \oa \|_{C^0} \leq K t^{\kappa}$ and\, $\| \tilde{\Oa} - \Oa \|_{C^0} \leq K t^{\kappa}$. Moreover, if $H^1 (M, \BR) = 0$, then the cohomology classes satisfy $[\textnormal{Re}(\Oa)]$ = $[\textnormal{Re}(\tilde{\Oa})] \in H^3 (M, \BR)$\, and\, $[\oa]$ = $c\,[\tilde{\oa}] \in H^2 (M, \BR)$ for some $c > 0$. Here the connection $\nab$ and all norms are computed with respect to $g_M$. \\
\end{thm}

\begin{center}
\renewcommand{\thesubsubsection}{\normalfont{2.3}}
\subsubsection{\normalfont{\textsc{\CY cones, \CY manifolds with conical singularities and Asymptotically Conical \CY manifolds}}}
\label{subsub:\CY cones\CY manifolds with conical}
\end{center}

\hfill

\par Here we summarize \S4 of \cite{Chan1} in which \CY cones, \CY manifolds with a conical singularity and Asymptotically Conical (AC) \CY manifolds are defined. Note that in this paper we shall extend the situation in \cite{Chan1} by considering \CY manifolds with finitely many conical singularities, and we always use $n$ to denote the number of these singular points. Moreover, we shall make the following definitions in $m$ complex dimensions for now; later $m$ will be 3. \\

\par We start with the definition of \CY cones \cite[Def. 4.1]{Chan1}: \\

\renewcommand{\thedefinition}{2.6}
\begin{definition}
\label{def:CYcone}
\quad \textnormal{Let $\Ga$ be a compact ($2m-1$)-dimensional smooth manifold, and let $V = \{0\} \cup V'$ where $V' = \Ga \times (0, \infty)$. Write points on $V'$ as $(\ga, r)$. $V$ is called a \emph{\CY cone} if $V'$ is a \CY $m$-fold with a \CY structure $(J_{V}, \oa_{V} , \Oa_{V})$ and its associated \CY metric $g_{V}$ satisfying 
\begin{align}
g_{V} =\, &r^2 g_{V} |_{\Ga \times \{1\}} + dr^2, \quad \oa_{V} = r^2 \oa_{V} |_{\Ga \times \{1\}} + rdr \we \al \notag \\[-0.3cm]
\tag{2.12} \label{eqn:CYconeforms} \\[-0.3cm]
&\textnormal{and}\quad \Oa_{V} = r^m \Oa_{V} |_{\Ga \times \{1\}} + r^{m-1}dr \we \ba. \notag
\end{align}
Here we identify $\Ga$ with $\Ga \times \{1\}$, and $\al$ is a real 1-form and $\ba$ a complex $(m-1)$-form on $\Ga$. We also call $\Ga$ the \emph{link} of $V$.
} \\
\end{definition}

\par In terms of \emph{Sasaki--Einstein geometry}, a Riemannian manifold $(M,g)$ of dimension $(2m-1)$ is Sasaki--Einstein if and only if the cone over $M$ with metric $r^2g + dr^2$ is Calabi--Yau, i.e. a \CY cone. Thus in our case, $V$ is a \CY cone is equivalent to $\Ga$ being Sasaki--Einstein. The canonical example of a Sasaki--Einstein manifold is given by an odd-dimensional sphere $S^{2m-1}$. The \CY cone over it is then $\BC^m$ with its flat metric. There has been considerable interest recently in Sasaki--Einstein geometry due to a new construction of an infinite family of explicit Sasaki--Einstein metrics on five dimensions, particularly on $S^2 \times S^3$ \cite{Gauntlett}. Much work has been done by Boyer and Galicki on Sasaki--Einstein and 3-Sasakian geometry, see for examples \cite{BoyerGalicki1} and \cite{BoyerGalicki2}. \\

\par One can define the radial vector field $X$ on the \CY cone $V$, the contact 1-form $\al$ on the link $\Ga$ and the Reeb vector field $Z = J_V X$ on $\Ga$ as in \S4.1 of \cite{Chan1}. The flow of $Z$ generates the diffeomorphism $\textnormal{exp} (\ta Z)$ on $\Ga$ for each $\ta \in \BR$. Thus for each $\ta \in \BR$ and $t>0$, there is a complex dilation $\lda$ on $V$ given by $\lda(0) = 0$ and $\lda(\ga, r) = (\textnormal{exp} (\ta Z)(\ga),\, tr)$. The scaling behaviour of $\lda$ is given in \cite[Lem. 4.2]{Chan1}. \\

\par Here is a nontrivial example of a \CY cone, taken from \cite[Example 4.3]{Chan1}: \\

\renewcommand{\theeg}{2.7}
\begin{eg} 
\label{eg:Cm/G}
\quad \textnormal{Let $G$ be a finite subgroup of SU($m$) acting freely on $\BC^m \setminus \{0\}$, then the quotient singularity $\BC^m /G$ is a \CY cone. An example of this type is given by the $\BZ_m$-action: define an action generated by $\zeta$ on $\BC^m$ by 
\begin{align*}
\zeta^k \cdot (z_1 , \dots , z_m) \ = \ (\zeta^k \,z_1 , \dots , \zeta^k \,z_m)
\end{align*}
where $\zeta = e^{2 \pi i / m}$ and $0 \leq k \leq m-1$. Note that $\zeta^m = 1$, so $\BZ_m = \{ 1, \zeta, \dots, \zeta^{m-1} \}$ is a subgroup of SU($m$) and acts freely on $\BC^m \setminus\!\{0\}$. Then $\BC^m  / \mathbb{Z}_m$ is a \CY cone, with the Sasaki--Einstein link $S^{2m-1}/\BZ_m$.} \\
\end{eg}

\par Next we define \CY $m$-folds with finitely many conical singularities \cite[Def. 4.6]{Chan1}. As we have mentioned, we shall consider the case when there are not only one but $n$ conical singular points. Note that in this paper we always assume the existence of \CY metrics on such kind of singular manifolds. A class of \CY $m$-folds with conical singularities are given by orbifolds, in which case the existence of such singular \CY metrics is known (see \cite[Thm. 6.5.6]{Joyce1}). \\ 

\renewcommand{\thedefinition}{2.8}
\begin{definition}
\label{def:CYmfoldcs}
\quad \textnormal{Let $(M_0, J_0, \oa_0, \Oa_0)$ be a singular \CY $m$-fold with isolated singularities $x_1, \dots, x_n \in M_0$, and no other singularities. We say that $M_0$ is a \emph{\CY $m$-fold with conical singularities $x_i$} for $i=1, \dots,n$ with rate $\nu > 0$ modelled on \CY cones $(V_i, J_{V_i}, \oa_{V_i} , \Oa_{V_i})$ if there exist a small $\ep > 0$, a small open neighbourhood $S_i$ of $x_i$ in $M_0$, and a diffeomorphism $\Phi_i : \Ga_i \times (0, \ep) \longrightarrow S_i \setminus \{x_i\}$ for each $i$ such that 
\begin{align}
|\nab^k ( \Phi_i^* ( \oa_0 ) - \oa_{V_i} ) |_{g_{V_i}} &= O(r^{\nu - k}), \ \textnormal{and} \tag{2.13} \label{eqn:defofcsforoa} \\
|\nab^k ( \Phi_i^* ( \Oa_0 ) - \Oa_{V_i} ) |_{g_{V_i}} &= O(r^{\nu - k}) \ \ \textnormal{as } r \rightarrow 0 \ \textnormal{and for all } k \geq 0. \tag{2.14} \label{eqn:defofcsforOa}
\end{align}
Here $\Ga_i$ is the link of $V_i$, and $\nab$, $|\cdot|_{g_{V_i}}$ are computed using the cone metric $g_{V_i}$.} \\
\end{definition}

\par Note that the asymptotic conditions on $g_0$ and $J_0$ follow from those on $\oa_0$ and $\Oa_0$, namely,
\begin{align}
|\nab^k ( \Phi_i^* ( g_0 ) - g_{V_i} ) |_{g_{V_i}} &= O(r^{\nu - k}), \ \textnormal{and} \tag{2.15} \label{eqn:defofcsforg} \\
|\nab^k ( \Phi_i^* ( J_0 ) - J_{V_i} ) |_{g_{V_i}} &= O(r^{\nu - k}) \ \ \textnormal{as } r \rightarrow 0 \ \textnormal{and for all } k \geq 0, \tag{2.16} \label{eqn:defofcsforJ}
\end{align}
and so it is enough to just assume asymptotic conditions on $\oa_0$ and $\Oa_0$. \\

\par We will usually assume that $M_0$ is compact. The point of the definition is that $M_0$ is locally modelled on $\Ga_i \times (0, \ep)$ near $x_i$, and as $r \rightarrow 0$, all the structures $g_0, J_0, \oa_0$ and $\Oa_0$ on $M_0$ converge to the cone structures $g_{V_i}, J_{V_i}, \oa_{V_i}$ and $\Oa_{V_i}$ with rate $\nu$ and with all their derivatives. \\

\par Two diffeomorphisms, or two coordinate systems, $\Phi_i$ and $\Phi_i'$ are equivalent if and only if the following relation holds: 
\begin{align*}
| \nab^k ( \Phi_i - \Phi_i' ) |_{g_{V_i}} = O(r^{\nu + 1 - k}) \quad \textnormal{as $r \rightarrow 0$ and for all $k \geq 0$.} \tag{2.17} \label{eqn:equivalentCYcoord}
\end{align*}
Here we interpret the difference between $\Phi_i$ and $\Phi_i'$ using local coordinates on the image $S_i \setminus \{ x_i\}$. Thus if $\Phi_i$ and $\Phi_i'$ are equivalent, we have
\begin{align*}
|\nab^k ( \Phi_i^* ( \oa_0 ) - \oa_{V_i} ) |_{g_{V_i}} &\leq | \nab^k ( \Phi_i^* (\oa_0) - (\Phi_i')^* (\oa_0) ) |_{g_{V_i}} \, + \, |\nab^k ( (\Phi_i')^* ( \oa_0 ) - \oa_{V_i} ) |_{g_{V_i}} \\
&= O(r^{\nu-k}) \, + \, |\nab^k ( (\Phi_i')^* ( \oa_0 ) - \oa_{V_i} ) |_{g_{V_i}} 
\end{align*}
and we see that $\Phi_i$ satisfies (\ref{eqn:defofcsforoa}) (and similarly (\ref{eqn:defofcsforOa})) if and only if $\Phi_i'$ does. \\

\par The 2-forms $\Phi_i^* ( \oa_0 )$ and $\oa_{V_i}$ are closed on $\Ga_i \times (0,\ep)$ and so $\Phi_i^* ( \oa_0 ) - \oa_{V_i}$ represents a cohomology class in $H^2 (\Ga_i \times (0,\ep), \BR) \cong H^2 (\Ga_i, \BR)$. Similarly, $\Phi_i^* ( \Oa_0 ) - \Oa_{V_i}$ represents a cohomology class in $H^m (\Ga_i \times (0,\ep), \BC) \cong H^m (\Ga_i, \BC)$. It follows from \cite[Lem. 4.7]{Chan1} that in the conically singular case, these two classes $[\Phi_i^* ( \oa_0 ) - \oa_{V_i}]$ and $[\Phi_i^* ( \Oa_0 ) - \Oa_{V_i}]$ are automatically zero. Moreover, there is a Darboux theorem in  \cite[Thm. 4.9]{Chan1} asserting that locally the symplectic form $\oa_0$ near the conical singular point $x_i$ is symplectomorphic to the cone form $\oa_{V_i}$ near the origin. \\

\par Now we proceed to define Asymptotically Conical (AC) \CY $m$-folds \cite[Def. 4.10]{Chan1}.  \\

\renewcommand{\thedefinition}{2.9}
\begin{definition}
\label{def:ACCYmfold}
\quad \textnormal{Let $(V, J_{V}, \oa_{V}, \Oa_{V})$ be a \CY cone of complex dimension $m$ with link $\Ga$ which is a compact Sasaki--Einstein $(2m-1)$-fold. Let $(Y, J_{Y}, \oa_{Y}, \Oa_{Y})$ be a complete, nonsingular \CY $m$-fold. Then $Y$ is an \emph{Asymptotically Conical (AC) \CY $m$-fold} with rate $\lda < 0$ modelled on $(V, J_{V}, \oa_{V}, \Oa_{V})$ if there exist a compact subset $K \subset Y$, and a diffeomorphism $\Ups: \Ga \times (R, \infty) \longrightarrow Y \setminus K$ for some $R > 0$ such that
\begin{align}
|\nab^k ( \Ups^* ( \oa_{Y} ) - \oa_{V} ) |_{g_{V}} &= O(r^{\lda - k}), \ \textnormal{and} \tag{2.18} \label{eqn:ACdefofoa} \\ 
|\nab^k ( \Ups^* ( \Oa_{Y} ) - \Oa_{V} ) |_{g_{V}} &= O(r^{\lda - k}) \ \ \textnormal{as } r \rightarrow \infty \ \textnormal{and for all } k \geq 0. \tag{2.19} \label{eqn:ACdefofOa}
\end{align}
Here $\nab$ and $| \cdot |$ are computed using the cone metric $g_{V}$.} \\
\end{definition}
 
\par Similar asymptotic conditions on $g_{Y}$ and $J_{Y}$ can be deduced from (\ref{eqn:ACdefofoa}) and (\ref{eqn:ACdefofOa}). The coordinates $\Ups$ and $\Ups'$ are equivalent if and only if $| \nab^k ( \Ups - \Ups' ) |_{g_{V}} = O(r^{\lda + 1 - k})$ as $r \rightarrow \infty$ and for all $k \geq 0$. \\

\renewcommand{\theremark}{2.10}
\begin{remark}
\label{rmk:ACCYmfold}
\quad \textnormal{If $Y$ is an AC \CY $m$-fold which is not a $\BC^m$, then $Y$ can only have one end, or equivalently, the link $\Ga$ is connected. One can show this by using the \emph{Cheeger-Gromoll splitting theorem} (see for example \cite[\S6.G]{Besse}). Suppose $Y$ has more than one end. As $Y$ is complete and Ricci-flat, Cheeger-Gromoll splitting theorem tells us that we can always split off a line so that $Y$ is isometric to a product $N \times \BC$, where $\BC$ carries the Euclidean metric. Now, either $N$ is a flat $\BC^{m-1}$, in which case $Y = \BC^m$, or $N$ has nonzero curvature at some $p \in N$, in which case the curvature of $N \times \BC$ is of order $O(1)$ as we go to infinity in $\{ p \} \times \BC$. But then this contradicts the AC condition which requires the curvature to decay at $O(r^{-2})$. Therefore $Y$ cannot have more than one end, and so from now on, 
we shall always take the link $\Ga$ to be a \emph{compact, connected $\it (2m-1)$-dimensional Sasaki--Einstein manifold.}}    \\
\end{remark}

\par Unlike the conical singularity case, $[ \Ups^* ( \oa_{Y} ) - \oa_{V} ]$ and $[ \Ups^* ( \Oa_{Y} ) - \Oa_{V} ]$ need not be zero cohomology classes. The result in \cite[Lem. 4.11]{Chan1} shows that if $\lda < -2$ or $H^2 (\Ga, \BR) = 0$, then $[ \Ups^* ( \oa_{Y} ) - \oa_{V} ] = 0$; if $\lda < -m$ or $H^m (\Ga, \BC) = 0$, then $[ \Ups^* ( \Oa_{Y} ) - \Oa_{V} ] = 0$. There is an analogue of Darboux Theorem \cite[Thm. 4.12]{Chan1} for AC \CY $m$-folds with rate $\lda < -2$. \\

\par In this paper we shall consider the simplest case $\lda < -3$ so that $\Ups^* ( \oa_{Y} ) - \oa_{V}$ and $\Ups^* ( \Oa_{Y} ) - \Oa_{V}$ are always exact for AC \CY 3-folds. \\

\par Here is an example of AC \CY $m$-fold, taken from \cite[Example. 4.13]{Chan1}: \\

\renewcommand{\theeg}{2.11}
\begin{eg} 
\label{eg:ACforCm/G}
\quad \textnormal{Let $G$ be a finite subgroup of SU($m$) acting freely on $\BC^m \setminus \{0\}$, and $(X, \pi)$ a crepant resolution of the \CY cone $V = \BC^m / G$ given in Example \ref{eg:Cm/G}. Then in each \ka class of ALE \ka metrics on $X$ there is a unique ALE Ricci-flat \ka metric (see Joyce \cite{Joyce1}, Chapter 8) and $X$ is then an AC \CY $m$-fold asymptotic to the cone $\BC^m / G$. In this case, it follows from \cite[Thm. 8.2.3]{Joyce1} that the rate $\lda$ is $-2m$.} \\
\par \textnormal{If we take $G = \mathbb{Z}_m$ acting on $\BC^m$ as in Example \ref{eg:Cm/G}, then a crepant resolution is given by the blow-up of $\BC^m / \mathbb{Z}_m$ at 0, which is also the total space of the canonical line bundle $K_{\BC \BP^{m-1}}$ over $\BC \BP^{m-1}$. An explicit ALE Ricci-flat \ka metric is given in \cite[p.284-5]{Calabi} and also in \cite[Example 8.2.5]{Joyce1}. } \\
\end{eg}

\begin{center}
\renewcommand{\thesubsubsection}{\normalfont{2.4}}
\subsubsection{\normalfont{\textsc{Theorem on desingularizing \CY 3-folds}}}
\label{subsub:Theorem on desingularizing \CY $m$-folds}
\end{center}

\hfill

\par We study the desingularization of a compact \CY 3-fold $M_0$ with conical singularities using an AC \CY 3-fold $Y_i$ with rate $\lda_i < -3$ for $i = 1, \dots, n$. We explicitly construct a 1-parameter family of diffeomorphic, nonsingular compact 6-folds $M_t$ for small $t$, and construct a real closed 2-form $\oa_t$ and a complex closed 3-form $\Oa_t$ on $M_t$ so that $(\oa_t, \Oa_t)$ gives a nearly \CY structure on $M_t$ for small enough $t$. Finally we state the main result on \CY desingularizations \cite[Thm. 5.2]{Chan1}, asserting that the nearly \CY structure $(\oa_t, \Oa_t)$ on $M_t$ can be deformed to a genuine \CY structure $(\tilde{\oa}_t, \tilde{\Oa}_t)$ for small $t$, using Theorem \ref{thm:existencegenuineCY}. \\

\par Let $(M_0 , J_0 , \oa_0 , \Oa_0 )$ be a compact \CY 3-fold with conical singularities $x_i$ with rate $\nu$ modelled on \CY cones $(V_i, J_{V_i} , \oa_{V_i}, \Oa_{V_i})$ for $i = 1, \dots, n$. Then there exists an $\ep > 0$, a small open neighbourhood $S_i$ of $x_i$ in $M_0$ and a diffeomorphism $\Phi_i : \Ga_i \times (0, \ep) \longrightarrow S_i \setminus \{x_i\}$ for each $i$ such that $\Phi_i^* (\oa_0 ) = \oa_{V_i}$. \\

\par Let $(Y_i, J_{Y_i} , \oa_{Y_i}, \Oa_{Y_i})$ be an AC \CY 3-fold with rate $\lda_i < -3$ modelled on the same \CY cone $V_i$. Then there is a diffeomorphism $ \Ups_i : \Ga_i \times (R, \infty) \longrightarrow Y_i \setminus\!K_i$ for some $R > 0$ such that 
$$\Ups_i^* (\oa_{Y_i}) = \oa_{V_i} \quad \textnormal{and} \quad |\nab^k (\Ups_i^* (\Oa_{Y_i}) - \Oa_{V_i} ) |_{g_{V_i}} = O(r^{\lda_i - k})$$ 
as $r \rightarrow \infty$ for all $k \geq 0$. We then apply a homothety to $Y_i$ such that 
$$ (Y_i, J_{Y_i} , \oa_{Y_i}, \Oa_{Y_i}) \longmapsto (Y_i, J_{Y_i} , t^2 \oa_{Y_i},  t^3 \Oa_{Y_i}). $$
Then $(Y_i, J_{Y_i} , t^2 \oa_{Y_i},  t^3 \Oa_{Y_i})$ is also an AC \CY 3-fold, with the diffeomorphism $\Ups_{t,i} : \Ga_i \times (tR, \infty) \longrightarrow Y_i \setminus\!K_i$ given by  
$$ \Ups_{t,i} (\ga, r) = \Ups_i (\ga, t^{-1}r).$$

\par Fix $\al \in (0,1)$ and let $t > 0$ be small enough that $tR < t^{\al} < 2t^{\al} < \ep$. Define
\begin{align*}
P_{t,i} &= K_i \cup \Ups_{t,i} ( \Ga_i \times (tR, 2t^{\al}) ) \subset Y_i \quad \textnormal{and} \\
Q_t &= M_0 \setminus \bigcup_{i=1}^n \Phi_i(\Ga_i \times (0, t^{\al}) ) \subset M_0
\end{align*} 
The diffeomorphism $\Phi_i \circ \Ups_{t,i}^{-1}$ identifies $\Ups_{t,i} ( \Ga_i \times (t^{\al}, 2t^{\al}) ) \subset P_{t,i}$ and $\Phi_i(\Ga_i \times (t^{\al}, 2t^{\al}) ) \subset Q_t$, and we define the intersection $P_{t,i} \cap Q_t$ to be the region $\Ups_{t,i} ( \Ga_i \times (t^{\al}, 2t^{\al}) ) \cong \Phi_i(\Ga_i \times (t^{\al}, 2t^{\al}) ) \cong \Ga_i \times (t^{\al}, 2t^{\al})$. Define $M_t$ to be the quotient space of the union $(\bigcup_{i=1}^n P_{t,i}) \cup Q_t$ under the equivalence relation identifying the two annuli $\Ups_{t,i} ( \Ga_i \times (t^{\al}, 2t^{\al}) )$ and $\Phi_i(\Ga_i \times (t^{\al}, 2t^{\al}) )$ for $i = 1, \dots, n$. Then $M_t$ is a smooth nonsingular compact 6-fold for each $t$. \\

\par Now we construct on $M_t$ a real closed 2-form $\oa_t$ and a complex closed 3-form $\Oa_t$, and show that they together give nearly \CY structures on $M_t$ for small enough $t$. Define 
\begin{center}
$\oa_t$ = 
$\begin{cases}$
$\oa_0$ \ \ on\ \ $Q_t$, $\\$
$t^2 \oa_{Y_i}$\ \ on\ \ $P_{t,i} \ \textnormal{for } i = 1, \dots, n$.
$\end{cases}$
\end{center}
This is well-defined as $\Phi_i^* (\oa_0) = \oa_{V_i} = \Ups_{t,i}^{*} (t^2 \oa_{Y_i})$ on each intersection $P_{t,i} \cap Q_t$. Thus $\oa_t$ gives a symplectic form on $M_t$. \\

\par Let $F: \BR \longrightarrow [0,1]$ be a smooth, increasing function with $F(s) = 0$ for $s \leq 1$ and $F(s) = 1$ for $s \geq 2$. Then for $r \in (tR, \ep)$, $F(t^{-\al}r) = 0$ for $tR < r \leq t^{\al}$ and $F(t^{-\al}r) = 1$ for $2t^{\al} \leq r < \ep$. We now define a complex 3-form on $M_t$. From (\ref{eqn:defofcsforOa}), we have $|\Phi_i^* (\Oa_0) - \Oa_{V_i} |_{g_{V_i}} = O(r^{\nu})$. As $\nu > 0$, it follows that $\Phi_i^* (\Oa_0) - \Oa_{V_i}$ is exact, and we can write 
\begin{align}
\Phi_i^* (\Oa_0) = \Oa_{V_i} + dA_i \tag{2.20}
\label{eqn:Phi*(Oa0)=OaV+dA}
\end{align}
for some complex 2-form $A_i(\ga, r)$ on $\Ga_i \times (0, \ep)$ satisfying 
\begin{align}
| \nab^k A_i(\ga, r)|_{g_{V_i}} = O(r^{\nu + 1 -k}) \quad \textnormal{as } r \rightarrow 0 \  \textnormal{ for all } k \geq 0.  \tag{2.21}
\label{eqn:asympbehavefornabkA}
\end{align}
Similarly, as we have assumed $\lda_i < -3$ to simplify the problem, the 3-form $\Ups_i^* (\Oa_{Y_i}) - \Oa_{V_i}$ is exact and we can write
$$\Ups_i^* (\Oa_{Y_i}) = \Oa_{V_i} + dB_i $$
for some complex 2-form $B_i(\ga, r)$ on $\Ga_i \times (R, \infty)$ satisfying 
$$| \nab^k B_i(\ga, r)|_{g_{V_i}} = O(r^{\lda_i + 1 -k}) \quad \textnormal{as } r \rightarrow \infty \textnormal{ and for all } k \geq 0. $$ 
Then we apply a homothety to $Y_i$ and rescale the forms to get $B_i(\ga, t^{-1}r)$ on $\Ga_i \times (tR, \infty)$ such that 
\begin{align}
\Ups_{t,i}^{*} (t^3 \Oa_{Y_i}) = \Oa_{V_i} + t^3 dB_i(\ga, t^{-1}r) \tag{2.22}
\label{eqn:Upst*(t3OaY)=OaV+t3dB}
\end{align}
and
\begin{align}
| \nab^k B_i(\ga, t^{-1} r)|_{g_{V_i}} = O(t^{-\lda_i-3} r^{\lda_i + 1 -k}) \quad \textnormal{for } r > tR \textnormal{ and for all } k \geq 0. \tag{2.23}
\label{eqn:asympbehavefornabkB}
\end{align}
Define a smooth, complex closed 3-form $\Oa_t$ on $M_t$ by \\[-0.3cm]
\begin{align}
\Oa_t = 
\begin{cases}
\Oa_0$ \ \ on\ \ $Q_t \setminus \big[ (\bigcup_{i=1}^n P_{t,i}) \cap Q_t \big], \notag \\
\Oa_{V_i} + d\big[ F(t^{-\al}r) A_i(\ga, r) + t^3 (1 - F(t^{-\al}r)) B_i(\ga, t^{-1}r)\big]$ on $P_{t,i} \cap Q_t, \textnormal{for } i = 1, \dots, n, \tag{2.24} \label{eqn:defofOat} \\ 
t^3 \Oa_{Y_i}$\ \ on\ \ $P_{t,i} \setminus (P_{t,i} \cap Q_t)\ \textnormal{for } i = 1, \dots, n.   
\end{cases} \\[-0.5cm]  \notag
\end{align} 
Note that when $2t^{\al} \leq r < \ep$ we have $F(t^{-\al}r) = 1$ so that $\Oa_t = \Phi_i^* (\Oa_0)$ by (\ref{eqn:Phi*(Oa0)=OaV+dA}), and when $tR < r \leq t^{\al}$ we have $F(t^{-\al}r) = 0$, so that $\Oa_t = \Ups_{t,i}^{*} (t^3 \Oa_{Y_i})$ by (\ref{eqn:Upst*(t3OaY)=OaV+t3dB}). Therefore, $\Oa_t$ interpolates between $\Phi_i^* (\Oa_0)$ near $r = \ep$ and $\Ups_{t,i}^{*} (t^3 \Oa_{Y_i})$ near $r = tR$. \\

\par Recall that if we get a real closed 2-form $\oa$ and a complex 3-form $\Oa$ which are sufficiently close to the \ka form $\tilde{\oa}$ and the holomorphic volume form $\tilde{\Oa}$ of a \CY structure respectively, then Proposition \ref{prop:closetogenuineisncy} tells us that $(\oa, \Oa)$ gives a nearly \CY structure on the manifold. Making use of this idea, we have the following \cite[Prop. 5.1]{Chan1}: \\

\renewcommand{\theprop}{2.12} 
\begin{prop}
\label{prop:NCYforMt}
\quad Let $M_t$, $\oa_t$ and $\Oa_t$ be defined as above. Then $(\oa_t, \Oa_t)$ gives a nearly \CY structure on $M_t$ for sufficiently small t. \\
\end{prop}

\par Therefore we can associate an almost complex structure $J_t$ and a real 3-form $\ta_{2,t}'$ such that $\Oa_t' :=$ Re($\Oa_t$) + $i \ta_{2,t}'$ is a (3,0)-form w.r.t. $J_t$. Moreover, we have the 2-form $\oa_t'$, which is the rescaled (1,1)-part of $\oa_t$ w.r.t. $J_t$, and the associated metric $g_t$ on $M_t$. Following similar arguments to Proposition \ref{prop:closetogenuineisncy}, we conclude that $|g_t - g_{V_i}|_{g_{V_i}} = O(t^{-\lda_i(1 - \al)}) + O(t^{\al \nu}) = |g^{-1}_{t} - g^{-1}_{V}|_{g_{V_i}}$ . \\

\par Here is our main result \cite[Thm. 5.2]{Chan1} on the desingularization of compact \CY 3-folds $M_0$ with conical singularities in the simplest case $\lda_i < -3$, assuming the existence of singular \CY metrics on them: \\

\renewcommand{\thethm}{2.13} 
\begin{thm}
\label{thm:desingthmforlda<-3}
\quad Suppose $(M_0, J_0, \oa_0, \Oa_0)$ is a compact \CY 3-fold with conical singularities $x_i$ with rate $\nu > 0$ modelled on \CY cones $(V_i, J_{V_i} , \oa_{V_i}, \Oa_{V_i})$ for $i = 1, \dots, n$. Let $(Y_i, J_{Y_i} , \oa_{Y_i}, \Oa_{Y_i})$ be an AC \CY 3-fold with rate $\lda_i < -3$ modelled on the same \CY cone $V_i$. Define a family $(M_t, \oa_t, \Oa_t)$ of nonsingular compact nearly \CY 3-folds, with the associated metrics $g_t$ as above. \\[-0.2cm] 
\par Then $M_t$ admits a \CY structure $(\tilde{J}_t, \tilde{\oa}_t, \tilde{\Oa}_t)$ such that $\| \tilde{\oa}_t - \oa_t \|_{C^0} \leq K t^{\kappa}$ and\, $\| \tilde{\Oa}_t - \Oa_t \|_{C^0} \leq K t^{\kappa}$ for some $\kappa, K > 0$ and for sufficiently small $t$. The cohomology classes satisfy $[\textnormal{Re}(\Oa_t)]$ = $[\textnormal{Re}(\tilde{\Oa}_t)] \in H^3 (M_t, \BR)$\, and\, $[\oa_t]$ = $c_t\,[\tilde{\oa}_t] \in H^2 (M_t, \BR)$ for some $c_t > 0$. Here all norms are computed with respect to $g_t$. \\
\end{thm}

\par We conclude by applying the above result to the examples given in \S\ref{subsub:\CY cones\CY manifolds with conical}. Consider the situation in Example \ref{eg:ACforCm/G} and take $m = 3$. Then the crepant resolution $X$ of the \CY cone $\BC^3 / G$ is an AC \CY 3-fold with rate $-6$. Thus Theorem \ref{thm:desingthmforlda<-3} applies and we can desingularize any compact \CY 3-fold $M_0$ with conical singularities modelled on $\BC^3 / G$, or equivalently, any \CY 3-orbifold with isolated singularities, by the gluing process. Note that in general we have to assume the existence of singular \CY metrics on manifolds with conical singularities, but in the orbifold case, there is a result asserting the existence of \CY metrics: if $M$ is a compact \ka orbifold with $c_1 (M) = 0$, then there is a unique Ricci-flat \ka metric in each \ka class on $M$ (see for instance \cite[Thm. 6.5.6]{Joyce1}). \\

\par If we take $G = \mathbb{Z}_3$, a standard example of compact \CY 3-orbifold with isolated singularities is given in \cite[Example 6.6.3]{Joyce1}. Define a lattice $\Lda$ in $\BC^3$ by 
$$ \Lda = \BZ^3 \oplus \za \BZ^3 = \{ (a_1 + b_1 \za,\, a_2 + b_2 \za,\, a_3 + b_3 \za) : a_j, b_j \in \BZ \} $$ 
where $\za = -\frac{1}{2} + i\frac{\sqrt{3}}{2} = e^{2\pi i/3}$ denotes the cube root of unity. Let $T^6$ be the quotient $\BC^3/\Lda$, with a flat \CY structure $(J,\oa,\Oa)$. Write points on $T^6$ as $(z_1, z_2, z_3) + \Lda$ for $(z_1, z_2, z_3) \in \BC^3$. We can also regard $T^6$ as the product of three $T^2$'s where each $T^2$ is the quotient $\BC/(\BZ \oplus \za \BZ)$. \\

\par Define an action generated by $\za$ on $T^6$ by 
$$ \za \, : \, (z_1, z_2, z_3) + \Lda \longmapsto (\za z_1, \za z_2, \za z_3) + \Lda. $$
This $\za$-action is well-defined, as $\za \cdot \Lda = \Lda$. The group $\BZ_3 = \{ 1, \za, \za^2 \}$ is a finite group of automorphisms of $T^6$, preserving the flat \CY structure on it. Thus the toroidal orbifold $T^6/\BZ_3$ is a \CY 3-orbifold, which can also be expressed as $\BC^3/A$, where $A$ is the group generated rotations by $\za$ and translations in $\Lda$. Write points on $T^6/\BZ_3$ as $\BZ_3 \cdot (z_1, z_2, z_3) + \Lda$. \\

\par In each $T^2$ there are three fixed points of $\za$ located at 0, $\frac{1}{2} + \frac{i}{2\sqrt{3}}$, $\frac{i}{\sqrt{3}}$. The element $\za^2 = \za^{-1}$ clearly has the same fixed points. Altogether the orbifold $T^6/\BZ_3$ has then 27 isolated singularities. Note that $\frac{1}{2} + \frac{i}{2\sqrt{3}} = \frac{2i}{\sqrt{3}} - \za $, so we write the 27 fixed points on $T^6$ as 
$$ \left\{ (c_1, c_2, c_3) + \Lda : c_1, c_2, c_3 \in \big\{ 0, \frac{i}{\sqrt{3}}, \frac{2i}{\sqrt{3}} \big\} \right\}. $$ 
Now these singular points are locally modelled on the \CY cone $\BC^3/\BZ_3$, thus making the orbifold $T^6/\BZ_3$ a \CY 3-fold with conical singularities. Applying Theorem \ref{thm:desingthmforlda<-3}, we can desingularize $T^6/\BZ_3$ by gluing in AC \CY 3-folds $K_{\BC \BP^2}$ (with rate $-6$) at the singular points, obtaining a \CY desingularization of $T^6/\BZ_3$. \\

\par Now the \emph{Schlessinger Rigidity Theorem} \cite{Sch} states that if $G$ is a finite subgroup of GL($m$,$\BC$) and the singularities of $\BC^m/G$ are all of codimension at least three, then $\BC^m/G$ is rigid, i.e. it admits no nontrivial deformations. It can then be shown by using this rigidity theorem that if we desingularize a \CY 3-orbifold with isolated singularities modelled on $\BC^3 / \BZ_3$ by gluing, we shall obtain a crepant resolution of the original orbifold. \\

\par On the crepant resolution the existence of \CY metrics is guaranteed by Yau's solution to the Calabi conjecture \cite{Yau}. However, it does not provide a way to write down the \CY metrics explicitly, and so in general we do not know much about what the \CY metrics are like. But in the orbifold case, our result tells a bit more by giving a quantitative description of these \CY metrics, showing that these metrics locally look like the metrics obtained by gluing the orbifold metrics and the ALE metrics on the crepant resolution of $\BC^3 / G$.  \\


\begin{center}
\renewcommand{\thesubsection}{\normalfont{3}}
\subsection{\normalfont{\textsc{Special Lagrangian cones and their Lagrangian Neighbourhoods}}}
\label{sec:SpecialLagrangianConesAndTheirLagrangianNeighbourhoods}
\end{center}

\par We shall keep on using the definitions and notation established in \S\ref{sec:review}. Based on these settings, the corresponding kind of special Lagrangian submanifolds are defined. In this section we consider special Lagrangian cones in \CY cones. Again, we make the definitions in \S\ref{sec:SpecialLagrangianConesAndTheirLagrangianNeighbourhoods}, \ref{sec:SLMFoldsWithConicalSingularities} and \ref{sec:ACSLMFolds} in $m$ complex dimensions and later $m$ will be 3.  \\

\renewcommand{\thedefinition}{3.1}
\begin{definition}
\label{def:SLconesinCYcones}
\textnormal{\ Let $C$ be an SL $m$-fold, which is closed and nonsingular except at $0$, in the \CY cone $(V, J_{V}, \oa_{V},$ $ \Oa_{V})$. Then $C$ is an \emph{SL cone} in $V$ if $C' = C \setminus \{ 0\}$ can be written as $\Si \times (0, \infty)$ for some compact, nonsingular $(m-1)$-dimensional submanifold $\Si$ of the link $\Ga$ of $V$. Let $g_{C}$ be the restriction of the \CY cone metric $g_{V}$ to $C$.}  \\[-0.2cm]
\end{definition}

\par We call $\Si$ the link of the SL cone $C$. Here we can allow $\Si$ to be disconnected, or equivalently $C'$ disconnected, though we have assumed $\Ga$ is connected. \\

\par The fact that the SL cone $C' = \Si \times (0, \infty)$ is Lagrangian in the \CY cone $V' = \Ga \times (0, \infty)$ implies $\Si$ is a $(m-1)$-dimensional Legendrian submanifold in $\Ga$. On the \CY cone $V$, there is a complex dilation $\lda_{t, \ta} : V \longrightarrow V$, given by $\lda_{t, \ta} (\ga, r) = (\textnormal{exp}(\ta Z) (\ga),\, tr)$ for $\ta \in \BR$ and $t > 0$, where $Z = J_{V} X$ and $X = r\frac{\pa}{\pa r}$. When $\ta = 0$,  $\lda_{t, 0}$ gives a ``real dilation" $(\ga, r) \mapsto (\ga, tr)$, and the cone $C$ is therefore invariant under this real dilation, i.e. $C = \lda_{t,0} (C)$ for all $t > 0$. We note that $\lda_{t,0} (C)$ is still special Lagrangian in $V$ since 
$$ 0 \,=\, \lda_{t,0}^* ( \oa_{V}|_{C} ) \,=\, t^2 \oa_{V}|_{\lda_{t,0}(C)}, \quad \textnormal{and} \quad 0 \,=\, \lda_{t,0}^* ( \textnormal{Im}(\Oa_{V})|_{C} ) \,=\, t^m \textnormal{Im}(\Oa_{V})|_{\lda_{t,0}(C)}. $$ \\[-1cm]

\par Let $\iota : \Si \times (0, \infty) \longrightarrow \Ga \times (0, \infty)$ be the inclusion map given by $\iota (\si, r) = (\si, r)$. We identify $\Si \cong \Si \times \{ 1 \}$. Let $(\si, r) \in \Si \times (0, \infty)$. A 1-form on $\Si \times (0, \infty)$ at the point $(\si,r)$ can be expressed as $\eta + c\,dr$, where $\eta \in T^*_{\si} \Si$ and $c \in \BR$. Use $(\si, r, \eta, c)$ to denote a point in $T^*_{(\si,r)} (\Si \times (0, \infty))$. Identify $\Si \times (0, \infty)$ as the zero section $\{ (\si, r, \eta, c) : \eta = c = 0 \}$ of the cotangent bundle $T^* (\Si \times (0, \infty))$. Define a dilation action on $T^* (\Si \times (0, \infty))$ by
\begin{align*}
t \cdot (\si, r, \eta, c) = (\si, tr, t^2\,\eta, tc), \tag{3.1}
\label{eqn:taction}
\end{align*}
where $t \in \BR_+$. This $t$-action restricts to the usual dilation on the cone $\Si \times (0, \infty)$, and the pullback of the canonical symplectic form $\oa_{\textnormal{can}}$ on $T^* (\Si \times (0, \infty))$ by $t$ satisfies $t^*(\oa_{\textnormal{can}}) = t^2 \oa_{\textnormal{can}}$.  \\

\par The Lagrangian Neighbourhood Theorem \cite[Thm. 3.33]{McDuff} shows that any compact Lagrangian submanifold $N$ in a symplectic manifold looks locally like the zero section in $T^*N$. We are going to extend the Lagrangian Neighbourhood Theorem to special Lagrangian cones $C$ in the \CY cones $V$: \\

\renewcommand{\thethm}{3.2}
\begin{thm}
\label{thm:LNTforCi}
\ With the above notations, there exist an open tubular neighbourhood $U_{C}$ of the zero section $\Si \times (0, \infty)$ in $T^* (\Si \times (0, \infty))$, which is invariant under the $t$-action, and an embedding $\Psi_{C} : U_{C} \longrightarrow V' \cong \Ga \times (0, \infty)$ such that 
$$ \Psi_{C} |_{\Si_i \times (0, \infty)} = \iota_i, \quad \Psi_{C}^* (\oa_{V}) = \oa_{\textnormal{can}} \quad \textnormal{and} \quad \Psi_{C} \circ t = \lda_{t,0} \circ \Psi_{C} $$
for $t \in \BR_+$, where $t$ acts on $U_{C}$ as in (\ref{eqn:taction}), and $\lda_{t,0}$ is the dilation on the \CY cone $V$. \\
\end{thm}

\par Theorem \ref{thm:LNTforCi} can be proved by arguing in the same way as in the proof of \cite[Thm. 4.3]{Joycedesing1}, which applies \cite[Thm. 4.2]{Joycedesing1}, a version of a result of Weinstein on Lagrangian foliations. Here is the rough idea. In order to use Theorem 4.2 of \cite{Joycedesing1}, we need a smooth family $L_{(\si,r)}$ of noncompact Lagrangian submanifolds in $\Ga \times (0, \infty)$ containing the point $(\si,r)$ and transverse to $\Si \times (0, \infty)$ at $(\si,r)$, i.e. $T_{(\si,r)}L_{(\si,r)} \cap T_{(\si,r)}(\Si \times (0, \infty)) = \{ 0 \}$, with $L_{(\si,tr)} = \lda_{t,0}\,(L_{(\si,r)})$. Since $\Si \cong \Si \times \{ 1\}$ is compact, it is possible to choose such a family for $r=1$, i.e. we can choose noncompact Lagrangian submanifolds $L_{(\si,1)}$ for all $\si \in \Si$, transverse to $\Si \times (0, \infty)$ at $(\si, 1)$. Define $L_{(\si,r)} = \lda_{r,0}\,(L_{(\si,1)}) $, then the family $\{ L_{(\si,r)} : (\si, r) \in \Si_i \times (0, \infty) \}$ is what we need. \\

\par As we shall see, the Lagrangian neighbourhoods for $C$, together with the Lagrangian neighbourhoods for $N_0$ (SL 3-fold with conical singularities) and $L$ (AC SL 3-fold), are useful in the analytic result in Theorem \ref{thm:IIIthm5.3} as we will glue all these neighbourhoods together. \\

\begin{center}
\renewcommand{\thesubsection}{\normalfont{4}}
\subsection{\normalfont{\textsc{SL $m$-folds with conical singularities}}}
\label{sec:SLMFoldsWithConicalSingularities}
\end{center}



\par After defining SL cones in \CY cones, we now define SL $m$-folds $N_0$ with conical singularities. Essentially, $N_0$ is an SL $m$-fold in some \CY $m$-fold $M_0$ with conical singularities at $x_1, \dots, x_n$, and it is asymptotic at $x_1, \dots, x_n$ to the SL cones $C_1, \dots, C_n$ in the \CY cones $V_1, \dots, V_n$. The idea is to define $N_0$ as a graph of some 1-form $\za_i$ on $C_i$ near $x_i$ for $i= 1, \dots, n$. The fact that $N_0$ is Lagrangian implies $\za_i$ is a closed 1-form. Moreover, for $N_0$ to approach $C_i$ near $x_i$, $\za_i$ should decay at least like $O(r)$ which then implies $\za_i$ is in fact exact. As a result, we are able to express $N_0$ as a graph of some exact 1-form $da_i$ near $x_i$. Here is the precise definition of SL $m$-folds with conical singularities: \\


\renewcommand{\thedefinition}{4.1}
\begin{definition}
\label{def:SLcsinCYcs}
\textnormal{\ Let $N_0$ be a singular SL $m$-fold in a \CY $m$-fold $M_0$ with conical singularities at $x_1, \dots, x_n$ and with rate $\nu > 0$. Suppose $N_0$ is also singular at $x_1, \dots, x_n$ and has no other singularities. Then $N_0$ is an \emph{SL $m$-fold with conical singularities at} $x_1, \dots, x_n$ with rate $\mu \in (1,\,\nu + 1)$, modelled on SL cones $C_1, \dots, C_n$, if there exist an open neighbourhood $T_i \subset S_i$ of $x_i$ in $N_0$, and a smooth function $a_i$ on $\Si_i \times (0, \ep')$ for $i = 1, \dots, n$ and $\ep' < \ep$, satisfying}
\begin{align*}
| \nab^k a_i |_{g_{C_i}} =\, O(r^{\mu + 1 - k}) \quad \textnormal{as $r \rightarrow 0$ and for all $k \geq 0$,} \tag{4.1}
\label{eqn:sizeforai}
\end{align*}
\textnormal{computing $\nab$ and $|\cdot|_{g_{C_i}}$ using the cone metric $g_{C_i}$, such that} 
\begin{align*}
T_i \setminus \{ x_i \} =\, \Phi_i \circ \Psi_{C_i} ( \Ga ( d\,a_i ) )  \tag{4.2}
\label{eqn:Ti=graphdai}
\end{align*}
\textnormal{where $\Psi_{C_i} : U_{C_i} \longrightarrow V_i'$ is the embedding from the Lagrangian neighbourhood $U_{C_i}$ to the \CY cone $V_i'$, and $\Ga(d\,a_i)$ is the graph of the 1-form $d\,a_i$. From (\ref{eqn:sizeforai}), $d\,a_i$ has rate $O(r^{\mu})$ as $r \rightarrow 0$, which shows it is a small 1-form, and by making $\ep'$ smaller if necessary, the graph $\Ga(d\,a_i)$ of $d\,a_i$ lies in $U_{C_i}$.} \\ 
\end{definition}

\par We require the rate $\mu$ to be greater than 1 to ensure $N_0$ approaches the cone $C_i$ near $x_i$. For the upper bound, we choose $\mu$ to be less than $\nu + 1$ so that whether $N_0$ is an SL $m$-fold with conical singularities with rate $\mu$ is independent of the choice of the diffeomorphisms or coordinates $\Phi_i$ amongst equivalent coordinates. Recall that two coordinates $\Phi_i$ and $\Phi_i'$ are equivalent if and only if 
(\ref{eqn:equivalentCYcoord}) holds. 
From (\ref{eqn:Ti=graphdai}), $ T_i \setminus \{ x_i \} $ is written as $\Phi_i \circ \Psi_{C_i} ( \Ga ( d\,a_i ) )$. If $a_i'$ is another smooth function such that $ T_i \setminus \{ x_i \} = \Phi_i' \circ \Psi_{C_i} ( \Ga ( d\,a_i' ) )$ under the coordinate $\Phi_i'$, then 
\begin{align*}
| \nab^{k+1} a_i' |_{g_{C_i}} &= | \nab^k ( \Phi_i - \Phi_i' ) |_{g_{V_i}}  + | \nab^{k+1} a_i |_{g_{C_i}} \ = \ O(r^{\nu + 1 - k}) + O(r^{\mu - k}) \quad \textnormal{from (\ref{eqn:sizeforai})} \\
&= O(r^{\mu - k}) \quad \textnormal{provided $\mu \leq \nu + 1$ } 
\end{align*}
for $k \geq 0$. Integration then gives $|a_i'|_{g_{C_i}} = O(r^{\mu + 1})$. Hence (\ref{eqn:sizeforai}) for $a_i$ is equivalent to (\ref{eqn:sizeforai}) for $a_i'$, and the definition of SL $m$-folds with conical singularities with rate $\mu$ is therefore independent of the choice of $\Phi_i$. \\

\par Now we give a result on the construction of a Lagrangian neighbourhood for $N_0$, compatible with the Lagrangian neighbourhoods $U_{C_i}$ of $C_i$ and with the embeddings $\Psi_{C_i} : U_{C_i} \longrightarrow V_i' \cong \Ga_i \times (0, \infty)$ for $i=1, \dots, n$ in Theorem \ref{thm:LNTforCi}. The proof is similar to that of \cite[Thm. 4.6]{Joycedesing1}. \\

\renewcommand{\thethm}{4.2}
\begin{thm}
\label{thm:LNTforN0}
\ With the above notations, there exists an open tubular neighbourhood $U_{N_0}$ of the zero section $N_0$ in $T^* N_0$ such that 
$$ d \Phi_i |_{U_{C_i} \cap T^* (\Si_i \times (0, \ep'))} \, ( U_{C_i} \cap T^* (\Si_i \times (0, \ep')) ) \ = \ T^* (T_i \setminus \{ x_i \} ) \cap U_{N_0} \quad \text{for $i = 1, \dots, n$,} $$ 
and there exists an embedding $\Psi_{N_0} : U_{N_0} \longrightarrow M_0$ with $\Psi_{N_0}|_{N_0} = Id$, $\Psi_{N_0}^* (\oa_0) = \oa_{T^* N_0}$, where Id is the identity map on $N_0$ and $\oa_{T^* N_0}$ the canonical symplectic structure on $T^* N_0$, such that 
$$ \Psi_{N_0} \, \circ \, d \Phi_i |_{U_{C_i} \cap T^* (\Si_i \times (0, \ep'))} (\si, r, \eta, c) \ = \ \Phi_i \circ \Psi_{C_i} (\si, r, \eta + da^1_i (\si,r), c + da^2_i (\si,r)) $$ 
for $i = 1, \dots, n$, $(\si, r, \eta, c) \in U_{C_i} \cap T^* (\Si_i \times (0, \ep'))$, and for $da_i (\si,r) = da^1_i (\si,r) + da^2_i (\si,r) dr $ with $da^1_i (\si,r) \in T^*_{\si} \Si_i$ and $da^2_i (\si,r) \in \BR$.  \\
\end{thm}

\par Let us focus on the 
zero section $\{ \eta = c = 0\} \cong \Si_i \times (0, \ep')$. The first part of the theorem gives 
$$ d \Phi_i |_{U_{C_i} \cap T^* (\Si_i \times (0, \ep'))} (\si, r, 0, 0)  \ = \ T_i \setminus \{ x_i \}. $$ 
This is consistent with the second part of the theorem, as from the left hand side we have
$$ \Psi_{N_0} \, \circ \, d \Phi_i |_{U_{C_i} \cap T^* (\Si_i \times (0, \ep'))} (\si, r, 0, 0) \ = \ \Psi_{N_0} (T_i \setminus \{ x_i \}) \ = \ T_i \setminus \{ x_i \}, $$
using $\Psi_{N_0}|_{N_0}$ = Id, and the right hand side gives 
$$ \Phi_i \circ \Psi_{C_i} (\si, r, da^1_i (\si,r), da^2_i (\si,r)) \ = \ \Phi_i \circ \Psi_{C_i} (\Ga(da_i)) \ = \ T_i \setminus \{ x_i \} \qquad \textnormal{by (\ref{eqn:Ti=graphdai})}. $$  \\[-0.3cm]

\begin{center}
\renewcommand{\thesubsection}{\normalfont{5}}
\subsection{\normalfont{\textsc{AC SL $m$-folds}}}
\label{sec:ACSLMFolds}
\end{center}


\par We now study AC SL $m$-folds $L$ in AC \CY $m$-folds $Y$. \\

\centerline{5.1\ \textsc{\ Definition and the Lagrangian Neighbourhood Theorem}}
\bigskip



\par Similar to the conical singularities case, we want to define $L$ as a graph of some closed 1-form $\chi$ near infinity. The condition for $L$ to converge to the SL cone $C$ at infinity is that $\chi$ has size $O(r^{\kp})$ for large $r$ and $\kp < 1$. \\


\renewcommand{\thedefinition}{5.1}
\begin{definition}
\label{def:ACSLinACCY}
\textnormal{\ Let $L$ be a nonsingular SL $m$-fold in an AC \CY $m$-fold $Y$ with rate $\lda < -m$. Then $L$ is an \emph{AC SL m-fold} with rate $\kp \in (\lda + 1, 1)$, modelled on SL cones $C$ if there exist a compact subset $H \subset L$ and a smooth closed 1-form $\chi$ on $\Si \times (R', \infty)$ for $R' > R$, satisfying}
\begin{align*}
| \nab^k \chi |_{g_{C}} =\, O(r^{\kp - k}) \quad \textnormal{as $r \rightarrow \infty$ and for all $k \geq 0$,} \tag{5.1}
\label{eqn:sizeforchii}
\end{align*}
\textnormal{computing $\nab$ and $|\cdot|_{g_{C}}$ using the cone metric $g_{C}$, such that} 
\begin{align*}
L \setminus H =\, \Ups \circ \Psi_{C} ( \Ga ( \chi ) ).  \tag{5.2}
\label{eqn:Li-Hi=graphchii}
\end{align*}
\textnormal{Equation (\ref{eqn:sizeforchii}) implies $\chi$ has rate $O(r^{\kp})$ as $r \rightarrow \infty$, and by making $R'$ larger if necessary, the graph $\Ga(\chi)$ of $\chi$ lies in $U_{C}$.} \\ 
\end{definition}

\par Analogous to the upper bound for the rate $\mu$ in the conical singularities case, we require $\kp > \lda + 1$ so that the definition of AC SL $m$-folds with rate $\kp$ does not depend on the choice of the coordinate $\Ups$. \\

\par Here we work with closed 1-forms with rate $\kp < 1$ in the definition. However, assuming $\chi$ to be closed is not enough for our purpose, as we hope to express $L$ as a graph of an exact 1-form near infinity. The reason for that is if the 1-form $\chi$ was not exact, we shall come across global topological obstructions (cf. the obstructed case in \cite{Joycedesing4}, or the $\lda_i = -3$ case for \CY 3-folds in \cite{Chan2}) which we do not want to deal with. \\

\par Note that if we assume $\kp < -1$, then $\chi$ is automatically exact, and it can be written as $d\,b$ for some smooth function $b$ on $\Si \times (R', \infty)$. We can construct $b$ by integration: Write $\chi (\si,r) = \chi^1 (\si,r) + \chi^2 (\si,r)\,dr$ where $\chi^1 (\si,r) \in T^*_{\si} \Si$ and $\chi^2 (\si,r) \in \BR$. Similar to the argument we used before in defining the 1-form in \cite[Thm. 4.9]{Chan1}, we define $b (\si,r) = \int^{\infty}_{r} \chi^2 (\si,s) ds $, which is well-defined if $\kp < -1$, and it satisfies $b(\si,r) = O(r^{\kp + 1})$.  Then $\chi (\si,r) = d\,b (\si,r)$. It follows that assuming $\kp < -1$ will suit our purpose, and in fact we shall see in \S\ref{sec:JoyceSDesingularizationTheory} that we need to assume $\kp < -3/2$ to apply Theorem \ref{thm:IIIthm5.3} (for the case $m=3$). 
Thus from now on, we adjust the rate $\kp_i$ of AC SL $m$-folds to be less than $-1$, so that $\kp \in (\lda + 1, -1)$, and equations (\ref{eqn:sizeforchii}) and (\ref{eqn:Li-Hi=graphchii}) in the definition become respectively 
\begin{align*}
| \nab^k b |_{g_{C}} =\, O(r^{\kp + 1 - k}) \quad \textnormal{as $r \rightarrow \infty$ and for all $k \geq 0$,} \tag{5.3}
\label{eqn:sizeforbi}
\end{align*}
and
\begin{align*}
L \setminus H =\, \Ups \circ \Psi_{C} ( \Ga ( d\,b ) ).  \tag{5.4}
\label{eqn:Li-Hi=graphdbi}
\end{align*} \\[-0.5cm]

\par In the last part of this section, we give the Lagrangian Neighbourhood Theorem for AC SL $m$-folds $L$ (compare to \cite[Thm. 7.5]{Joycedesing1}), which is an analogue of Theorem \ref{thm:LNTforN0}: \\

\renewcommand{\thethm}{5.2}
\begin{thm}
\label{thm:LNTforLi}
\ With the above notations, there exists an open tubular neighbourhood $U_{L}$ of the zero section $L$ in $T^* L$ such that 
$$ d \Ups |_{U_{C} \cap T^* (\Si \times (R', \infty))} \, ( U_{C} \cap T^* (\Si \times (R', \infty)) ) \ = \ T^* (L \setminus H) \cap U_{L}, $$ 
and there exists an embedding $\Psi_{L} : U_{L} \longrightarrow Y$ with $\Psi_{L}|_{L} = Id$, $\Psi_{L}^* (\oa_{Y}) = \oa_{T^* L}$, where Id is the identity map on $L$ and $\oa_{T^* L}$ the canonical symplectic structure on $T^* L$, such that 
$$ \Psi_{L} \, \circ \, d \Ups |_{U_{C} \cap T^* (\Si \times (R', \infty))} (\si, r, \eta, c) \ = \ \Ups \circ \Psi_{C} (\si, r, \eta + db^1(\si,r), c + db^2(\si,r)) $$ 
for $(\si, r, \eta, c) \in U_{C} \cap T^* (\Si \times (R', \infty))$ and for $db(\si,r) = db^1(\si,r) + db^2(\si,r) dr$ with $db^1(\si,r) \in T^*_{\si} \Si$ and $db^2(\si,r) \in \BR$. \\
\end{thm}

\begin{center}
\renewcommand{\thesubsubsection}{\normalfont{5.2}}
\subsubsection{\normalfont{\textsc{Some examples of AC SL $m$-folds}}}
\label{sec:SomeExamplesOfACSLMFolds}
\end{center}



\par Some examples of AC SL $m$-folds in the complex Euclidean space $\BC^m$ can be found in Harvey and Lawson \cite[III.3]{HarLaw} and Joyce \cite[\S6.4]{Joycedesing5}. Here we give some examples of AC SL $m$-folds in the crepant resolution of the \CY cone $\BC^m / \BZ_m$, or equivalently, the total space of the canonical line bundle $K_{\BC \BP^{m-1}}$ over $\BC \BP^{m-1}$ endowed with Calabi's metric, as described in Example \ref{eg:ACforCm/G}. \\ 


\renewcommand{\theeg}{5.3}
\begin{eg}
\label{eg:fpsinKCPm-1} 
\textnormal{\ Our first example of an AC SL $m$-fold will be given by the fixed point set of an antiholomorphic isometric involution on the AC \CY $m$-fold $K_{\BC \BP^{m-1}}$. Consider the usual complex conjugation $\si_0 : \BC^m \longrightarrow \BC^m$ having $\BR^m$ as its fixed point set. Since $\si_0 \circ \ga \circ \si_{0}^{-1} = \ga^{-1} $ for any $\ga \in \BZ_m$, $\si_0$ induces an involution on $\BC^m / \BZ_m$, and then lifts to $\si$ on the crepant resolution of $\BC^m / \BZ_m$, or $K_{\BC \BP^{m-1}}$. This map $\si$ is actually an antiholomorphic isometric involution on $K_{\BC \BP^{m-1}}$: since $\si$ is induced by the complex conjugation, it satisfies $\si^2 = \textnormal{Id}$ and $\si^* (J) = -J$ where $J$ is the complex structure on $K_{\BC \BP^{m-1}}$. We know that $\si^* (r) = r$, where $r$ is the radius function on $\BC^m / \BZ_m \setminus \{0\} $, and hence $\si^* (f) = f$ as the \ka potential $f$ defined by Calabi is a function of $r^2$. It follows that $\si^* (g) = g$ and hence $\si^* (\oa) = -\oa$. Furthermore, $\si^* (\Oa) = \bar{\Oa} $ on $K_{\BC \BP^{m-1}}$ since the holomorphic volume form on $K_{\BC \BP^{m-1}}$ is just the same as that on $\BC^m$. Consequently, $\si : K_{\BC \BP^{m-1}} \longrightarrow K_{\BC \BP^{m-1}}$ is an antiholomorphic isometric involution on $K_{\BC \BP^{m-1}}$, and the fixed point set of $\si$ is an SL $m$-fold of $K_{\BC \BP^{m-1}}$.} \\
\end{eg}

\par Next we investigate the fixed point set of $\si$ in $K_{\BC \BP^{m-1}}$, which is just the union of the fixed points set $L_0$ of $\si_0 : \BC^m / \BZ_m \longrightarrow \BC^m / \BZ_m $ with the origin removed and the fixed points set of $\si|_{\BC \BP^{m-1}} : \BC \BP^{m-1} \longrightarrow \BC \BP^{m-1}$ in $\BC \BP^{m-1}$. Now a point $\BZ_m \cdot (z_1,\, \dots , z_m)$ in $\BC^m / \BZ_m$ is fixed by $\si_0$ whenever $(\bar{z}_1,\, \dots , \bar{z}_m) = (z_1 \, e^{2\pi i k/m},\, \dots , z_m \, e^{2\pi ik/m})$ for some $0 \leq k \leq m-1$. It follows that $z_j = r_j \, e^{\pi ik/m}$ where $r_j \in \BR$, $j = 1, \dots , m$ and hence the fixed point set $L_0$ in $\BC^m / \BZ_m$ is given by
$$L_0 = \{  \BZ_m \cdot (r_1 \, e^{\pi ik/m},\, \dots , r_m \, e^{\pi ik/m} ) : r_j \in \BR, \ k = 0,1, \dots , m-1 \}. $$
Observe that $\BZ_m \cdot (r_1 \, e^{\pi ik/m},\, \dots , r_m \, e^{\pi ik/m} )$ is equal to $\BZ_m \cdot (r_1,\, \dots , r_m )$ for $k$ even and $\BZ_m \cdot (r_1 \, e^{\pi i/m},\, \dots , r_m \, e^{\pi i/m} )$ for $k$ odd. Thus $L_0$ has two components, 
\begin{align*}
L_0 = (\BZ_m \cdot \BR^m) / \BZ_m  \, \cup \, e^{\pi i/m}\,((\BZ_m \cdot \BR^m)/ \BZ_m). \tag{5.8}
\label{eqn:fixpointofsi0} 
\end{align*}
\begin{itemize}
\item[(i)] When $m$ is even, then $-1 \in \BZ_m$, and $\BZ_2$ is the subgroup of $\BZ_m$ fixing $\BR^m$. In this case, $L_0$ is topologically a union of two copies of $\BR^m / \BZ_2$, i.e. two cones on $\BR \BP^{m-1}$ meeting at 0; 
\item[(ii)] When $m$ is odd, then $\frac{m-1}{2} \in \BZ$, and hence
\begin{align*}
e^{\pi i/m}\,\BZ_m \cdot (r_1,\, \dots , r_m) &= e^{\pi i/m}\, e^{\frac{2\pi i}{m} (\frac{m-1}{2})}\,\BZ_m \cdot (r_1,\, \dots , r_m) \\
&= e^{\pi i}\,\BZ_m \cdot (r_1,\, \dots , r_m) \\
&= \BZ_m \cdot (-r_1,\, \dots , -r_m).
\end{align*}
It follows from (\ref{eqn:fixpointofsi0}) that $L_0 = (\BZ_m \cdot \BR^m) / \BZ_m $, which is topologically a copy of $\BR^m $, or equivalently, a cone on $S^{m-1}$. 
\end{itemize}
Together with the fixed point set of $\si|_{\BC \BP^{m-1}} : \BC \BP^{m-1} \longrightarrow \BC \BP^{m-1}$, which is just $\BR \BP^{m-1}$, this yields the fixed point set $L$ of $\si$, i.e. an AC SL $m$-fold in $K_{\BC \BP^{m-1}}$: 
\begin{itemize}
\item[(i)] When $m$ is even, $L$ is homeomorphic to $\BR \BP^{m-1} \times \BR$; 
\item[(ii)] When $m$ is odd, $L$ is homeomorphic to $S^{m-1} \times (0, \infty) \ \cup \ \BR \BP^{m-1}$. We may regard $L$ as the quotient $(S^{m-1} \times \BR) / \BZ_2$, where $\BZ_2$ acts freely on $S^{m-1} \times \BR$. 
\end{itemize}
In both cases, $L$ is the canonical line bundle $K_{\BR \BP^{m-1}}$ over $\BR \BP^{m-1}$. In fact, $K_{\BR \BP^{m-1}}$ is trivial (nontrivial) if $m$ is even (odd), as $\BR \BP^{m-1}$ is oriented (non-oriented).  \\

\par We are particularly interested in the case $m=3$. From the above analysis we obtain an AC SL 3-fold $K_{\BR \BP^2}$ as the fixed point set of some antiholomorphic isometric involution in the AC \CY 3-fold $K_{\BC \BP^2}$. Observe that $K_{\BR \BP^2}$ admits a double cover which is diffeomorphic to $S^2 \times \BR$. We remark that this AC   
SL 3-fold $K_{\BR \BP^2}$ has rate $``\kp = -\infty"$, since there will be coordinates in which $K_{\BR \BP^2}$ is a cone. As we have discussed before in the definition of AC SL $m$-folds, we require $1 + \textnormal{rate for } K_{\BC \BP^2} < \kp < -1$. Thus we could say the rate for $K_{\BR \BP^2}$ is any $\kp \in (-5, -1)$, as $K_{\BC \BP^2}$ has rate $-6$. \\

\renewcommand{\theeg}{5.4}
\begin{eg}
\label{eg:SO3invarinKCPm-1} 
\textnormal{\ This example constructs AC SL 3-folds in $K_{\BC \BP^2}$ invariant under the standard SO(3)-action. From the calculation of the moment map of the U($m$)-action on $\BC^m$ w.r.t. Calabi's metric in \cite[\S5.1.1]{Chan}, it can be seen that the moment map $\mu : K_{\BC \BP^2} \longrightarrow \mathfrak{so}(3)^* \cong \BR^3 $
of the SO(3)-action on $K_{\BC \BP^2}$ is given by
$$ \mu = f'(r^2)\,( \textnormal{Im}(z_1 \bar{z}_2 ), \textnormal{Im}(z_2 \bar{z}_3 ),\, \textnormal{Im}(z_3 \bar{z}_1 ) ). $$}
\end{eg}
Again, $f$ denotes the \ka potential defined by Calabi. Since $Z(\mathfrak{so}(3)^*) = \{0\}$, it follows that any SO(3)-invariant SL 3-fold in $K_{\BC \BP^2}$ must lie in the level set $\mu^{-1} (0)$. Using the same construction as in \cite[Example 5.2]{JoyceSLsym} or \cite[Example 4.7]{Chan}, we obtain a family of SL 3-folds $L_c$ in $K_{\BC \BP^2}$ which has the same form as in \cite[Example 5.2]{JoyceSLsym} or \cite[Example 4.7]{Chan}, since the level sets of both moment maps coincide. Thus for $c \ne 0$, $L_c$ is an AC SL 3-fold in $K_{\BC \BP^2}$ with Calabi's metric and is diffeomorphic to $S^2 \times \BR$. As the cone $(\BZ_3 \cdot \BR^3) / \BZ_3$ is identified with the cone $e^{i\pi/3} \cdot (\BZ_3 \cdot \BR^3) / \BZ_3$ in $\BC^3 / \BZ_3$, $L_c$ converges to two copies of $(\BZ_3 \cdot \BR^3) / \BZ_3$ in $\BC^3 / \BZ_3$ from \cite[Thm. A]{Haskins} or \cite[Thm. 6.4]{JoyceSLsym}. We have seen that the rate for the SO(3)-invariant SL 3-folds in $\BC^3$ is $-2$, thus $L_c$ also has rate $-2$. With the discussion in Example \ref{eg:fpsinKCPm-1} the fixed point set $K_{\BR \BP^2}$ admits a double cover diffeomorphic to $S^2 \times \BR$, which is just one possible $L_c$, and hence we have obtained a family of deformations of a double cover of Example \ref{eg:fpsinKCPm-1} for $m=3$. One can generalize this example to higher dimensions and obtain SO($m$)-invariant AC SL $m$-folds in $K_{\BC \BP^{m-1}}$. \\

\renewcommand{\theeg}{5.5}
\begin{eg}
\label{eg:Tm-1invarinKCPm-1} 
\textnormal{\ Finally we briefly describe an example of a $T^{m-1}$-invariant AC SL $m$-fold in $K_{\BC \BP^{m-1}}$ taken from \cite[\S5.1.1]{Chan}. The idea of the construction is to use the method of moment maps, similar to \cite[Prop. 4.4]{Chan}. The $T^{m-1}$-action on $K_{\BC \BP^{m-1}}$ is just the $T^{m-1}$-action on $K_{\BC \BP^{m-1}} \setminus \BC \BP^{m-1} \cong \BC^m / \BZ_m \setminus \{ 0 \}$ given by 
$$ (e^{i\ta_1}, \dots , e^{i\ta_m}) \cdot (z_1, \dots , z_m) = (e^{i\ta_1} z_1, \dots , e^{i\ta_m} z_m) $$}
\end{eg} 
where $\ta_1 + \cdots + \ta_m = 0$, and the $T^{m-1}$-action on $\BC \BP^{m-1}$. According to the calculation in \cite[\S5.1.1]{Chan}, the moment map $\mu$ for the $T^{m-1}$-action on $K_{\BC \BP^{m-1}}$ is given by: 
\begin{align*}
\mu (z_1 , \dots , z_m) = f'(r^2)\,\left(|z_1|^2 - \frac{1}{m}\,r^2 ,\, \dots , |z_m|^2 - \frac{1}{m}\,r^2 \right) 
\end{align*}
on the bit $K_{\BC \BP^{m-1}} \setminus \BC \BP^{m-1} \cong \BC^m / \BZ_m \setminus \{ 0 \}$. Here $f$ is the \ka potential defined by Calabi. On $\BC \BP^{m-1}$, we have 
\begin{align*}
\mu([w_1, \dots, w_m]) = \left( |w_1|^2 - \frac{1}{m}\,,\, \dots , |w_m|^2 - \frac{1}{m} \right) 
\end{align*}
where $w_1, \dots , w_m$ are normalized such that $|w_1|^2 + \cdots + |w_m|^2 = 1$. Up to a constant factor, this is just the same as the moment map of $T^{m-1}$ acting on $\BC \BP^{m-1}$ with the Fubini-Study metric. \\

\par Following \cite[\S5.1.1]{Chan} we have the following: Let $c = (c_1, \dots , c_m, c')$ where $c_1,\, \dots , c_m, c' \in \BR$ and $c_1 + \cdots + c_m = 0$. Define 
\begin{align*}
L_c = \Big\{ &(z_1,\, \dots , z_m ) \in \BC^m \setminus \{ 0 \} : f'(r^2)\,\big( |z_j|^2 - \frac{1}{m}\,r^2 \big) = c_j\,,\ \textnormal{for } j = 1, \dots , m,  \\
&\quad \textnormal{and} \ \textnormal{Re}(z_1 \cdots z_m) = c' \ \textnormal{if $m$ is even, }\ \textnormal{Im}(z_1 \cdots z_m) = c' \ \textnormal{if $m$ is odd} \Big\}.
\end{align*} 
Then $L_c$ is a $T^{m-1}$-invariant SL $m$-fold in $\BC^m \setminus \{ 0 \}$ with respect to Calabi's metric. Since $L_c$ is invariant under the $\BZ_m$-action on $\BC^m \setminus \{ 0 \}$, it follows that the quotient $L_c/\BZ_m$ is a $T^{m-1}$-invariant SL $m$-fold in $\BC^m / \BZ_m \setminus \{ 0 \}$, and converges to the cone 
\begin{align*}
\Big\{ \BZ_m \cdot (z_1,\, \dots , z_m ) : |z_1|^2 = \cdots = |z_m|^2 = \frac{1}{m}\,r^2,\ \textnormal{and} \ &\textnormal{Re}(z_1 \cdots z_m) = c' \ \textnormal{if $m$ is even, } \\
&\textnormal{Im}(z_1 \cdots z_m) = c' \ \textnormal{if $m$ is odd} \Big\}
\end{align*} 
in $\BC^m / \BZ_m$. \\

\par Now suppose $c' \ne 0$, then $L_c/\BZ_m$ will not intersect $\BC \BP^{m-1}$ in $K_{\BC \BP^{m-1}}$. 
Consequently, $L_c/\BZ_m$ is a $T^{m-1}$-invariant SL $m$-fold in $K_{\BC \BP^{m-1}}$ for any $c_1,\, \dots , c_m, c' \in \BR$ with $c_1 + \cdots + c_m = 0$ and $c' \ne 0$. \\

\par For the case $c' = 0$, if $(c_1,\, \dots , c_m)$ does not lie in the convex hull of the $m$ points $\mu(p_1),\, \dots, \mu(p_m)$ in $\BR^{m-1}$, then $L_c/\BZ_m$ is a $T^{m-1}$-invariant SL $m$-fold in $K_{\BC \BP^{m-1}}$. If the point $(c_1,\, \dots , c_m)$ is on a $k$-dimensional face of the simplex for some $k = 0, 1,\,\dots , m-1$, then we have to include a $k$-torus $T^k$ in $\BC \BP^{m-1}$, and hence in this case we obtain a $T^{m-1}$-invariant SL $m$-fold $L_c/\BZ_m \cup T^k$ in $K_{\BC \BP^{m-1}}$. \\

\begin{center}
\renewcommand{\thesubsection}{\normalfont{6}}
\subsection{\normalfont{\textsc{Joyce's desingularization theory}}}
\label{sec:JoyceSDesingularizationTheory}
\end{center}


\par Joyce has developed a comprehensive desingularization theory of SL $m$-folds with isolated conical singularities in \CY $m$-folds (and more generally in \emph{almost} \CY $m$-folds). His approach is to glue in appropriate AC SL $m$-folds in $\BC^m$ which are asymptotic to some SL cones, thus obtaining a 1-parameter family of compact nonsingular Lagrangian $m$-folds. Then he proves using analysis that the Lagrangian $m$-folds which are close to being special Lagrangian can actually be deformed to SL $m$-folds in the \CY $m$-fold. The whole programme on SL $m$-folds with isolated conical singularities is given in the series of his papers \cite{Joycedesing1, Joycedesing2, Joycedesing3, Joycedesing4, Joycedesing5}. \\

\par We shall now fix $m=3$ to fit into our situation for desingularizing SL 3-folds in \CY 3-folds.  \\

\par Before stating Joyce's analytic result, we need to establish the necessary notations, which can be found in Definition 5.2 of \cite{Joycedesing3}. However we shall only consider the following particular case of his definition: \\

\par Let $(M, J, \oa, \Oa)$ be a \CY 3-fold, with \CY metric $g$. Since we only deal with \CY manifolds, not the more general class of almost \CY manifolds, we can take the smooth function $\psi$ in Definition 5.2 of \cite{Joycedesing3} to be 1 on $M$, so that condition (ii) of Theorem 5.3 in \cite{Joycedesing3} becomes trivial in our case. \\
  
\par Suppose $N \subset M$ is a Lagrangian 3-fold. Restricting the \CY metric $g$ on $N$, we obtain a metric $h = g|_N$, with volume form $dV$. The holomorphic (3,0)-form $\Oa$ restricts to a 3-form $\Oa|_N$ on $N$. As $N$ is Lagrangian, we can write
$$ \Oa|_N = e^{i\ta}\,dV $$    
for some phase function $\ta$ on $N$, which equals zero if and only if $N$ is special Lagrangian (with phase 1). Suppose [Im($\Oa$)$|_N$] = 0 in $H^3 (N, \BR)$, or equivalently,
$$ \int_N \textnormal{Im}(\Oa) = \int_N \textnormal{sin}\ta\,dV = 0. $$ 
Clearly, this is a necessary condition for $N$ to be special Lagrangian in $M$. Note that this can be satisfied by choosing the phase of $\Oa$ appropriately, which means this is actually a fairly mild restriction. \\

\par In (iii) of Theorem 5.3 in \cite{Joycedesing3}, we need a 3-form $\ba$ and a connection $\acute{\nab}$ on $T^* N$ ($\hat{\nab}$ in the original notation). For $r > 0$, define $\mathcal{B}_r \subset T^* N$ to be $\mathcal{B}_r = \{ \al \in C^{\infty} (T^* N) : |\al|_{h^{-1}} < r \}$, where $|\cdot|_{h^{-1}}$ is computed using the metric $h^{-1}$. The Levi-Civita connection $\nab$ of $h$ on $TN$ induces a splitting $T\mathcal{B}_r = H \oplus V$, where $H \cong TN$ and $V \cong T^*N$ are the horizontal and vertical subbundles of $T^*N$. Define $\acute{h}$ on $\mathcal{B}_r$ such that $H$ and $V$ are orthogonal w.r.t. $\acute{h}$, and $\acute{h}|_{H} = \pi^* (h)$, $\acute{h}|_{V} = \pi^* (h^{-1})$ where $\pi : \mathcal{B}_r \rightarrow N$ denotes the natural projection. Let $\acute{\nab}$ be the connection on $T\mathcal{B}_r$ given by lift of $\nab$ on $N$ in the horizontal directions $H$, and by partial differentiation in the vertical directions $V$. Since $N$ is Lagrangian in $M$, for small $r > 0$, the Lagrangian Neighbourhood Theorem gives an embedding $\Psi : \mathcal{B}_r \rightarrow M$ such that $\Psi^* (\oa) = \oa_{\textnormal{can}}$ and $\Psi|_N = \textnormal{Id} $, where $\oa_{\textnormal{can}}$ is the natural symplectic structure on $\mathcal{B}_r \subset T^* N$. Finally, we define a 3-form $\ba$ on $\mathcal{B}_r$ by $\ba = \Psi^* (\textnormal{Im}\,(\Oa))$, the pullback of the imaginary part of the holomorphic (3,0)-form $\Oa$ on $M$ by the Lagrangian embedding.    \\

\par The finite dimensional vector space $W$ in Definition 5.2 of \cite{Joycedesing3} has to do with the number of connected components of $N_0 \setminus \{ x_1, \dots, x_n \}$ where $N_0$ is a SL 3-fold with conical singularities at $x_1, \dots, x_n$. We will make the nonsingular Lagrangian 3-fold $N_t$ (which is our $N$ defined above) by gluing AC SL 3-folds $L_1, \dots, L_n$ into $N_0$. If $N_0 \setminus \{ x_1, \dots, x_n \}$ is not connected, then each of the $L_i$'s is connected but may contain more than one end, so that $N_t$ consists of several components of $N_0$ joined by ``small necks'' from $L_i$'s. The vector space $W$ will then be a space of functions which is approximately constant on each component of $N_0 \setminus \{ x_1, \dots, x_n \}$ and changes on small necks. The dimension dim$W$ will be the number of connected component of $N_0 \setminus \{ x_1, \dots, x_n \}$. For the sake of simplicity, we only study the case when $N_0 \setminus \{ x_1, \dots, x_n \}$ is connected, which means we can take $W = \left< 1 \right> $, the space of constant functions. As a result, condition (vii) of Theorem 5.3 of \cite{Joycedesing3} is trivial and we can drop it entirely. The natural projection $\pi_W: L^2(N) \rightarrow W$ is now given by
$$ \pi_W (v) = \textnormal{vol}(N)^{-1} \int_N v\,dV $$ 
for $W = \left< 1 \right> $, and we have $\pi_W (v) = 0 \Longleftrightarrow \int_N v\,dV = 0$. It follows that the last inequality of Theorem 5.3 (i) in \cite{Joycedesing3} holds automatically, as we have assumed $\int_N \textnormal{sin}\ta\,dV = 0$. Moreover, we can replace $\pi_W (v) = 0$ by $\int_N v\,dV = 0$ in (vi) of the theorem to leave $W$ out of our definition. \\

\par We are now ready to state the following analytic existence result for SL 3-folds, adapted from \cite[Thm. 5.3]{Joycedesing3}: \\

\renewcommand{\thethm}{6.1}
\begin{thm}
\label{thm:IIIthm5.3}
\ Let $\kp' > 1$ and $A_1, \dots, A_6 > 0$. Then there exist $\ep, K > 0$ depending only on $\kp', A_1, \dots, A_6$ such that the following holds. \\[-0.3cm]
\par Refer to the notation in $\S$\ref{sec:JoyceSDesingularizationTheory}. \emph{Suppose} $0 < t \leq \ep$ \emph{and} $r = A_1 t$, \emph{and} 
\begin{itemize}
	\item[\textnormal{(i)}] $\| \textnormal{sin}\,\ta \|_{L^{6/5}} \leq A_2 t^{\kp' + 3/2}$, \quad $\| \textnormal{sin}\,\ta \|_{C^0} \leq A_2 t^{\kp' - 1}$ \quad and \quad $\| d\,\textnormal{sin}\,\ta \|_{L^{6}} \leq A_2 t^{\kp' - 3/2}$.	
	\item[\textnormal{(ii)}] $\| \acute{\nab}^k \ba \|_{C^0} \leq A_3 t^{-k}$ \ \text{for} $k$ = $0,\, 1,\, 2$ and $3$. 
	\item[\textnormal{(iii)}] The injectivity radius $\da (h)$ satisfies $\da (h) \geq A_4 t$.
	\item[\textnormal{(iv)}] The Riemann curvature $R(h)$ satisfies $\| R(h) \|_{C^0} \leq A_5 t^{-2}$.
	\item[\textnormal{(v)}] If $v \in L^2_1 (N)$ with $\int_N v\, dV = 0$, then $v \in L^6(N)$, and\, $\| v \|_{L^6} \leq A_6 \| dv \|_{L^2}$. 
\end{itemize}
Here all norms are computed using the metric $h$ on $N$ in \textnormal{(i), (iv)} and \textnormal{(v)}, and the metric $\acute{h}$ on $\mathcal{B}_{A_1 t}$ in \textnormal{(ii)}. Then there exists $f \in C^{\infty} (N)$ with $\int_N f\, dV = 0$, such that $\| df \|_{C^0} \leq Kt^{\kp'} < A_1 t$ and $\hat{N} = \Psi_* (\Ga (df))$ is an immersed special Lagrangian 3-fold in $(M, J, \oa, \Oa)$. \\ 
\end{thm}

\par For a small 1-form $\al \in C^{\infty} (T^*N)$, $\Psi_* (\Ga (\al))$ is a Lagrangian submanifold in $M$ if $\al$ is closed, and it is a special Lagrangian submanifold if $\al$ also satisfies $d^* (\textnormal{cos}\,\ta\, \al) = \textnormal{sin}\,\ta + Q(\al)$ where $Q$ is a smooth map with $Q(\al) = O( |\al|^2 + |\nab \al|^2)$ (see \cite[Lemma 5.7]{Joycedesing3}). Now if $f \in C^{\infty} (N)$ is a small function in $C^1(N)$, then $df$ is a small 1-form, and 
$\Psi_* (\Ga (df))$ is special Lagrangian if and only if $d^* (\textnormal{cos}\,\ta\, df) = \textnormal{sin}\,\ta + Q(df)$. Thus the function $f$ in the theorem is basically the solution to the above nonlinear elliptic equation. Joyce \cite[\S5.5]{Joycedesing3} solved this equation by constructing inductively a sequence $(f_k)_{k=0}^{\infty}$ in $C^{\infty}(N)$ with $f_0 = 0$ and $\int_N f_k\, dV = 0$ satisfying $d^* (\textnormal{cos}\,\ta df_k) = \textnormal{sin}\,\ta + Q(df_{k-1})$ for $k \geq 1$. Then he showed that this sequence converges in some Sobolev space to $f$ which satisfies the nonlinear elliptic equation and is smooth by elliptic regularity. \\

\begin{center}
\renewcommand{\thesubsection}{\normalfont{7}}
\subsection{\normalfont{\textsc{Construction of $N_t$}}}
\label{sec:ConstructionOfNT}
\end{center}


\par In this section we are going to define a 1-parameter family of compact nonsingular Lagrangian 3-folds $N_t$ in the nearly \CY 3-folds $M_t$ for small $t > 0$. Recall that we desingularize the \CY 3-fold $M_0$ with conical singularities by first applying a homothety to each AC \CY 3-fold $Y_i$, and then gluing it into $M_0$ at $x_i$ for $i=1, \dots, n$ to make the nonsingular $M_t$'s. For the special Lagrangians inside these Calabi--Yau's, we then desingularize $N_0$ (inside $(M_0, J_0, \oa_0, \Oa_0)$) by gluing $L_i$ (inside $(Y_i, J_{Y_i}, t^2\oa_{Y_i}, t^3\Oa_{Y_i})$) into $N_0$ at $x_i$ for each $i$. Note that after applying the homothety to $Y_i$, equations (\ref{eqn:sizeforbi}) and (\ref{eqn:Li-Hi=graphdbi}) now become 
\begin{align*}
| \nab^k b_i (\si, t^{-1}r) |_{g_{C_i}} =\, O(t^{-\kp_i - 1} r^{\kp_i + 1 - k}) \quad \textnormal{as $r \rightarrow \infty$ and for all $k \geq 0$,} \tag{7.1}
\label{eqn:sizeforbi(t-1r)}
\end{align*}
and 
\begin{align*}
L_i \setminus H_i =\, \Ups_{t,i} \circ \Psi_{C_i} ( \Ga ( t^2 d\,b_i (\si, t^{-1} r) ) )  \tag{7.2}
\label{eqn:Li-Hi=graphdbi(t-1r)}
\end{align*} 
where the diffeomorphism $ \Ups_{t,i} : \Ga_i \times (tR, \infty) \longrightarrow Y_i \setminus K_i $ is given by $\Ups_{t,i} (\ga, r) = \Ups_i (\ga, t^{-1} r) $. \\

\par Now for $i=1, \dots, n$, $\al \in (0,1)$ and small enough $t > 0$ with $tR < tR' < t^{\al} < 2t^{\al} < \ep' < \ep$, define a smooth function $u_{t,i}$ on $\Si_i \times (tR', \ep')$ by
\begin{align*}
u_{t,i}(\si, r) = F(t^{-\al}r)\,a_i(\si, r) + t^2 ( 1 - F(t^{-\al}r) )\, b_i(\si, t^{-1}r), \tag{7.3}
\label{eqn:defforuti}
\end{align*} 
where $F$ is the smooth increasing function we used before in defining $\Oa_t$ in \CY desingularizations. Thus we have $F(t^{-\al}r) = 1$ when $2t^{\al} \leq r < \ep'$, in which case $u_{t,i} (\si, r) = a_i (\si, r)$, and $F(t^{-\al}r) = 0$ when $tR' \leq r < t^{\al}$, in which case $u_{t,i} (\si, r) = t^2 b_i (\si, t^{-1}r)$. \\

\par Here we are going to use the same $\al \in (0,1)$ as in the previous \CY 3-fold desingularizations, and we shall require $\al$ to satisfy certain inequalities for this SL desingularization as well. \\

\par Define $ N_t$ to be the union of $(N_0 \setminus \bigcup_{i=1}^{n} T_i)$, $\bigcup_{i=1}^{n} \Psi_{C_i} (\Ga(du_{t,i}))$, and $\bigcup_{i=1}^{n} H_i $. Basically, we construct $N_t$ in the way that for $r > \ep'$, $N_t$ is just $(N_0 \setminus \bigcup_{i=1}^{n} T_i) \subset N_0$. For $tR' < r < \ep'$, $N_t$ is diffeomorphic to the union of the graphs $\Ga(du_{t,i})$ of the 1-forms $du_{t,i}$ which in fact interpolate between the graphs $\Ga(da_i)$ of $da_i$, i.e. part of $N_0$, for $2t^{\al} \leq r < \ep'$ and the graphs $\Ga(t^2 db_i)$ of $t^2 db_i$, i.e. part of $L_i$, for $tR' \leq r < t^{\al}$. Finally for $ r < tR'$, $N_t$ is the union of $H_i \subset L_i$. \\

\par Under our construction the boundary of each $H_i$, which is just the link $\Si_i$, joins smoothly onto $\bigcup_{i=1}^{n} \Psi_{C_i} (\Ga(du_{t,i}))$ at the $\Si_i \times \{ tR' \} $ end, and the boundary of $(N_0 \setminus \bigcup_{i=1}^{n} T_i)$, which is the disjoint union of the $\Si_i$, joins smoothly onto $\bigcup_{i=1}^{n} \Psi_{C_i} (\Ga(du_{t,i}))$ at the $\Si_i \times \{ \ep' \} $ end. Thus $N_t$ is a compact smooth manifold without boundary. More importantly, $N_t$ is in fact a Lagrangian submanifold: \\

\renewcommand{\theprop}{7.1}
\begin{prop}
\label{prop:NtLag}
\ $N_t$ is a Lagrangian 3-fold in the nearly \CY 3-fold $(M_t, \oa_t, \Oa_t)$ for sufficiently small $t > 0$. \\[-0.3cm]
\end{prop}
\pf \, We shall look at the symplectic form $\oa_t$ restricts to different regions of $N_t$. Since $N_t$ coincides with $N_0$, which is Lagrangian in $M_0$, in the component $N_0 \setminus \bigcup_{i=1}^{n} T_i$, and $\oa_t$ equals $\oa_0$ on this part, we see that $\oa_t \equiv 0$ on $N_0 \setminus \bigcup_{i=1}^{n} T_i$. In the same way, as $N_t$ coincides with the union of $L_i$ in the component $\bigcup_{i=1}^{n} H_i $, and $\oa_t$ is now $t^2 \oa_{Y_i}$ on each $L_i$, thus we have $\oa_t \equiv 0$ on $\bigcup_{i=1}^{n} H_i $, as $L_i$ is Lagrangian in $Y_i$. For the middle part, $N_t$ is given by $\bigcup_{i=1}^{n} \Psi_{C_i} (\Ga(du_{t,i}))$. As $\oa_t$ is equal to $\oa_{V_i}$ on each $\Psi_{C_i} (\Ga(du_{t,i}))$, and $\Psi_{C_i}^* (\oa_{V_i}) = \oa_{\textnormal{can}}$ from Theorem 5.2, where $\oa_{\textnormal{can}}$ is the canonical symplectic form on $T^* (\Si_i \times (0, \infty))$, we get $\oa_t \equiv 0$ on $\bigcup_{i=1}^{n} \Psi_{C_i} (\Ga(du_{t,i}))$ as $\Ga(du_{t,i})$ is the graph of a closed 1-form $du_{t,i}$ on $\Si_i \times (tR', \ep')$. It follows that $\oa_t|_{N_t} = 0$, and $N_t$ is then Lagrangian in $M_t$. \hfill $\Box$ \\

\par Now we deform the underlying nearly \CY structure on $M_t$ to a genuine \CY structure $(\tilde{J}_t, \tilde{\oa}_t, \tilde{\Oa}_t)$ for small $t$ by applying Theorem \ref{thm:desingthmforlda<-3}. As shown in the theorem, we have the relation $[\oa_t] = c_t\,[\tilde{\oa}_t] \in H^2 (M_t, \BR)$ for some $c_t > 0$ between the cohomology classes of the \ka forms. Thus $\oa_t$ and $c_t\,\tilde{\oa}_t$ are in the same cohomology class. Using Moser's type argument there is a diffeomorphism $\psi_t : M_t \longrightarrow M_t $ on $M_t$ satisfying $\psi_t^* (c_t\,\tilde{\oa}_t) = \oa_t$. Write $\hat{\oa}_t = \psi_t^* (\tilde{\oa}_t)$, $\hat{J}_t = \psi_t^* (\tilde{J}_t)$, $\hat{g}_t = \psi_t^* (\tilde{g}_t)$ and $\hat{\Oa}_t = \psi_t^* (\tilde{\Oa}_t)$ under the new coordinates. The fact that $\oa_t$ and $c_t\,\tilde{\oa}_t$ are close for small $t>0$ means that the diffeomorphism $\psi_t$ is close to identity, which then implies that the complex 3-forms are also close under the new coordinates, i.e. $\hat{\Oa}_t \approx \Oa_t$. We shall evaluate this difference in the next section. \\

\par As we have arranged $\oa_t = c_t\,\hat{\oa}_t$ by applying a diffeomorphism $\psi_t$, it follows that $\hat{\oa}_t|_{N_t} = 0$ since $c_t > 0$, and so we obtain: \\

\renewcommand{\theprop}{7.2}
\begin{prop}
\label{prop:NtLagingenuine}
\ $N_t$ is a Lagrangian 3-fold in the \CY 3-fold $(M_t, \hat{J}_t, \hat{\oa}_t, \hat{\Oa}_t)$ for sufficiently small $t >0$. \\
\end{prop}

\begin{center}
\renewcommand{\thesubsection}{\normalfont{8}}
\subsection{\normalfont{\textsc{Estimates of Im($\hat{\Oa}_t$)$|_{N_t}$}}}
\label{sec:EstimatesOfImTildeOaTNT}
\end{center}


\par We have constructed a family of compact, nonsingular Lagrangian 3-folds $N_t$ by gluing $L_i$ into $N_0$ at each $x_i$, and our next step is to apply Theorem \ref{thm:IIIthm5.3} to deform $N_t$ to a special Lagrangian 3-fold $\hat{N}_t$ in the \CY 3-fold $(M_t, \hat{J}_t, \hat{\oa}_t, \hat{\Oa}_t)$ for small enough $t >0$. This leads us to consider the estimates of various norms of Im($\hat{\Oa}_t$)$|_{N_t}$ for (i) of Theorem \ref{thm:IIIthm5.3}. \\

\par Let $h_t, \tilde{h}_t$ and $\hat{h}_t$ be the restrictions of $g_t, \tilde{g}_t$ and $\hat{g}_t$ to $N_t$ respectively. In view of Theorem \ref{thm:desingthmforlda<-3}, the metrics $g_t, \tilde{g}_t$, and hence $h_t, \tilde{h}_t$, are uniformly equivalent in $t$, so norms of any tensor on $N_t$ measuring with respect to $h_t, \tilde{h}_t$ only differ by a bounded factor independent of $t$. We shall see later the uniform equivalence between $\hat{h}_t$ and the other two.   \\

\par Here is the basic estimate, computing with respect to $\hat{h}_t$:
\begin{align*}
\big| \textnormal{Im}(\hat{\Oa}_t)|_{N_t} \big|_{\hat{h}_t} \, &\leq \, \big| (\hat{\Oa}_t - \tilde{\Oa}_t) |_{N_t} \big|_{\hat{h}_t} \, + \, \big| (\tilde{\Oa}_t - \Oa_t) |_{N_t} \big|_{\hat{h}_t} \, + \, \big| \textnormal{Im}(\Oa_t)|_{N_t} \big|_{\hat{h}_t}. \tag{8.1}
\label{eqn:basicestimate}
\end{align*}
We hope to arrange for this error to be small enough that we can deform $N_t$ to an SL 3-fold by using the analytic result in Theorem \ref{thm:IIIthm5.3}. The first term $(\hat{\Oa}_t - \tilde{\Oa}_t) |_{N_t}$ in the right side of (\ref{eqn:basicestimate}) is basically the error coming from changing coordinates on $M_t$, which can be estimated by considering the diffeomorphism $\psi_t$ from Moser's argument. The second term $(\tilde{\Oa}_t - \Oa_t) |_{N_t} $ is the error arising from deforming the nearly \CY structure $(\oa_t, \Oa_t)$ to the genuine \CY structure $(\tilde{J}_t, \tilde{\oa}_t, \tilde{\Oa}_t)$ on $M_t$. We already have the $C^0$-estimates from Theorem \ref{thm:desingthmforlda<-3}, but for part (i) of Theorem \ref{thm:IIIthm5.3} to hold we need to improve and get a better control of this term. We shall devote most of the section to achieving this. For the final term $\textnormal{Im}(\Oa_t)|_{N_t}$, we can estimate it by restricting $\Oa_t$ to different regions of $N_t$. \\

\par Let us first evaluate the last term Im($\Oa_t$) on various components of $N_t$. 
From (\ref{eqn:defofOat}) we have
\begin{align*}
\Oa_t = 
\begin{cases}
\Oa_0$ \ \ on\ \ $Q_t \setminus \big[ (\bigcup_{i=1}^n P_{t,i}) \cap Q_t \big], \notag \\
\Oa_{V_i} + d\big[ F(t^{-\al}r) A_i(\ga, r) + t^3 (1 - F(t^{-\al}r)) B_i(\ga, t^{-1}r)\big]$ on $P_{t,i} \cap Q_t, \textnormal{for } i = 1, \dots, n, \\ 
t^3 \Oa_{Y_i}$\ \ on\ \ $P_{t,i} \setminus (P_{t,i} \cap Q_t)\ \textnormal{for } i = 1, \dots, n.   
\end{cases} \\[-0.5cm]  
\end{align*}

\par On $Q_t \setminus [(\bigcup_{i=1}^{n} P_{t,i}) \cap Q_t ]$, $\Oa_t$ is given by $\Oa_0$, and $N_t$ is the union of $N_0 \setminus \bigcup_{i=1}^{n} T_i$ and $\bigcup_{i=1}^{n} \Psi_{C_i} (\Ga(da_i))$, which is a part of $N_0$. Thus Im($\Oa_t$) = 0 on this region of $N_t$, as $N_0$ is special Lagrangian in $M_0$. Similarly, on $P_{t,i} \setminus (P_{t,i} \cap Q_t) $ for each $i$, $\Oa_t$ is given by $t^3 \Oa_{Y_i}$, and $N_t$ is the union of $\bigcup_{i=1}^{n} H_i $ and $\bigcup_{i=1}^{n} \Psi_{C_i} (\Ga(t^2 db_i))$ which lies in $L_i$. It follows that Im($\Oa_t$) = 0 on this region of $N_t$, as $L_i$ is special Lagrangian in $Y_i$. For the annuli region, we have 
\begin{align*}
\textnormal{Im}(\Oa_t)|_{N_t} = \textnormal{Im} \big( \Oa_{V_i} + d\big[ F(t^{-\al}r) A_i(\ga, r) + t^3 (1 - 
F(t^{-\al}r)) B_i(\ga, t^{-1}r)\big] \big) \big|_{\Psi_{C_i} (\Ga(du_{t,i}))}. \tag{8.2} 
\label{eqn:annuliOatonNt}
\end{align*}  
Consider the term Im($\Oa_{V_i}$)$|_{\Psi_{C_i} (\Ga(du_{t,i}))}$. Regard $C_i$ as the zero section in $T^* (\Si_i \times (0, \infty))$. Then the difference between Im($\Oa_{V_i}$)$|_{\Psi_{C_i} (\Ga(du_{t,i}))}$ and Im($\Oa_{V_i}$)$|_{C_i}$ is given by 
$$ O( \, |\acute{\nab} \Oa_{V_i} |_{g_{V_i}} \cdot \, |du_{t,i}|_{g_{C_i}} \,) \, + \, O( \, |\Oa_{V_i} |_{g_{V_i}} \cdot \, |\nab du_{t,i}|_{g_{C_i}} \,), $$
where $\acute{\nab}$ denotes a connection on $T^* (\Si_i \times (0, \infty))$ (see $\S$\ref{sec:JoyceSDesingularizationTheory} for the construction of $\acute{\nab}$ on $T^* N$), and $\nab$ the connection on $C_i$ computed using the metric $g_{C_i}$. Roughly speaking, the first term is coming from moving base points, whereas the second term from changing tangent spaces. Now we have $|\acute{\nab} \Oa_{V_i} |_{g_{V_i}} = O(t^{-\al})$ on the annulus, and $|\Oa_{V_i} |_{g_{V_i}}$ is a constant. Together with the fact that $C_i$ is special Lagrangian in $V_i$, i.e. Im($\Oa_{V_i}$)$|_{C_i}$ = 0, we obtain the size for Im($\Oa_{V_i}$)$|_{\Psi_{C_i} (\Ga(du_{t,i}))}$: 
\begin{align*}
\big| \, \textnormal{Im}(\Oa_{V_i})|_{\Psi_{C_i} (\Ga(du_{t,i}))} \, \big|_{g_{C_i}} \, = \ O(\,t^{-\al} |du_{t,i}|_{g_{C_i}}\,) \, + \, O(\,|\nab du_{t,i}|_{g_{C_i}}\,) \tag{8.3}
\label{eqn:sizeofOaVionduti}
\end{align*}
for $r \in (t^{\al}, 2t^{\al})$. Using (\ref{eqn:sizeforai}) and (\ref{eqn:sizeforbi(t-1r)}), and the definition of $u_{t,i}$ in (\ref{eqn:defforuti}) we get 
\begin{align*}
| du_{t,i} |_{g_{C_i}} &= \, O(t^{\mu \al}) + O(t^{1 - \kp_i(1-\al)}), \ \ \textnormal{and} \\[-0.3cm]
\tag{8.4} \label{eqn:sizedutiandnabduti}  \\[-0.3cm]
| \nab du_{t,i} |_{g_{C_i}} &= \, O(t^{(\mu-1) \al}) + O(t^{(1-\kp_i)(1-\al)}) \ \ \textnormal{for } r \in (t^{\al}, 2t^{\al}).
\end{align*}  
Putting (\ref{eqn:sizedutiandnabduti}) into (\ref{eqn:sizeofOaVionduti}), and using the estimates for $A_i$ and $B_i$ from (\ref{eqn:asympbehavefornabkA}) and (\ref{eqn:asympbehavefornabkB}), we compute the size for (\ref{eqn:annuliOatonNt}):
\begin{align*}
\big| \, \textnormal{Im}(\Oa_t)|_{N_t} \, \big|_{g_{C_i}} &= \, O(t^{(\mu-1) \al}) + O(t^{(1-\kp_i)(1-\al)}) + O(t^{\al \nu}) + O(t^{-\lda_i(1-\al)}) \\
&= \, O(t^{(\mu-1) \al}) + O(t^{(1-\kp_i)(1-\al)})  \ \ \textnormal{for } r \in (t^{\al}, 2t^{\al}). \tag{8.5}
\label{eqn:3termestimate}
\end{align*}  
The term $O(t^{\al \nu})$ is absorbed into $O(t^{(\mu-1) \al})$ as we have chosen $\mu < \nu + 1$ in the definition of SL 3-folds with conical singularities, and similarly the term $O(t^{-\lda_i(1-\al)})$ is absorbed into $O(t^{(1-\kp_i)(1-\al)})$ as $\kp_i > \lda_i + 1$ in the definition of AC SL 3-folds. \\

\par Summing up all these, and using the uniform equivalence between metrics $g_{C_i}$ and $h_t$ (follows from that between $g_{V_i}$ and $g_t$), we see that \\

\renewcommand{\theprop}{8.1}
\begin{prop}
\label{prop:thirdtermofbasicestimate}
\ In the situation above, the error term $\textnormal{Im}(\Oa_t)|_{N_t}$ satisfies
\begin{align*}
\big| \, \textnormal{Im}(\Oa_t)|_{N_t} \, \big|_{h_t} = 
\begin{cases}
0$ \ \ on\ \ $N_t \cap \big( Q_t \setminus [(\bigcup_{i=1}^{n} P_{t,i}) \cap Q_t ] \big),  \\
O(t^{(\mu-1) \al}) + O(t^{(1-\kp_i)(1-\al)})$ \ \ on\ \ $N_t \cap \big( P_{t,i} \cap Q_t \big), $\ \ for\ \ $ i = 1, \dots, n,  \\ 
0$ \ \ on\ \ $N_t \cap \big( P_{t,i} \setminus (P_{t,i} \cap Q_t) \big), $\ \ for\ \ $ i = 1, \dots, n.   
\end{cases} \\[-0.5cm]  
\end{align*} 
\end{prop}
 



\par Next we estimate the term $(\tilde{\Oa}_t - \Oa_t)|_{N_t}$ in (\ref{eqn:basicestimate}), which comes from deforming the nearly \CY structure to the genuine \CY structure on $M_t$. From Theorem \ref{thm:desingthmforlda<-3}, we have the $C^0$-estimates for $\tilde{\Oa}_t - \Oa_t$ on the whole $M_t$ given by $\| \tilde{\Oa}_t - \Oa_t \|_{C^0} = O(t^{\kp})$ for some $\kp > 0$. This term then contributes $O(t^{\kp})$ to the basic estimate (\ref{eqn:basicestimate}), which in turn contributes $O(t^{\kp})$ to different norms of sin\,$\ta$ in Theorem \ref{thm:IIIthm5.3}. But for the $L^{6/5}$-estimate to hold, one needs $\kp \geq \, \kp' + 3/2 > \, 5/2$ as $\kp' > 1$. Since we have no a priori control of $\kp > 0$, this will put a strong restriction on $\kp$, and in turn on the rates $\nu$ and $\lda_i$ as well. \\

\par To resolve this problem we are going to improve the global $L^2$-estimate for $\tilde{\Oa}_t - \Oa_t$ to $C^0$-estimates locally by applying modified versions of \cite[Thm. 3.7 \& 3.8]{Chan1} (see also Theorems G1 and G2 in \cite[\S11.6]{Joyce1} for 7 dimensions). Note that if we look back on the construction of $\tilde{\Oa}_t$ in Theorem \ref{thm:existencegenuineCY}, the size of $\tilde{\Oa}_t - \Oa_t$ is of the same order as the size of $d\eta_t = \tilde{\varphi}_t - \varphi_t$, in other words, the error introduced when deforming the nearly \CY structure to the genuine \CY structure on the 6-fold $M_t$ is essentially the same as that introduced when deforming the $G_2$-structure to the torsion-free $G_2$-structure on the 7-fold $S^1 \times M_t$. It suggests that in order to get a better control of the $C^0$-norm of $\tilde{\Oa}_t - \Oa_t$ on $M_t$, one could consider improving the $C^0$-norm of $d\eta_t$ on $S^1 \times M_t$. \\

\par In \cite[\S11.6]{Joyce1} Joyce proved an existence result for torsion-free $G_2$-structures by constructing the 2-form $\eta$ upon solving a nonlinear elliptic p.d.e. (equation (11.33) of \cite{Joyce1}) in $\eta$. His method of solving the p.d.e. is to inductively construct sequences of 2-forms $\{ \eta_j \}_{j=0}^{\infty}$ and functions $\{ f_j \}_{j=0}^{\infty}$ with $\eta_0 = f_0 = 0$, and then he showed that these sequences converge in some Sobolev spaces to limits $\eta$ and $f$ which satisfy the p.d.e. The $C^0$-estimate of $d\eta$ is derived from the $C^0$-estimates of the sequence elements $d\eta_j$. So to improve $\| d\eta \|_{C^0}$, we need to improve Theorems G1 and G2 in \cite{Joyce1}, or more appropriately, the 6-dimensional version of them, i.e. Theorems 3.7 and 3.8 in \cite{Chan1}, in our situation. \\

\par Here is the modified version of \cite[Thm. 3.7]{Chan1}:

\renewcommand{\thethm}{8.2} 
\begin{thm}
\label{thm:6foldestimatesmodified}
\ Let $D_2, D_3 > 0$ be constants, and suppose $(M,g)$ is a complete Riemannian 6-fold with a continuous function $\rho$ having the following properties:
\begin{itemize}
	\item[(1)] the injectivity radius of geodesics $\da(g)_x$ of $(M,g)$ starting at $x$ satisfies $\da(g)_x \geq D_2 \rho(x)$, 
	\item[(2)] the Riemann curvature $R(g)$ satisfies $|R(g)|_g \leq D_3 \rho^{-2}$ on $M$, and 
	\item[(3)] for all $x \in M$, we have $1/2 \, \rho(x) \leq \rho \leq 2\rho(x)$ on balls $B_{D_2 \rho(x)} (x)$ of radius $D_2 \rho(x)$ about $x$. 
\end{itemize}
Then there exist $K_1, K_2 > 0$ depending only on $D_2$ and $D_3$, such that if \,$\chi \in L^{12}_{1} (\Lda^3 T^* M) \cap L^2 (\Lda^3 T^* M)$ then 
\begin{align*}
\| \rho^{7/2} \, \nab \chi \|_{L^{12}} &\leq K_1 \,( \, \| \rho^{7/2} \, d\chi \|_{L^{12}} + \| \rho^{7/2} \, d^* \chi \|_{L^{12}} + \| \chi \|_{L^2} \, ) \\[0.2cm]
and \quad \| \rho^3 \, \chi \|_{C^0} &\leq K_2 \,( \, \| \rho^{7/2} \nab \chi \|_{L^{12}} + \| \chi \|_{L^2} \, ). \\
\end{align*} 
\end{thm}

\par We shall call $\rho$ a \emph{local injectivity radius function} on $M$. Condition (3) ensures that $\rho$ does not change quickly, and we may treat it as constant on $B_{D_2 \rho(x)} (x)$. Moreover, (1) and (3) imply $\da(g)_y \geq 1/2 \, D_2 \rho(x)$ for all $y \in B_{D_2 \rho(x)} (x)$, whereas (2) and (3) give $|R(g)|_g \leq 4 D_3 \rho(x)^{-2}$ on $B_{D_2 \rho(x)} (x)$. The right hand sides of these inequalities are just constants, so that we get control of injectivity radius and Riemann curvature on balls about $x$ with radius at most $D_2 \rho(x)$, which can then be compared with Euclidean balls.  \\

\par As in \cite[Thm. 3.7]{Chan1}, we can prove Theorem \ref{thm:6foldestimatesmodified} using the same method of proof as that of Theorem G1 in \cite[p. 298]{Joyce1}, but we now use balls of radius $L \rho(x)$, where $0 < L < D_2$, about $x$ instead of $Lt$ in the proof of Theorem G1. Note that it is important to have the constants $D_2, D_3$ in the theorem independent of $t$, so that $K_1$ and $K_2$ are independent of $t$ as well. \\

\par Next we give the improved result for Theorem \cite[Thm. 3.8]{Chan1}: 

\renewcommand{\thethm}{8.3} 
\begin{thm}
\label{thm:6foldsequencemodified}
\ Let $\kappa > 0$ and $D_1, D_2, D_3, D_4, K_1, K_2 > 0$ be constants. Then there exist constants $\ep \in (0,1]$, $K_3$ and $K > 0$ such that whenever $0 < t \leq \ep$, the following is true. \\[-0.4cm]
\par Let $M$ be a compact 6-fold, with metric $g_M$ and a local injectivity radius function $\rho$ satisfying (1), (2) and (3) in Theorem \ref{thm:6foldestimatesmodified}. Suppose $\rho$ also satisfies $\rho \geq D_4 t > 0$ on $M$. Let $(\vi, g)$ be an $S^1$-invariant $G_2$-structure on $\SM$ with $d\vi = 0$. Suppose $\psi$ is an $S^1$-invariant smooth 3-form on the 7-fold $\SM$ with $d^* \psi = d^* \vi$, and
\begin{itemize}
	\item[\textnormal{(i)}] $\| \psi \|_{L^2} \leq D_1 t^{3+\kp}$, $\| \rho^3 \, \psi \|_{C^0} \leq D_1 t^{3+\kp}$\ \,and\ \,$\| \rho^{7/2} \, d^* \psi \|_{L^{12}} \leq D_1 t^{3+\kp}$,
	\item[\textnormal{(ii)}] if $\chi \in L^{12}_{1} (\Lda^3 T^* (\SM))$ is $S^1$-invariant, then $\| \rho^{7/2} \, \nab \chi \|_{L^{12}} \leq K_1 \,(\| \rho^{7/2} \, d\chi \|_{L^{12}} + \| \rho^{7/2} \, d^* \chi \|_{L^{12}} + \| \chi \|_{L^2} )$, 
	\item[\textnormal{(iii)}] if $\chi \in L^{12}_{1} (\Lda^3 T^* (\SM))$ is $S^1$-invariant, then $\| \rho^3 \, \chi \|_{C^0} \leq K_2 \,(  \| \rho^{7/2} \nab \chi \|_{L^{12}} + \| \chi \|_{L^2})$.
	\end{itemize} 
Let $\ep_1$ and $F$ be as in Definition 10.3.3 and Proposition 10.3.5 of Joyce \cite{Joyce1} respectively. Denote by $\pi_1$ the orthogonal projection from $\Lda^3 T^* (\SM)$ to the 1-dimensional component of the decomposition into irreducible representation of $G_2$. Then there exist sequences $\{\eta_j\}^{\infty}_{j=0}$ in $L^{12}_{2} (\Lda^2 T^*(\SM))$ and $\{f_j\}^{\infty}_{j=0}$ in $L^{12}_{1} (\SM)$ with $\eta_j, f_j$ being all $S^1$-invariant and $\eta_0 = f_0 = 0$, satisfying the equations 
$$ (dd^* + d^* d)\eta_j = d^* \psi + d^* (f_{j-1} \psi) + \ast \, dF(d\eta_{j-1}) \ and \ f_j \vi = \frac{7}{3} \pi_1 (d\eta_j) $$ 
for each $j >0$, and the inequalities
\begin{align*}
&(a)\ \| d\eta_j \|_{L^2} \leq 2 D_1 t^{3 + \kp},  &&(d)\ \| d\eta_j - d\eta_{j-1} \|_{L^2} \leq 2 D_1 2^{-j} t^{3 + \kappa},  \\ 
&(b)\ \| \rho^{7/2} \, \nab d\eta_j \|_{L^{12}} \leq K_3 t^{3 + \kp}, &&(e)\ \| \rho^{7/2} \, \nab (d\eta_j - d\eta_{j-1}) \|_{L^{12}} \leq K_3 2^{-j} t^{3 + \kp}, \\
&(c)\ \| \rho^3 \, d\eta_j \|_{C^0} \leq K t^{3+\kp}  \qquad and \ &&(f)\ \| \rho^3 \, (d\eta_j - d\eta_{j-1}) \|_{C^0} \leq K 2^{-j} t^{3+\kp}. 
\end{align*}
Here $\nab$ and $\| \cdot \|$ are computed using $g$ on $\SM$. \\
\end{thm} 

\par We can prove Theorem \ref{thm:6foldsequencemodified} by applying Theorem \ref{thm:6foldestimatesmodified} in place of \cite[Thm. 3.7]{Chan1}, and then follow the same arguments in the proof of Theorem G2 \cite[p. 299]{Joyce1}. The only extra issue here is that we need a lower bound for $\rho$: $\rho \geq D_4 t > 0$ on $M$. The inequality in part (c) implies $ | d\eta_j |_{g} \leq K t^{3+\kp} \rho^{-3} \leq K D_{4}^{-3} t^{\kp}$ if $\rho \geq D_4 t$. Thus assuming $\rho \geq D_4 t$ on $M$ gives $| d\eta_j |_{g} \leq \ep_1$ if $t$ is sufficiently small, where $\ep_1$ is the small positive constant defined in Definition 10.3.3 in \cite{Joyce1}. The lower bound for $\rho$ therefore ensures $d\eta_j$ is small enough for each $j$ which is needed to apply Proposition 10.3.5 in \cite{Joyce1}. We also make a remark here about the difference between (a), (d) and (b), (c), (e), (f) in the theorem: (a) and (d) are global estimates on the whole manifold, as we get $\| d\eta_j \|^{2}_{L^2}$ and $\| d\eta_j - d\eta_{j-1} \|^{2}_{L^2}$ from the elliptic equation in the theorem by integration by parts. On the other hand, (b), (c), (e) and (f) are local estimates on small balls, and so we are allowed to insert powers of $\rho$.   \\

\par Let us now return to our \CY 3-fold $M_t$. Define a function $\rho_t$ on $M_t$ by 

\begin{align*}
\rho_t = 
\begin{cases} \\[-1cm]
\ep \ \ &\textnormal{on}\ \ M_0 \setminus \bigcup_{i=1}^{n} S_i ,  \\
r \ \ &\textnormal{on}\ \ \Ga_i \times (tR, \ep), $\ \ for\ \ $ i = 1, \dots, n,  \tag{8.6} \label{eqn:defofrhotonMt} \\ 
tR \ \ &\textnormal{on}\ \ K_i \subset Y_i, $\ \ for\ \ $ i = 1, \dots, n.   
\end{cases} \\[-0.5cm]  
\end{align*} 
We claim that $\rho_t$ is a local injectivity radius function on $M_t$. To see properties (1) and (2), recall that the way we construct $g_t$ on $M_t$ is, on $M_0 \setminus \bigcup_{i=1}^{n} S_i$ it is equal to $g_0$, on the annulus $\Ga_i \times (tR, \ep)$ it is $g_{V_i}$, and on $K_i \subset Y_i$ it is $t^2 g_{Y_i}$. It follows that for $x \in M_0 \setminus \bigcup_{i=1}^{n} S_i$, we have $\da(g_t)_x = \da(g_0)_x \geq C_1$ and $| R(g_t) |_{g_t} = | R(g_0) |_{g_0} \leq C_2$ for some constant $C_1, C_2 > 0$, as the metric here is independent of $t$. For $x \in \Ga_i \times (tR, \ep)$, we have $\da(g_t)_x = \da(g_{V_i})_x \geq C_3 r$ and $| R(g_t) |_{g_t} = | R(g_{V_i}) |_{g_{V_i}} \leq C_4 r^{-2}$ for some constant $C_3, C_4 > 0$, as the length scale for the cone metric is given by $r$. Finally for $x \in K_i \subset Y_i$, we have $\da(g_t)_x = \da(t^2 g_{Y_i})_x \geq C_5 t$ and $| R(g_t) |_{g_t} = | R(t^2 g_{Y_i}) |_{t^2 g_{Y_i}} \leq C_6 t^{-2}$ for some constant $C_5, C_6 > 0$, as the length scale for the metric $t^2 g_{Y_i}$ is given by $t$. Thus from the explicit definition of $\rho_t$ in (\ref{eqn:defofrhotonMt}), (1) and (2) hold with $D_2 = \textnormal{min}(C_1 \ep^{-1}, C_3, C_5 R^{-1})$ and $D_3 = \textnormal{max}(C_2 \ep^2, C_4, C_6 R^2)$. \\

\par Condition (3) holds with small enough $D_2 << 1/2$, thus by making $D_2$ smaller if necessary, $\rho_t$ satisfies (3) as well. Therefore Theorem \ref{thm:6foldestimatesmodified} applies to $(M_t, g_t)$ and $\rho_t$. Let $D_4 = R$, then $\rho_t \geq D_4 t $ on $M_t$, and we thus have a lower bound for $\rho_t$. It follows that Theorem \ref{thm:6foldsequencemodified} also applies to $(M_t, g_t)$ and $\rho_t$. \\

\par As proved in Theorem G2, the sequence $\{ \eta_j \}$ converges to $\eta $ in some Sobolev space of $\Lda^2 T^* (\SM)$. From part (c) of Theorem \ref{thm:6foldsequencemodified}, we deduce that $\| \rho^3 \, d\eta \|_{C^0} \leq K t^{3+\kp}$, that is, $| d\eta |_{g} = O(t^{3+\kp}\, \rho^{-3}) $. Thus on our 7-fold $\SMt$, we have $| d\eta_t |_{g_{\vi_t}} = O(t^{3+\kp}\, \rho_t^{-3}) $, where the norm is measured by the metric $g_{\vi_t}$ associated to the $G_2$ 3-form $\vi_t = ds \we \oa_t + \textnormal{Re}(\Oa_t)$. Using (\ref{eqn:defofrhotonMt}) and the fact that the metrics $g_{\vi_t}$ and $ds^2 + g_t$ on $\SMt$ are uniformly equivalent (see Lemma \ref{lem:vigisG2}), we obtain for $(s,x) \in \SMt$, \\[-0.8cm]

\begin{align*}
\big| \, (d\eta_t)_{(s,x)} |_{T_x M_t}  \big|_{g_t} = 
\begin{cases}
O(t^{3+\kp}) \ \ &\textnormal{for}\ x \in M_0 \setminus \bigcup_{i=1}^{n} S_i ,  \\
O(t^{3+\kp} r^{-3}) \ \ &\textnormal{for}\ x \in \Ga_i \times (tR, \ep), $\ \ for\ \ $ i = 1, \dots, n,  \tag{8.7} \label{eqn:sizeofdetatsx} \\ 
O(t^{\kp}) \ \ &\textnormal{for}\ x \in K_i \subset Y_i, $\ \ for\ \ $ i = 1, \dots, n,   
\end{cases} \\[-0.5cm]  
\end{align*} 
which then implies the improved $C^0$-estimate of $\tilde{\Oa}_t - \Oa_t$ given by: \\

\renewcommand{\theprop}{8.4} 
\begin{prop}
\label{prop:secondtermofbasicestimate}
\ In the situation above, the error term $(\tilde{\Oa}_t - \Oa_t)|_{N_t}$ satisfies 

\begin{align*}
\big| \, (\tilde{\Oa}_t - \Oa_t)|_{N_t} \,  \big|_{h_t} = 
\begin{cases}
O(t^{3+\kp}) \ \ &on\ \ N_t \cap (M_0 \setminus \bigcup_{i=1}^{n} S_i) ,  \\
O(t^{3+\kp} r^{-3}) \ \ &on\ \ N_t \cap (\Ga_i \times (tR, \ep)), $\ \ for\ \ $ i = 1, \dots, n,   \\ 
O(t^{\kp}) \ \ &on\ \ N_t \cap K_i , $\ \ for\ \ $ i = 1, \dots, n.   
\end{cases} \\[-0.5cm]  
\end{align*} 
\end{prop}

\par To finish the basic estimate, it remains to compute the error term $(\hat{\Oa}_t - \tilde{\Oa}_t)|_{N_t}$. This term arises from changing the coordinates on $M_t$ by applying the diffeomorphism $\psi_t$ obtained from Moser's argument. We claim that this term is of the same order as the term $(\tilde{\Oa}_t - \Oa_t)|_{N_t}$. \\

\par Recall that we use Moser's argument to construct the diffeomorphism $\psi_t : M_t \longrightarrow M_t$ so that $\psi_t^* (c_t \tilde{\oa}_t) = \oa_t$, and we write $\hat{\Oa}_t = \psi_t^* (\tilde{\Oa}_t)$. Thus the difference between $\hat{\Oa}_t$ and $\tilde{\Oa}_t$ is essentially given by the term ``$\pa (\psi_t - \textnormal{Id})$'', where Id denotes the identity map on $M_t$. Here what we mean by the difference between $\psi_t$ and Id can be interpreted in terms of local coordinates on $M_t$, and we use $\pa$ to denote the usual partial differentiation. \\

\par Now as $c_t \tilde{\oa}_t $ and $\oa_t$ are cohomologous, we write $c_t \tilde{\oa}_t - \oa_t = d\si_t$ for some smooth 1-form $\si_t$. Note that $c_t \tilde{\oa}_t - \oa_t = \iota(\frac{\pa}{\pa s})(\tilde{\vi}_t - \vi_t) = \iota(\frac{\pa}{\pa s})d\eta_t $ is just a component of $d\eta_t$, and so $|d\si_t|_{g_t}$ is given by (\ref{eqn:sizeofdetatsx}). We claim that we can choose a small 1-form $\si_t$ uniquely on $M_t$, so that Moser's argument defines ``small'' vector fields $X_t$ (see also in the proof of Theorem 4.9 of \cite{Chan1}), and then constructs ``small'' diffeomorphisms $\psi_t$ by representing them as the flow of $X_t$ on $M_t$. Our technique is to adopt a kind of \emph{isoperimetric inequality} which is similar to the one in (v) of Theorem \ref{thm:IIIthm5.3}, but we are working with 1-forms on the real 6-folds $M_0'$ and $Y_i$. To control the length of this paper, we state here without proof the following two results: (the proof can be found in the author's thesis, and the essential tools are the theory of \emph{Weighted Sobolev spaces} on manifolds with ends due to Lockhart and McOwen \cite{LockMcOwen}, including the ``weighted" version of Sobolev embedding and elliptic regularity theorems) \\


\renewcommand{\theprop}{8.5} 
\begin{prop}
\label{prop:isoperonM0}
\ There exists a constant $C_1 > 0$ such that $\| \si \|_{L^3} \leq C_1 ( \| d\si \|_{L^2} + \| d^* \si \|_{L^2} )$ 
for all smooth compactly-supported 1-forms $\si$ on $M_0'$. \\[-0.3cm] 
\end{prop}

\par A similar result holds on AC \CY 3-fold $Y_i$:

\renewcommand{\theprop}{8.6} 
\begin{prop}
\label{prop:isoperonYi}
\ There exists a constant $C_2 > 0$ such that $\| \si \|_{L^3} \leq C_2 ( \| d\si \|_{L^2} + \| d^* \si \|_{L^2} )$ 
for all smooth compactly-supported 1-forms $\si$ on $Y_i$ for $i=1, \dots, n$. \\ 
\end{prop}

\par We remark here that the inequality in Proposition \ref{prop:isoperonYi} is invariant under homotheties, which means the inequality also holds on $(Y_i, J_{Y_i}, t^2 \oa_{Y_i}, t^3 \Oa_{Y_i})$ with the same constant. Now take $C$ = max($C_1, C_2$), we have $\| \si \|_{L^3} \leq C ( \| d\si \|_{L^2} + \| d^* \si \|_{L^2} )$ for 1-forms $\si$ on $M_0'$ and $Y_i$ for $i=1, \dots, n$.   \\



\par We now proceed to ``glue'' together the inequalities on $M_0'$ and $Y_i$ to obtain an inequality for 1-forms on $M_t$ for small $t>0$. Choose $u,v > 0$ with $v < u < \al$ such that $tR < t^{\al} < 2t^{\al} < t^u < t^v < \ep$ for small $t > 0$. Let $H: (0, \infty) \longrightarrow [0,1]$ be a smooth decreasing function so that $H(s) = 1$ for $s \in (0, v]$, and $H(s) = 0$ for $s \in [u, \infty)$. Define a function $G_t : M_t \longrightarrow [0,1]$ by $G_t(x) = 1 $ for  $x \in M_0 \setminus \bigcup_{i=1}^{n} S_i$, $G_t(x) = H(\textnormal{log }r/\textnormal{log }t) $ for  $x \in \Ga_i \times (tR, \ep)$, $i=1, \dots, n$, and $G_t(x) = 0 $ for  $x \in K_i \subset Y_i$, $i=1, \dots, n$. Observe that $G_t \equiv 0$ on $K_i$ and $\Ga_i \times (tR, t^u)$ for $i=1, \dots, n$, and $G_t \equiv 1$ on $\Ga_i \times (t^v, \ep)$ and $M_0 \setminus \bigcup_{i=1}^{n} S_i$ for $i=1, \dots, n$. \\

\par Let $\si_t$ be a smooth 1-form on $M_t$ with $d\si_t = c_t \tilde{\oa}_t - \oa_t$. Since $H^1(M_t, \BC) = 0$ from a general fact on compact \CY manifolds, so $\si_t$ is automatically orthogonal to the space of closed and coclosed 1-forms on $M_t$, and it follows that we can choose $\si_t$ uniquely by requiring $d^* \si_t = 0$, where $d^*$ is computed using the metric $g_t$. \\

\par Write $\si_t = G_t \si_t + (1-G_t)\si_t$. We can regard $G_t \si_t$ as a compactly-supported 1-form on $M_0'$, and $(1-G_t)\si_t$ a sum of compactly-supported 1-forms on $Y_i$. Applying Proposition \ref{prop:isoperonM0} to $G_t \si_t$, using the metric $g_0$ on the support of $G_t$ and putting the constant $C$ = max($C_1, C_2$) gives 
\begin{align*}
\| G_t \si_t \|_{L^3} \ &\leq \ C \big( \, \| d(G_t \si_t) \|_{L^2} \, + \, \| d^*(G_t \si_t) \|_{L^2} \big) \\
&\leq \ C \big( \, \| G_t \,d\si_t \|_{L^2} \, + \, \| dG_t \we \si_t \|_{L^2} \, + \, \| dG_t \|_{L^6} \cdot \| \si_t \|_{L^3} \, + \, \| G_t \,d^* \si_t \|_{L^2} \big) \\
&\leq \ C \big( \,\| G_t \,d\si_t \|_{L^2} \, + \, 2 \| dG_t \|_{L^6} \cdot \| \si_t \|_{L^3} \big),
\end{align*}
where we have used $d^* \si_t = 0$ and H$\ddot{\text{o}}$lder's inequality. The same inequality holds with the metric $g_t$, as it coincides with $g_0$ on the support of $G_t$. For the 1-form $(1-G_t)\si_t$, since the support of $1-G_t$ in $Y_i$ is $K_i \cup (\Ga_i \times (tR, t^{v}))$ for $i=1,\dots,n$, we apply Proposition \ref{prop:isoperonYi}, using the constant $C$, and get 
\begin{align*}
&\| (1-G_t)\si_t |_{K_i \cup (\Ga_i \times (tR, t^{v}))} \|_{L^3} \\
\leq \ &C \big( \, \| d((1-G_t)\si_t) |_{K_i \cup (\Ga_i \times (tR, t^{v}))} \|_{L^2} \, + \, \| d^* ((1-G_t)\si_t) |_{K_i \cup (\Ga_i \times (tR, t^{v}))} \|_{L^2} \big) \\
\leq \ &C \big( \, \| (1-G_t)d\si_t |_{K_i \cup (\Ga_i \times (tR, t^{v}))} \|_{L^2} \, + \, 2 \| dG_t |_{K_i \cup (\Ga_i \times (tR, t^{v}))} \|_{L^6} \cdot \| \si_t |_{K_i \cup (\Ga_i \times (tR, t^{v}))} \|_{L^3} \\ 
&+ \, \| (1-G_t)d^* \si_t |_{K_i \cup (\Ga_i \times (tR, t^{v}))} \|_{L^2} \big)
\end{align*}
using the metric $t^2 g_{Y_i}$. As the metric $t^2 g_{Y_i}$ coincides with $g_t$ for $r \leq t^{\al}$, and is close to it for $t^{\al} \leq r \leq t^v$, $d^* \si_t$ equals zero for $r \leq t^{\al}$, and is small for $t^{\al} \leq r \leq t^v$, computed using $t^2 g_{Y_i}$. Thus by increasing $C$, we have 
\begin{align*}
\| (1-G_t) \si_t |_{K_i \cup (\Ga_i \times (tR, t^v))} \|_{L^3} \ \leq \ C &\big( \, \| ((1-G_t) d\si_t) |_{K_i \cup (\Ga_i \times (tR, t^v))} \|_{L^2} \\
&+ \, 2  \| dG_t |_{K_i \cup (\Ga_i \times (tR, t^v))}  \|_{L^6} \cdot \| \si_t |_{K_i \cup (\Ga_i \times (tR, t^v))} \|_{L^3}   \big),
\end{align*}
computed using $g_t$. It follows that 
\begin{align*}
\| (1-G_t) \si_t \|_{L^3} \ &\leq \ \sum_{i=1}^{n}\,\| (1-G_t) \si_t |_{K_i \cup (\Ga_i \times (tR, t^v))} \|_{L^3} \\
&\leq \ C\,\sum_{i=1}^{n}\,\big( \, \| ((1-G_t) d\si_t) |_{K_i \cup (\Ga_i \times (tR, t^v))} \|_{L^2} \\ 
&\qquad \qquad + \, 2  \| dG_t |_{K_i \cup (\Ga_i \times (tR, t^v))}  \|_{L^6} \cdot \| \si_t |_{K_i \cup (\Ga_i \times (tR, t^v))} \|_{L^3}   \big)\\
&\leq \ C\,\sqrt{n}\,\big( \, \| (1-G_t) d\si_t \|_{L^2} \, + \, 2 \| dG_t \|_{L^6} \cdot \| \si_t \|_{L^3} \big), 
\end{align*}
where we used the inequality of arithmetic and geometric means on the last row. Consequently we have 
\begin{align*}
\| \si_t \|_{L^3} &\leq \| G_t \si_t \|_{L^3} \ + \ \| (1-G_t)\,\si_t \|_{L^3} \\
&\leq  C \big( \,\| G_t \,d\si_t \|_{L^2}  +  2 \| dG_t \|_{L^6} \cdot \| \si_t \|_{L^3} \big)  +  C\,\sqrt{n}\,\big( \| (1-G_t)\,d\si_t \|_{L^2}  +  2 \| dG_t \|_{L^6} \cdot \| \si_t \|_{L^3} \big),
\end{align*}
which implies
\begin{align*}
\big( 1 - 2C(1 + \sqrt{n})\,\| dG_t \|_{L^6} \big) \cdot \| \si_t \|_{L^3} \ &\leq \ C \| G_t \,d\si_t \|_{L^2} \ + \ C\,\sqrt{n}\,\| (1-G_t)\,d\si_t \|_{L^2} \\
&\leq C(1 + \sqrt{n})\,\| d\si_t \|_{L^2}
\end{align*}
as $\| G_t \,d\si_t \|_{L^2}$, $\| (1-G_t)\,d\si_t \|_{L^2} \ \leq \ \| d\si_t \|_{L^2}$. Calculation shows that the $L^6$-norm of $dG_t$ is given by $O(|\text{log }t|^{-5/6})$, which tends to zero as $t \rightarrow 0$. Thus for sufficiently small $t>0$, we can make      
$$ 2C(1 + \sqrt{n})\,\| dG_t \|_{L^6} \leq 1/2. $$
Therefore 
$$ \| \si_t \|_{L^3} \ \leq \ 2C(1 + \sqrt{n})\,\| d\si_t \|_{L^2}, $$
and hence we have proved: 

\renewcommand{\thethm}{8.7} 
\begin{thm}
\label{thm:isoperionMt}
\ Suppose $\si_t$ is a smooth 1-form on $M_t$ with $d\si_t = c_t \tilde{\oa}_t - \oa_t$ and $d^* \si_t = 0$. Then there exists a constant $K >0$, independent of $t$, such that
\begin{align*}
\| \si_t \|_{L^3} \ \leq \ K\,\| d\si_t \|_{L^2}
\end{align*}
for sufficiently small $t>0$.
\end{thm}

\par As we have seen earlier, $|d\si_t|_{g_t}$ is of the same order as $|(d\eta_t)_{(s,x)}|_{T_x M_t}|_{g_t}$, and so is given by (\ref{eqn:sizeofdetatsx}), i.e. $| d\si_t |_{g_t} = O(t^{3+\kp} \rho_t^{-3})$ on balls of radius $O(\rho_t (x))$ about $x \in M_t$, where $\rho_t$ is given in (\ref{eqn:defofrhotonMt}). Thus we have estimates of $(d+d^*)\si_t = d\si_t$ and all derivatives, given by $| \nab^l  d\si_t |_{g_t} = O(t^{3+\kp} \rho_t^{-3-l})$ for $l \geq 0$. Moreover, we have the global estimate $\| d\si_t \|_{L^2} = O(t^{3+\kp})$ on the whole $M_t$ as in Theorem \ref{thm:6foldsequencemodified}, which implies $ \| \si_t \|_{L^3} = O(t^{3+\kp})$ from Theorem \ref{thm:isoperionMt}. Now using similar arguments to the proof of Theorem \ref{thm:6foldsequencemodified}, the elliptic regularity for the operator $d+d^*$ and the global estimate $ \| \si_t \|_{L^3} = O(t^{3+\kp})$ give 
$$ | \nab^l  \si_t |_{g_t} = O(t^{3+\kp} \rho_t^{-2-l})  $$
for $l \geq 0$. It follows that we have the same estimates for $ | \nab^l  X_t |_{g_t}$, and hence for $ | \nab^l  (\psi_t - \text{Id}) |_{g_t}$ which can be interpreted using local coordinates. Now this diffeomorphism estimate for $l=1$ is sufficient to prove our expected size for the term $(\hat{\Oa}_t - \tilde{\Oa}_t)|_{N_t}$: \\

\renewcommand{\theprop}{8.8} 
\begin{prop}
\label{prop:firsttermofbasicestimate}
\ In the situation above, the error term $(\hat{\Oa}_t - \tilde{\Oa}_t )|_{N_t}$ satisfies 
\begin{align*}
\big| \, (\hat{\Oa}_t - \tilde{\Oa}_t )|_{N_t} \,  \big|_{h_t} = 
\begin{cases}
O(t^{3+\kp}) \ \ &on\ \ N_t \cap (M_0 \setminus \bigcup_{i=1}^{n} S_i) ,  \\
O(t^{3+\kp} r^{-3}) \ \ &on\ \ N_t \cap (\Ga_i \times (tR, \ep)), $\ \ for\ \ $ i = 1, \dots, n,   \\ 
O(t^{\kp}) \ \ &on\ \ N_t \cap K_i , $\ \ for\ \ $ i = 1, \dots, n.   
\end{cases} \\[-0.5cm]  
\end{align*} 
\end{prop}

\par Before proceeding to combining the errors to get the basic estimate in (\ref{eqn:basicestimate}), let us return to the issue on the uniform equivalence between the metrics $\hat{h}_t$ and $h_t$. We already got the size for $(\tilde{\Oa}_t - \Oa_t )|_{N_t}$ and $(\tilde{\oa}_t - \oa_t )|_{N_t}$, both have same order. The size for $(\hat{\Oa}_t - \tilde{\Oa}_t )|_{N_t}$ and $(\hat{\oa}_t - \tilde{\oa}_t )|_{N_t}$ are essentially given by the ``difference'' between the diffeomorphism $\psi_t$ and the identity, and we have shown that it is of the same order as the size for $(\tilde{\Oa}_t - \Oa_t )|_{N_t}$ or $(\tilde{\oa}_t - \oa_t )|_{N_t}$. It follows that the size for $(\hat{\Oa}_t - \Oa_t )|_{N_t}$ and $(\hat{\oa}_t - \oa_t )|_{N_t}$ has the same order as $(\tilde{\Oa}_t - \Oa_t )|_{N_t}$ or $(\tilde{\oa}_t - \oa_t )|_{N_t}$. This implies the metrics $\hat{h}_t$ and $h_t$ are uniformly equivalent in $t$, and therefore Propositions \ref{prop:thirdtermofbasicestimate}, \ref{prop:secondtermofbasicestimate} and \ref{prop:firsttermofbasicestimate} also hold for $\hat{h}_t$. \\

\par We summarize the above estimates from Propositions \ref{prop:thirdtermofbasicestimate}, \ref{prop:secondtermofbasicestimate} and \ref{prop:firsttermofbasicestimate} in the following table, measuring w.r.t $\hat{h}_t$: \\

\begin{tabular}{l|c c c}
&	$(\hat{\Oa}_t - \tilde{\Oa}_t )|_{N_t}$ & \ \ $(\tilde{\Oa}_t - \Oa_t)|_{N_t}$ & $\textnormal{Im}(\Oa_t)|_{N_t}$ \\
\hline \\[-0.3cm]
$N_t \cap (M_0 \setminus \bigcup_{i=1}^{n} S_i)$ & $O(t^{3+\kp})$ & $O(t^{3+\kp})$ & 0 \\[0.2cm]
$N_t \cap (\Ga_i \times (2t^{\al}, \ep))$ &  $O(t^{3+\kp} r^{-3})$ & $O(t^{3+\kp} r^{-3})$ & 0  \\[0.2cm]
$N_t \cap (\Ga_i \times (t^{\al}, 2t^{\al}))$ & $O(t^{3(1-\al)+\kp})$ & $O(t^{3(1-\al)+\kp})$ & $O(t^{(\mu-1) \al}) + O(t^{(1-\kp_i)(1-\al)})$ \\[0.2cm] 
$N_t \cap (\Ga_i \times (tR, t^{\al}))$ & $O(t^{3+\kp} r^{-3})$ & $O(t^{3+\kp} r^{-3})$ & 0 \\[0.2cm]
$N_t \cap K_i$ & $O(t^{\kp})$ & $O(t^{\kp})$ & 0  
\end{tabular} \\
\begin{center}
Table 1. \ The estimate (\ref{eqn:basicestimate}) on different regions of $N_t$ \\[0.6cm]
\end{center}



\par As in \S\ref{sec:JoyceSDesingularizationTheory} we may write $\hat{\Oa}_t|_{N_t} = e^{i\ta_t}dV_t$ for some phase function $e^{i\ta_t}$ on $N_t$. Here $dV_t$ is the volume form induced by the metric $\hat{h}_t$. Then $\textnormal{Im}(\hat{\Oa}_t)|_{N_t} = \textnormal{sin}\,\ta_t\,dV_t$ and from the table we have 
\begin{align*}
| \textnormal{sin}\,\ta_t |_{\hat{h}_t} = 
\begin{cases}
O(t^{3+\kp}) \, &\textnormal{on } N_t \cap (M_0 \setminus \bigcup_{i=1}^{n} S_i), \\
O(t^{3+\kp}r^{-3}) \, &\textnormal{on } N_t \cap (\Ga_i \times (2t^{\al}, \ep)), \\
O(t^{3(1-\al)+\kp}) \,+ \, O(t^{(\mu-1) \al}) \, + \, O(t^{(1-\kp_i)(1-\al)}) &\textnormal{on } N_t \cap (\Ga_i \times (t^{\al}, 2t^{\al})), \tag{8.8} \label{eqn:sintatonNt}  \\
O(t^{3+\kp}r^{-3}) \, &\textnormal{on } N_t \cap (\Ga_i \times (tR, t^{\al})), \\
O(t^{\kp}) &\textnormal{on } N_t \cap K_i,
\end{cases} 
\end{align*} 
for all $i = 1, \dots, n$. \\


\par We also need to estimate the derivative $d\,\textnormal{sin}\,\ta_t$. Using similar arguments it can be deduced that
\begin{align*}
| d \, \textnormal{sin}\,\ta_t |_{\hat{h}_t} = 
\begin{cases}
O(t^{3+\kp}) \, &\textnormal{on } N_t \cap (M_0 \setminus \bigcup_{i=1}^{n} S_i), \\
O(t^{3+\kp} r^{-4}) \, &\textnormal{on } N_t \cap (\Ga_i \times (2t^{\al}, \ep)), \\
O(t^{3(1-\al)+\kp-\al}) + O(t^{(\mu-1) \al -\al}) +  O(t^{(1-\kp_i)(1-\al)-\al}) \! \! \! \! \! &\textnormal{on } N_t \cap (\Ga_i \times (t^{\al}, 2t^{\al})),  \tag{8.9} \label{eqn:dsintatonNt-3}   \\
O(t^{3+\kp}r^{-4}) \, &\textnormal{on } N_t \cap (\Ga_i \times (tR, t^{\al})), \\
O(t^{\kp-1}) &\textnormal{on } N_t \cap K_i,
\end{cases}
\end{align*}  
for all $i = 1, \dots, n$. Here we used equations (\ref{eqn:sizeforai}) 
and (\ref{eqn:sizeforbi(t-1r)}) to obtain the bound on $N_t \cap (\Ga_i \times (t^{\al}, 2t^{\al}))$. \\

\par As in part (i) of Theorem \ref{thm:IIIthm5.3}, we need bounds for $\| \textnormal{sin}\,\ta_t \|_{L^{6/5}}$,  $\| \textnormal{sin}\,\ta_t \|_{C^0}$ and  $\| d\,\textnormal{sin}\,\ta_t \|_{L^6}$, computing norms w.r.t. $\hat{h}_t$. Since vol($N_t \cap K_i$) = $O(t^3)$ and vol($N_t \cap (\Ga_i \times (t^{\al}, 2t^{\al}))$) = $O(t^{3\al})$, and vol($N_t \cap (M_0 \setminus \bigcup_{i=1}^{n} S_i)$) = $O(1)$, it follows that 

\begin{align*}
\| \textnormal{sin}\,\ta_t \|_{L^{6/5}} \, &= \, O(1)^{5/6} \cdot O(t^{3+\kp}) \, + \, O \left( \sum_{i=1}^{n}\,\textnormal{vol}(\Si_i)^{5/6} \left( \int_{2t^{\al}}^{\ep} (t^{3+\kp} r^{-3} )^{6/5}\, r^2 dr \right)^{5/6} \right) \\ 
&\quad \, + \, O(t^{3\al})^{5/6} \cdot \big( \, O(t^{3(1-\al)+\kp}) \, + \, O(t^{(\mu-1) \al}) \, + \, \sum_{i=1}^{n} O(t^{(1-\kp_i)(1-\al)}) \, \big) \\
&\quad \, + \, O \left( \sum_{i=1}^{n}\,\textnormal{vol}(\Si_i)^{5/6} \left( \int_{tR}^{t^{\al}} (t^{3+\kp} r^{-3} )^{6/5}\, r^2 dr \right)^{5/6} \right) \, + \, O(t^3)^{5/6} \cdot O(t^{\kp}) \\
&= \, O(t^{3(1-\al)+\kp+5\al/2}) \, + \, O(t^{(\mu-1) \al + 5\al/2}) \, + \, \sum_{i=1}^{n} O(t^{(1-\kp_i)(1-\al) + 5\al/2})  \\ 
&\ \quad + \, O(t^{\kp + 5/2})  \tag{8.10} \label{eqn:L6/5normofsintat}
\end{align*} 
using (\ref{eqn:sintatonNt}). Similarly, we have 
\begin{align*}
\| \textnormal{sin}\,\ta_t \|_{C^0} \, &= \, O(t^{3+\kp}) \, + \, O(t^{3+\kp} (2t^{\al})^{-3} ) \, + \, O(t^{3(1-\al)+\kp}) \, + \, O(t^{(\mu-1) \al }) \\
&\quad \, + \, \sum_{i=1}^{n} O(t^{(1-\kp_i)(1-\al) }) \, + \, O(t^{3+\kp} (tR)^{-3} ) \, + \, O(t^{\kp}) \\
&= \, O(t^{3(1-\al)+\kp}) \, + \, O(t^{(\mu-1) \al }) \, + \, \sum_{i=1}^{n} O(t^{(1-\kp_i)(1-\al) }) \, + \, O(t^{\kp}). \tag{8.11}
\label{eqn:C0normofsintat}
\end{align*} \\
Using the estimate (\ref{eqn:dsintatonNt-3}) for the derivative, we have 
\begin{align*}
\| d \, \textnormal{sin}\,\ta_t \|_{L^6} \, &= \, O(1)^{1/6} \cdot O(t^{3+\kp}) \, + \, O \left( \sum_{i=1}^{n}\,\textnormal{vol}(\Si_i)^{1/6} \left( \int_{2t^{\al}}^{\ep} (t^{3+\kp} r^{-4} )^6 \, r^2 dr \right)^{1/6} \right) \\ 
&\quad \, + \, O(t^{3\al})^{1/6} \cdot \big( \, O(t^{3(1-\al)+\kp-\al}) \, + \, O(t^{(\mu-1) \al -\al}) \, + \, \sum_{i=1}^{n} O(t^{(1-\kp_i)(1-\al)-\al}) \, \big) \\
&\quad \, + \, O \left( \sum_{i=1}^{n}\,\textnormal{vol}(\Si_i)^{1/6} \left( \int_{tR}^{t^{\al}} (t^{3+\kp} r^{-4} )^6 \, r^2 dr \right)^{1/6} \right) \, + \, O(t^3)^{1/6} \cdot O(t^{\kp-1}) \\
&= \, O(t^{3(1-\al)+\kp-\al/2}) \, + \, O(t^{(\mu-1) \al - \al/2}) \, + \, \sum_{i=1}^{n} O(t^{(1-\kp_i)(1-\al) - \al/2})  +  O(t^{\kp - 1/2}).   \tag{8.12}
\label{eqn:L6normofsintat}
\end{align*}  
Now for part (i) of Theorem \ref{thm:IIIthm5.3} to hold, we need:
\begin{align*}
\begin{cases}
\kp + 5/2 \geq \kp' + 3/2,   &3(1-\al) + \kp + 5\al/2 \geq \kp' + 3/2, \\[-0.3cm]
\tag{8.13} \label{eqn:checkineqL6/5} \\[-0.3cm]
(\mu-1)\al + 5\al/2 \geq \kp' + 3/2, \qquad \textnormal{and } &(1-\kp_i)(1-\al) + 5\al/2 \geq \kp' + 3/2
\end{cases}
\end{align*}  
from (\ref{eqn:L6/5normofsintat}), 
\begin{align*}
\begin{cases}
\kp \geq \kp' - 1,   &3(1-\al) + \kp \geq \kp' - 1, \\[-0.3cm]
\tag{8.14} \label{eqn:checkineqC0} \\[-0.3cm]
(\mu-1)\al \geq \kp' - 1, \qquad \textnormal{and } &(1-\kp_i)(1-\al) \geq \kp' - 1
\end{cases}
\end{align*}  
from (\ref{eqn:C0normofsintat}), and 
\begin{align*}
\begin{cases}
\kp - 1/2 \geq \kp' - 3/2,   &3(1-\al) + \kp - \al/2 \geq \kp' - 3/2, \\[-0.3cm]
\tag{8.15} \label{eqn:checkineqL6} \\[-0.3cm]
(\mu-1)\al - \al/2 \geq \kp' - 3/2, \qquad \textnormal{and } &(1-\kp_i)(1-\al) - \al/2 \geq \kp' - 3/2
\end{cases}
\end{align*} 
from (\ref{eqn:L6normofsintat}). \\

\par Calculations show that given $\kp > 0$, $\mu > 1$ and $\kp_i < -3/2$, we can choose $\kp' > 1$ close to 1, and $\al \in (0,1)$ close to 1, such that (\ref{eqn:checkineqL6/5}), (\ref{eqn:checkineqC0}) and (\ref{eqn:checkineqL6}) hold. Here is the place where we need to assume the rate $\kp_i$ of AC SL 3-folds $L_i$ to be less than $-3/2$. As we have fixed the rate $\lda_i$ of the AC \CY 3-fold $Y_i$ to satisfy $\lda_i \leq -3$, and also we require $\kp_i > \lda_i + 1$ in Definition \ref{def:ACSLinACCY}, so that assuming $\kp_i < -3/2$ is still possible. \\

\par Therefore, we have shown that there exist $\kp' > 1$ and $A_2 > 0$ such that $\| \textnormal{sin}\,\ta_t \|_{L^{6/5}} \leq A_2 t^{\kp'+3/2} $, $\| \textnormal{sin}\,\ta_t \|_{C^0} \leq A_2 t^{\kp' - 1} $ and $\| d\,\textnormal{sin}\,\ta_t \|_{L^6} \leq A_2 t^{\kp'-3/2} $ for sufficiently small $t > 0$, i.e. (i) of Theorem \ref{thm:IIIthm5.3} holds for $N_t$. \\

\begin{center}
\renewcommand{\thesubsection}{\normalfont{9}}
\subsection{\normalfont{\textsc{Desingularizations of $N_0$ }}}
\label{sec:DesingularizationsOfN0}
\end{center}


\par This section gives the main result of the chapter, the desingularizations of SL 3-folds $N_0$ with conical singularities. The proof of it is based on an analytic existence theorem for SL 3-folds, Theorem \ref{thm:IIIthm5.3}, which is adapted from Joyce's result \cite[Thm. 5.3]{Joycedesing3}. We have already verified part (i) of Theorem \ref{thm:IIIthm5.3} in \S\ref{sec:EstimatesOfImTildeOaTNT}, and it remains to check (ii) to (v) hold for the Lagrangian 3-folds $N_t$ we constructed. \\

\renewcommand{\thethm}{9.1} 
\begin{thm}
\label{thm:mainthmonSLdesing}
\ Suppose $(M_0, J_0, \oa_0, \Oa_0)$ is a compact \CY 3-fold with finitely many conical singularities at $x_1, \dots, x_n$ with rate $\nu > 0$ modelled on \CY cones $V_1, \dots , V_n$. Let $(Y_1, J_{Y_1}, \oa_{Y_1}, \Oa_{Y_1}), \dots, (Y_n, J_{Y_n}, \oa_{Y_n}, \Oa_{Y_n})$ be AC \CY 3-folds with rates $\lda_1, \dots,$ $\lda_n < -3$ modelled on the same \CY cones $V_1, \dots , V_n$. \\

\par Then Theorem \ref{thm:desingthmforlda<-3} gives a family of \CY 3-folds $(M_t, \tilde{J}_t, \tilde{\oa}_t, \tilde{\Oa}_t)$ for sufficiently small $t>0$. Apply the diffeomorphism $\psi_t : M_t \longrightarrow M_t$ on $M_t$ to get \CY structures $(\hat{J}_t, \hat{\oa}_t, \hat{\Oa}_t)$, as in $\S\ref{sec:ConstructionOfNT}$. \\

\par Let $N_0$ be a compact SL 3-fold in $M_0$ with the same conical singularities at $x_1, \dots, x_n$ with rate $\mu \in (1, \nu + 1)$ modelled on SL cones $C_1, \dots, C_n$. Suppose $N_0 \setminus \{ x_1, \dots, x_n  \}$ is connected. Let $L_1, \dots, L_n$ be AC SL 3-folds in $Y_1, \dots, Y_n$ with rates $\kp_1 \in (\lda_1 + 1, -3/2), \dots, \kp_n \in (\lda_n + 1, -3/2)$ modelled on the same SL cones $C_1, \dots, C_n$. \\

\par Then there exists a family of compact nonsingular SL 3-folds $\hat{N}_t$ in $(M_t, \hat{J}_t, \hat{\oa}_t, \hat{\Oa}_t)$ for sufficiently small $t$, such that $\hat{N}_t$ is constructed by deforming the Lagrangian 3-fold $N_t$ which is made by gluing $L_i$ into $N_0$ at $x_i$ for $i = 1, \dots, n$. \\
\end{thm}
\pf \ First of all, we have to check that $N_t$ satisfies the conditions in \S\ref{sec:JoyceSDesingularizationTheory}. 
Let us start with evaluating the integral $\int_{N_t} \textnormal{Im}(\hat{\Oa}_t)$. Calculation using (\ref{eqn:sintatonNt}) shows that 
$$ \int_{N_t} \textnormal{Im}(\hat{\Oa}_t) \ = \ O(t^{3+\tau}) $$
for some $\tau \in (0, \kp)$. 
As mentioned in \S\ref{sec:JoyceSDesingularizationTheory}, we can rescale the phase for $\hat{\Oa}_t$ by $\hat{\Oa}_t \mapsto e^{i\zeta_t} \, \hat{\Oa}_t$ such that $\int_{N_t} \textnormal{Im}(e^{i\zeta_t} \hat{\Oa}_t) = 0$. Thus we have $\textnormal{sin}\,\ta_t \mapsto \textnormal{sin}\,(\ta_t + \zeta_t) \approx \textnormal{sin}\,(\ta_t) + \zeta_t $, and the size for the term $\zeta_t$ is approximately given by the ratio between $\int_{N_t} \textnormal{Im}(\hat{\Oa}_t)$ and $\int_{N_t} \textnormal{Re}(\hat{\Oa}_t)$. Now since $\int_{N_t} \textnormal{Re}(\hat{\Oa}_t) \approx \textnormal{vol}(N_t) = O(1) $, the correction term $\zeta_t$ essentially contributes $O(t^{3+\tau})$ 
to $\textnormal{sin}(\ta_t)$. As we have shown that $\| \textnormal{sin}(\ta_t) \|_{L^{6/5}} = O(t^{\kp'+ 3/2})$ for some $\kp' > 1$, then $\| \textnormal{sin}(\ta_t + \zeta_t) \|_{L^{6/5}} = O(t^{\kp'+ 3/2}) + O(t^{3+\tau})$. 
The term $O(t^{\kp'+ 3/2}) $ will then be dominant if $\kp'$ is close to 1. As a result, the rescaling of phases does not affect the $L^{6/5}$-estimates in (i) of Theorem \ref{thm:IIIthm5.3} at all. For the terms $\| \textnormal{sin}(\ta_t) \|_{C^0}$ and $\| d\,\textnormal{sin}(\ta_t) \|_{L^6}$, calculation shows that the rescaling of phases does not affect the estimates as well. \\

\par As we have assumed $N_0 \setminus \{ x_1, \dots, x_n  \}$ is connected, we can take the finite dimensional vector space $W$ to be the space of constant functions, i.e. $W = \left< 1 \right>$, as in \S\ref{sec:JoyceSDesingularizationTheory}. \\

\par Under our construction, $N$, $h$ and $\acute{h}$ in \S\ref{sec:JoyceSDesingularizationTheory} are replaced by $N_t$, $\hat{h}_t$ and $\acute{\hat{h}}_t$ respectively, and we thus need to show (i), (iii), (iv) and (v) hold using the metric $\hat{h}_t$ on $N_t$, and (ii) holds using the metric $\acute{\hat{h}}_t$ on $T^* N_t$. Basically the proof for (iii) and (iv) using the metric $h_t$ and for (ii) using the metric $\acute{h}_t$ can be found in \cite[Thm. 6.8]{Joycedesing3}, and the proof for (v) using the metric $h_t$ is given in \cite[Thm. 6.12]{Joycedesing3}. Thus our approach to showing (ii)-(v) in Theorem \ref{thm:IIIthm5.3} is to apply Theorems 6.8 and 6.12 in \cite{Joycedesing3} together with the uniform equivalence between the metrics $h_t$ and $\hat{h}_t$. \\

\par We have shown in \S\ref{sec:EstimatesOfImTildeOaTNT} that given $\kp > 0$, $\mu > 1$ and $\kp_i < -3/2$ for $i = 1, \dots, n$, there exists $\kp' > 1$ and $A_2 > 0$ such that (i) of Theorem \ref{thm:IIIthm5.3} holds for sufficiently small $t > 0$, measuring w.r.t. the metric $\hat{h}_t$.   \\

\par For part (v), Theorem 6.12 in \cite{Joycedesing3} shows that there exists $A_6 > 0$ such that (v) holds using the metric $h_t$. Note that the assumption on the connectedness of $N_0 \setminus \{ x_1, \dots, x_n  \}$ is used here. The fact that $v \in L^6 (N_t)$ follows from $L^2_1 (N_t) \hookrightarrow L^6 (N_t)$ by the Sobolev Embedding Theorem (see for example \cite[Thm. 1.2.1]{Joyce1}). The idea of proving the inequality for the metric $h_t$ on $N_t$ for small $t$ is to combine the Sobolev embedding inequalities on $N_0 \setminus \{ x_1, \dots, x_n  \}$ and $L_i$. Now as $h_t$ and $\hat{h}_t$ are uniformly equivalent metrics, so (v) is true for $h_t$ if and only if it is true for $\hat{h}_t$. As a result, by making $A_6$ larger if necessary, (v) holds with the metric $\hat{h}_t$. \\

\par To deduce (iii) and (iv) for $\hat{h}_t$, we first apply Theorem 6.8 in \cite{Joycedesing3} to show they are true for $h_t$ for some $A_4, A_5 > 0$. The idea of which is to consider the behaviour of the metric $h_t$ for small $t$. Since $h_t$ is $t^2 g_{Y_i}|_{L_i}$ on $H_i$ and on $\Psi_{C_i} (\Ga (du_{t,i}))$ near $\Si_i \times \{ tR' \}$ for each $i$, we have $\da(t^2 g_{Y_i}|_{L_i}) = t \da( g_{Y_i}|_{L_i} )$ and $\| R(t^2 g_{Y_i}|_{L_i}) \|_{C^0} = t^{-2} \| R(g_{Y_i}|_{L_i}) \|_{C^0} $. For small $t > 0$, the dominant contributions to $\da(h_t)$ and $\| R(h_t) \|_{C^0}$ come from $\da(t^2 g_{Y_i}|_{L_i})$ and $\| R(t^2 g_{Y_i}|_{L_i}) \|_{C^0}$ for some $i$, and hence we have $\da(h_t) = O(t)$ and $\| R(h_t) \|_{C^0} = O(t^{-2}) $. Now we prove (iii) and (iv) also hold, increasing $A_4, A_5$ if necessary, for $\hat{h}_t$ by showing the metrics $t^{-2} \hat{h}_t$ and $t^{-2} h_t$ are $C^2$-close w.r.t. $t^{-2} h_t$ (compare to the similar argument in the proof of Theorem 3.10 in \cite{Chan1}). From the estimate we know using elliptic regularity on balls of radius $O(t)$, we have $| (d\eta_t)_{(s,x)}|_{T_x M_t}  |_{g_t} = O(t^{\kp} )$ for $(s,x) \in S^1 \times M_t$. Here we do not need to use the improved estimate for $d\eta_t$ as in \S\ref{sec:EstimatesOfImTildeOaTNT}. Then we have $| (\nab^{g_t})^l \, (d\eta_t)_{(s,x)}|_{T_x M_t}  |_{g_t} = O(t^{\kp-l} )$ for $l \geq 0$, where $\nab^{g_t}$ denotes the Levi-Civita connection of $g_t$. This implies 
$$ | (\nab^{g_t})^l ( \tilde{\Oa}_t - \Oa_t ) |_{g_t} \ = \ O(t^{\kp-l} ) \ = \ | (\nab^{g_t})^l ( \tilde{\oa}_t - \oa_t ) |_{g_t} \quad \text{for $l \geq 0$,} $$
and from Moser's argument in \S\ref{sec:EstimatesOfImTildeOaTNT}, we also have 
$$ | (\nab^{g_t})^l ( \hat{\Oa}_t - \tilde{\Oa}_t ) |_{g_t} \ = \ O(t^{\kp-l} ) \ = \ | (\nab^{g_t})^l ( \hat{\oa}_t - \tilde{\oa}_t ) |_{g_t} \quad \text{for $l \geq 0$.} $$
Putting together implies 
\begin{align*}
| (\nab^{g_t})^l ( \hat{g}_t - g_t ) |_{g_t} \ = \ O(t^{\kp-l} ) \quad \text{for $l \geq 0$.} \tag{9.1}
\label{eqn:nablhatgt-gt}
\end{align*}
Denote by $\nab^{h_t}$ the Levi-Civita connection of $h_t = g_t|_{N_t}$. Then we have 
$$ \nab^{g_t} ( \hat{g}_t - g_t ) |_{N_t} \, = \, \nab^{h_t} ( \hat{h}_t - h_t ) \, + \ \text{bilinear terms in } (\hat{g}_t - g_t )|_{N_t} \, \text{and }\, T , $$
where $T$ is the second fundamental form of $N_t$ in $M_t$ w.r.t. $g_t$. The largest contribution to $|T|_{h_t}$ comes from the second fundamental form $T_{L_i}$ of $L_i$ in $Y_i$ w.r.t. $t^2 g_{Y_i}$. As $|T_{L_i}|_{g_{Y_i}|_{L_i} }$ is bounded on $L_i$, conformal rescaling then shows $|T_{L_i}|_{t^2 g_{Y_i}|_{L_i} } = O(t^{-1})$. Thus we have $|T|_{h_t} = O(t^{-1})$, and more generally $|\,(\nab^{h_t})^l T \,|_{h_t} = O(t^{-l-1})$ for $l \geq 0$. The estimate for the $l^{\text{th}}$ derivative of $\hat{h}_t - h_t$ then follows from the relation 
\begin{align*}
\left| (\nab^{g_t})^l ( \hat{g}_t - g_t )|_{N_t} \right|_{h_t} \! &= \! \left| (\nab^{h_t})^l ( \hat{h}_t - h_t ) \right|_{h_t} \! + \! O \left( \sum |\,(\nab^{h_t})^j T \,|_{h_t} \cdot \left|(\nab^{g_t})^{l-j-1} ( \hat{g}_t - g_t )|_{N_t} \right|_{h_t} \right) \\
&\ +  O \left( \sum |\,(\nab^{h_t})^{j_1} T \,|_{h_t} \cdot |\,(\nab^{h_t})^{j_2} T \,|_{h_t} \cdot \left|(\nab^{g_t})^{l- j_1 - j_2 -2} ( \hat{g}_t - g_t )|_{N_t} \right|_{h_t} \right) \\
&\ + \dots + O \big( |T|_{h_t} \cdots |T|_{h_t} \cdot ( \hat{g}_t - g_t )|_{N_t} |_{h_t} \big). 
\end{align*}
Note that the terms $O(\cdot)$ all have size $O(t^{\kp-l})$, and therefore by (\ref{eqn:nablhatgt-gt}) we see that  
\begin{align*}
\left| (\nab^{h_t})^l ( \hat{h}_t - h_t ) \right|_{h_t} \ = \ O(t^{\kp-l} ) \quad \text{for $l \geq 0$.} \tag{9.2} 
\label{eqn:nablhatht-ht}
\end{align*} 
In particular, (\ref{eqn:nablhatht-ht}) shows $| \hat{h}_t - h_t |_{h_t} $, $t\,| \nab^{h_t} (\hat{h}_t - h_t) |_{h_t} $ and $t^2 \,| (\nab^{h_t})^2 (\hat{h}_t - h_t) |_{h_t} $ are all of size $O(t^{\kp})$. It follows that $| t^{-2} \hat{h}_t - t^{-2} h_t |_{t^{-2} h_t} $, $| \nab^{t^{-2} h_t} (t^{-2} \hat{h}_t - t^{-2} h_t) |_{t^{-2} h_t} $ and $| (\nab^{t^{-2} h_t})^2 (t^{-2} \hat{h}_t - t^{-2} h_t) |_{t^{-2} h_t} $ are all of the same size $O(t^{\kp})$, where $\nab^{t^{-2} h_t}$ and $| \cdot | _{t^{-2} h_t}$ are computed using $t^{-2} h_t$. Therefore the metrics $ t^{-2} \hat{h}_t$ and $t^{-2} h_t$ are $C^2$-close w.r.t. $t^{-2} h_t$ for small $t$, and hence (iii) and (iv) are true for $\hat{h}_t$ as well. \\

\par We remain to show (ii), using the metric $\acute{\hat{h}}_t$ and the connection $\nab^{\acute{\hat{h}}_t}$ on $T^* N_t$. Here we recall the construction of $\acute{\hat{h}}_t$ and $\nab^{\acute{\hat{h}}_t}$, as in \S\ref{sec:JoyceSDesingularizationTheory}. Write $T(T^* N_t) = H_t \oplus V_t$, where $H_t \cong TN_t$ and $V_t \cong T^* N_t$ are the horizontal and vertical subbundles w.r.t. $\nab^{\hat{h}_t}$, and define $\acute{\hat{h}}_t |_{H_t} = \hat{h}_t $ and $\acute{\hat{h}}_t |_{V_t} = \hat{h}_t^{-1} $. The connection $\nab^{\acute{\hat{h}}_t}$ is given by the lift of the Levi-Civita connection $\nab^{\hat{h}_t}$ of $\hat{h}_t$ in $H_t$, and by partial differentiation in $V_t$. Following the steps in Definition 6.7 in \cite{Joycedesing3}, we define Lagrangian neighbourhoods $U_{N_t}$, $\Psi_{N_t}$ for $N_t$ by gluing together the Lagrangian neighbourhoods $U_{N_0}$, $\Psi_{N_0}$ for $N_0$ from Theorem \ref{thm:LNTforN0} and Lagrangian neighbourhoods $U_{L_i}$, $\Psi_{L_i}$ for $L_i$ from Theorem \ref{thm:LNTforLi}. The neighbourhood $U_{N_t}$ is an open tubular neighbourhood of $N_t$ in $T^* N_t$, and $\Psi_{N_t} : U_{N_t} \longrightarrow M_t $ is an embedding with $\Psi_{N_t}|_{N_t} = \textnormal{Id}$ and $\Psi_{N_t}^* (\oa_t) = \oa_{T^* N_t}$ where $\oa_{T^* N_t}$ is the canonical symplectic structure on $T^* N_t$. Recall that we have $c_t \hat{\oa}_t = \oa_t$, so $\Psi_{N_t}^* (c_t \hat{\oa}_t) = \oa_{T^* N_t}$. Now we define 3-forms $\ba_t$ and $\hat{\ba}_t$ by $\ba_t = \Psi_{N_t}^* (\textnormal{Im}(\Oa_t))$ and $\hat{\ba}_t = \Psi_{N_t}^* (\textnormal{Im}(c_t^{3/2} \hat{\Oa}_t))$. \\

\par Using arguments in Theorem 6.8 in \cite{Joycedesing3}, we have $ \| \, (\nab^{\acute{h}_t})^l \ba_t \, \|_{C^0} \leq A_3\, t^{-l} $ for $l = 0, 1, 2, 3$ on $\mathcal{B}_{A_1 t} \subset U_{N_t}$ for some $A_1, A_3 > 0$, where the norm is measuring w.r.t. $\acute{h}_t$. To prove (ii) in our case, we try to get from estimates on $\ba_t$ to estimates on $\hat{\ba}_t$. Note that $\hat{\ba}_t$ and $\ba_t$ are $C^l$-close w.r.t. $t^2 g_t$. This follows from $| (\nab^{g_t})^l ( \hat{\Oa}_t - \Oa_t ) |_{g_t} \ = \ O(t^{\kp-l} ) \ = \ | (\nab^{g_t})^l ( \hat{\oa}_t - \oa_t ) |_{g_t}$ for $l \geq 0$, which we have discussed earlier. Combining the $C^{l+1}$-closeness of $t^{-2} \hat{h}_t$ and $t^{-2} h_t$ from (\ref{eqn:nablhatht-ht}), we get a similar estimate for $\hat{\ba}_t$, using the metric $\acute{\hat{h}}_t$. Thus making $A_1$ smaller and $A_3$ larger if necessary, we obtain $\| \, (\nab^{\acute{\hat{h}}_t})^l \hat{\ba}_t \, \|_{C^0} \leq A_3\, t^{-l}$ for $l = 0, 1, 2, 3$ on $\mathcal{B}_{A_1 t} \subset U_{N_t}$, where the norm is computed using $ \acute{\hat{h}}_t$. \\

\par The theorem now follows from Theorem \ref{thm:IIIthm5.3} which shows that for sufficiently small $t >0$ we can deform $N_t$ to a nearby special Lagrangian 3-fold $\hat{N}_t = (\Psi_{N_t})_* (\Ga (df_t))$ for some $f_t \in C^{\infty} (N_t)$ with $\int_{N_t} f_t dV_t = 0$ and $\| df_t \|_{C^0} \leq Kt^{\kp'} \leq A_1 t$. This completes the proof of Theorem \ref{thm:mainthmonSLdesing}.  \hfill $\Box$ \\

\begin{center}
\renewcommand{\thesubsection}{\normalfont{10}}
\subsection{\normalfont{\textsc{An application of the desingularization theory}}}
\label{sec:ApplicationsOfTheDesingularizationTheory}
\end{center}


\par We conclude with applying the result of \S\ref{sec:DesingularizationsOfN0} to the case where the ambient \CY 3-fold $M_0$ is the \CY 3-orbifold $T^6/\BZ_3$, as described in the last part of \S\ref{subsub:Theorem on desingularizing \CY $m$-folds}. We shall take the SL 3-folds $N_0$ with conical singularities as the fixed point set of some antiholomorphic involutions on those \CY 3-folds, whereas the AC SL 3-folds $L_i$ will be taken from examples in \S\ref{sec:SomeExamplesOfACSLMFolds} inside the corresponding AC \CY 3-folds. \\ \\
\textbf{Example 10.1} \ Take the \CY 3-fold $M_0$ with conical singularities to be the \CY 3-orbifold $T^6/\BZ_3$. Applying our desingularization result in Theorem \ref{thm:desingthmforlda<-3}, we can desingularize $T^6/\BZ_3$ by gluing in AC \CY 3-folds $K_{\BC \BP^2}$ at the singular points, obtaining the crepant resolution of $T^6/\BZ_3$. \\

\par Now we produce examples of SL 3-folds $N_0$ with conical singularities in $T^6/\BZ_3$ by using the fixed point set of an antiholomorphic isometric involution, a well-known way of producing special Lagrangians in \CY manifolds. Recall that an antiholomorphic isometric involution of a \CY manifold $(M, J, \oa, \Oa)$ is a diffeomorphism $\si : M \longrightarrow M$ such that $\si^2 = \text{Id}$, $\si^*(J) = -J$, $\si^*(\oa) = -\oa$, $\si^*(\Oa) = \bar{\Oa}$ and $\si^*(g) = g$, where $g$ is the associated \CY metric. \\

\par Let $\si_0 : T^6 \longrightarrow T^6$ be the complex conjugation given by
$$ \si_0 :  (z_1, z_2, z_3) + \Lda \longmapsto (\bar{z}_1, \bar{z}_2, \bar{z}_3) + \Lda, $$
which is well-defined as $\bar{\Lda} = \Lda$. The fixed points of $\si_0$ satisfy $z_j = \bar{z}_j + a_j + b_j \za$ for some $a_j, b_j \in \BZ$. Write $z_j = x_j + y_j \za$ for $x_j, y_j \in \BR$. It follows that $z_j = x_j + a_j \za$ and $b_j = 2a_j$, and hence the fixed point set of $\si_0$ is given by
$$  \{ (x_1, x_2, x_3) + \Lda : x_j \in \BR \}, $$ 
and is then topologically a $T^3$. \\

\par Since $\si_0 \cdot \za \cdot \si_0^{-1} = \za^{-1}$, the map $\si_0$ on $T^6$ induces a conjugation $\si$ on $T^6/\BZ_3$ which is given by
$$ \si : \BZ_3 \cdot (z_1, z_2, z_3) + \Lda \longmapsto \BZ_3 \cdot (\bar{z}_1, \bar{z}_2, \bar{z}_3) + \Lda. $$ 

\par Observe that $\si_0$ swaps $\za T^3$ and $\za^2 T^3$, and fixes $T^3$ as above. But in the orbifold level, $T^3$, $\za T^3$ and $\za^2 T^3$ are the same. Moreover, it is not hard to see that the map $T^6 \longrightarrow T^6/\BZ_3$ is injective when restricted to $T^3 \subset T^6$, which means the image of $T^3$ in $T^6/\BZ_3$ is homeomorphic to $T^3$. As a result, the fixed point set of $\si$ is topologically a $T^3$, and is given by $(\BZ_3 \cdot T^3) / \BZ_3$, which is then our SL 3-fold $N_0$ in $M_0 = T^6/\BZ_3$.  \\ 

\par It is worth knowing how $\si$ acts on those 27 orbifold singular points, and see how many of them are being fixed by $\si$, which will then be the singular points of the SL 3-fold $N_0$. It turns out that $\si$ only fixes the point $\BZ_3 \cdot (0,0,0) + \Lda$, i.e. 0 in $T^6/\BZ_3$, and swaps the other 26 points in pairs, for example, $\BZ_3 \cdot (0, \frac{i}{\sqrt{3}}, \frac{2i}{\sqrt{3}}) + \Lda \longleftrightarrow \BZ_3 \cdot (0, \frac{2i}{\sqrt{3}}, \frac{i}{\sqrt{3}}) + \Lda $. This means that we have constructed an SL 3-fold $N_0$, which is topologically a $T^3$, with one singular point at $0$ in $T^6/\BZ_3$, modelled on the SL cone $(\BZ_3 \cdot \BR^3) / \BZ_3$ in $\BC^3 / \BZ_3$. \\

\par To desingularize this $N_0$ we glue in at the singular point some appropriate pieces of AC SL 3-folds in the AC \CY 3-folds $K_{\BC \BP^2}$, the canonical bundle over $\BC \BP^2$. As we have discussed in \S\ref{subsub:Theorem on desingularizing \CY $m$-folds}, the \CY desingularization we obtain is the crepant resolution of orbifold $T^6/\BZ_3$, and so the nonsingular SL 3-folds we constructed will sit inside this crepant resolution. Our first example of an AC SL 3-fold will be taken from Example \ref{eg:fpsinKCPm-1} in which the real line bundle $K_{\BR \BP^2}$ over $\BR \BP^2$ is constructed as the fixed point set of an antiholomorphic isometric involution. By gluing this $K_{\BR \BP^2}$ into $N_0 = (\BZ_3 \cdot T^3) / \BZ_3$ at the singular point, we obtain a nonsingular SL 3-fold in the crepant resolution of $M_0 = T^6/\BZ_3$. Topologically, what we obtain will be a real blow-up of $T^3$ at a point, i.e. replacing a point by an $\BR \BP^2$, which can also be interpreted as a $T^3 \# \BR \BP^3$. As we have discussed in Example \ref{eg:fpsinKCPm-1}, the AC SL 3-fold $K_{\BR \BP^2}$ has rate $\kp = -\infty$, and in order to fit into our desingularization theorem, we could choose the rate for $K_{\BR \BP^2}$ to be any $\kp \in (-5, -3/2)$. \\

\par The next example of AC SL 3-folds in $K_{\BC \BP^2}$ is given by Example \ref{eg:SO3invarinKCPm-1}. There we have constructed a family of SO(3)-invariant SL 3-folds $L_c$ diffeomorphic to $S^2 \times \BR$ which converges to two copies of the cone $((\BZ_3 \cdot \BR^3) / \BZ_3)$ in $\BC^3 / \BZ_3$. If we now take $N_0$ to be a connected double cover of $T^3$ so that we have two singular points in the same place in $M_0 = T^6/\BZ_3$ (and $N_0 \setminus \{0\}$ is connected), we can desingularize this $N_0$ by doing the connected sum with an $S^2 \times \BR$, obtaining an SL 3-fold which is homeomorphic to a $T^3 \# (S^1 \times S^2)$ in the crepant resolution of $T^6/\BZ_3$. We mentioned in Example \ref{eg:SO3invarinKCPm-1} that one possible $L_c$ will be a double cover of $K_{\BR \BP^2}$, which means the family of nonsingular SL 3-folds we constructed here will be deformations of a double cover of the SL 3-fold $T^3 \# \BR \BP^3$ in the first example. Notice also that each $L_c$ has rate $-2$, as discussed in Example \ref{eg:SO3invarinKCPm-1}, and so our desingularization theorem works. \\ \\

\renewcommand\refname{\begin{center} \normalfont{\textsc{References}}\end{center}}

\hfill

\small 
\noindent 
\textsc{Yat-Ming Chan} \\
\textsc{Current address: Imperial College, London} \\
\textsc{E-mail}: \texttt{yatming.chan@ic.ac.uk}

\end{document}